\documentclass{article}
\usepackage[dvips]{graphicx}
\usepackage{booktabs}
\usepackage{amsmath,amssymb}
\usepackage{amscd}
\usepackage[all]{xy}

\newtheorem{Def}{Definition}[section]
\newtheorem{Thm}[Def]{Theorem}
\newtheorem{Lem}[Def]{Lemma}
\newtheorem{Pro}[Def]{Proposition}
\newtheorem{Cor}[Def]{Corollary}
\newtheorem{Ex}[Def]{Example}

\newcommand{\To}{\Rightarrow}

\newcommand{\Tot}{\Leftrightarrow}

\renewcommand{\[}{\begin{eqnarray*}}
\renewcommand{\]}{\end{eqnarray*}}

\begin{document}
\title{
Theory of Interface: Category Theory, Directed Networks and Evolution of Biological Networks
}

\author{Taichi Haruna$\ ^{\rm 1}$ \\
\footnotesize{$\ ^{\rm 1}$ Graduate School of Science, Kobe University} \\
\footnotesize{1-1 Rokkodaicho, Nada, Kobe 657-8501, JAPAN} \\
\footnotesize{E-mail: tharuna@penguin.kobe-u.ac.jp} 
}
\date{}

\maketitle

\begin{abstract}
Biological networks have two modes. The first mode is static: a network is a passage 
on which something flows. The second mode is dynamic: a network is a pattern constructed 
by gluing functions of entities constituting the network. In this paper, first we discuss that these two modes can be 
associated with the category theoretic duality (adjunction) and derive a natural network structure (a path notion) 
for each mode by appealing to the category theoretic universality. The path notion corresponding to the static mode is just the 
usual directed path. The path notion for the dynamic mode is called lateral path which is 
the alternating path considered on the set of arcs. 
Their general functionalities in a network are transport and coherence, respectively. 
Second, we introduce a betweenness centrality of arcs 
for each mode and see how the two modes are embedded in various real biological network data. 
We find that there is a trade-off relationship between the two centralities: if the value of one is large then 
the value of the other is small. This can be seen as a kind of division of labor in a network into 
transport on the network and coherence of the network. Finally, we propose an optimization model 
of networks based on a quality function involving intensities of the two modes in order to see how networks 
with the above trade-off relationship can emerge through evolution. We show that the trade-off relationship can be observed in 
the evolved networks only when the dynamic mode is dominant in the quality function by numerical simulations. 
We also show that the evolved networks have features qualitatively similar to real biological networks 
by standard complex network analysis. 
\end{abstract}

\section{Introduction}
In this decade, large interaction network data on biological, social and technological systems have 
become available. Science of complex networks has attempted to reveal structures and functions underlying 
these network data by proposing various mathematical indices and models \cite{Albert2002,Newman2003,Boccaletti2006}. 
For example, the notions of small-world property, scale-free property and modular organization have become 
inevitable tools to understand complex systems. On the other hand, 
we are aware of the criticism for purely graph-theoretic analysis forgetting meanings of networks \cite{Arita2004}. 
However, it makes mathematical analysis difficult if we stick to a meaning of each individual network too much. 
We might need a mathematical language that makes discussion on meaning of networks in a large sense possible. 
This paper is an attempt to discuss a comprehensive meaning of nodes and arcs in directed biological networks 
by appealing to category theory \cite{MacLane1998}. 

Applications of category theory to biology originates from papers by R. Rosen published in the late 1950s 
\cite{Rosen1958a,Rosen1958b}. In the beginning, Rosen proposed a model of the maintenance mechanism of metabolic networks 
in terms of category theory. However, it seems that the viewpoint of `network' became implicit as his 
theory of the metabolism-repair system developed. There are several other attempts to describe functions of 
biological systems by category theory (for example, \cite{Ehresmann1987,Wolkenhauer2007}). 
Since these studies propose highly abstract models of biological systems, it is not so clear 
how they can be applied to real world data. 
In this paper, we theoretically extend the research direction which we have sought in 
recent years \cite{Haruna2007,Haruna2008,Haruna2009,Haruna2010,Haruna2011} and try to make a bridge to 
real world data analysis. 

Our starting point is two modes of biological networks. One is static and the other is dynamic. 
In the static mode, a network is a passage on which something flows. In the dynamic mode, 
a network is a pattern constructed by gluing functions of entities constituting the network. 
For example, let us consider a neuronal network: nodes are neurons and arcs are synaptic connections 
between neurons. In the static mode, the neuronal network is a passage on which electric or chemical signals flow. 
On the other hand, in the dynamic mode, each node (a neuron) is an information processing entity that 
receives signals from other nodes, modifies them and sends its response signals to other nodes. 
Similar pictures can hold for other biological networks such as ecological flow networks and gene regulation 
networks. Note that our dynamic mode does not consider change along the time parameter directly. 
Hence, it is different from both dynamics on networks 
such as percolation \cite{Dorogovtsev2008}, synchronization \cite{Arenas2008} and games \cite{Szabo2007} and 
dynamic structural change in network structure as in temporal networks \cite{Holme2012}, although there may be 
conceptual links with them. 

Now, let us consider the following question: can we associate a natural network structure with each mode? 
For the static mode, the answer seems to be simple and intuitive. If we assume that something flows along 
the direction of arcs, then the notion of directed path may be a natural network structure corresponding to 
the static mode. On the other hand, it seems difficult to derive a natural network structure for the dynamic mode 
intuitively. In this paper, we make use of category theoretic universality to solve this problem. 

Before closing this section, we sketch the story without using category theoretic terminology. 
The directed path corresponding to the static mode is generated by the network transformation $R$
\footnote{For a directed network $G$, $R(G)$ is so-called \textbf{line graph} of $G$. We use the two terms 
$network$ and $graph$ as synonymous words. } defined below 
which sends arcs to nodes in the following sense: 
let $G=(A,N,\partial_0,\partial_1)$ be a directed network, where $A$ is the set of arcs, $N$ is the set of nodes 
and $\partial_0, \partial_1$ are maps from $A$ to $N$ sending each arc to its source or target, respectively. 
A directed network $R(G)=(A^*,N^*,\partial_0^*,\partial_1^*)$ is defined by putting 
$A^*=\{ (f,g) \in A \times A | \partial_1(f)=\partial_0(g) \}$, $N^*=A$ and 
$\partial_0^*(f,g)=f$, $\partial_1^*(f,g)=g$ for $(f,g) \in A^*$. 
The set of arcs for $R(G)$ is the set of directed paths of length 2 in $G$. 
If we apply $R$ to $G$ twice, then we obtain the set of directed paths of length 3 in $G$ as the set of arcs for $R^2(G)$. 
In general, the set of arcs for $R^n(G)$ is the set of directed paths of length $n+1$ in $G$ for any $n \geq 0$. 

For the network transformation $R$, we have the dual network transformation $L$ (in category theoretic terminology, 
both $R$ and $L$ can be extended to endofunctors on the category of directed networks and $L$ is the left adjoint functor to $R$). 
The network transformation $L$ sends each node to an arc: for any directed network $G=(A,N,\partial_0,\partial_1)$, 
$L(G)=(A_*,N_*,\partial_{0*},\partial_{1*})$ is defined by putting 
$A_*=N$, $N_*=N \times \{0,1\} / \sim$ and $\partial_{0*}(x)=[(x,0)]$, $\partial_{1*}(x)=[(x,1)]$ for $x \in A_*$. 
Here, $\sim$ is the equivalence relation on the set $N \times \{0,1\}$ generated by the relation $r$ defined by 
$(x,1) r (y,0) :\Tot$ there exists $f \in A$ such that $\partial_0(f)=x$ and $\partial_1(f)=y$. 
$[(x,i)]$ is the equivalence class containing $(x,i)$. As we will discuss in Section 2, the network transformation 
$L$ can be associated with the dynamic mode. From this, we can derive a path notion called lateral path 
that can be seen as a path notion dual to the directed path. 

Let us consider the meaning of the network transformation $L$. When we apply $L$ to a directed network $G$, each 
node is mapped to an arc. We regard this arc as representing \textit{function} of the node, namely, 
the arc in $L(G)$ to which a node in $G$ is mapped is thought of representing a process occurring on the node. 
On the other hand, we can regard each arc $f$ in $G$ as being sent to the node $[(\partial_0(f),1)](=[(\partial_1(f),0)])$ 
connecting two arcs $\partial_0(f)$ and $\partial_1(f)$ in $L(G)$. Namely, \textit{interaction is interface between 
functions}. This idea is materialized as a mathematical entity by the map 
$\varphi:A \to N_*$ defined by $f \mapsto [(\partial_0(f),1)]$. 
For arcs $f,g$ in $G$, a necessary and sufficient condition for the equality $\varphi(f)=\varphi(g)$ is 
the existence of an alternating sequence of arcs in $G$ connecting $f$ and $g$ such as 
\[
\xymatrix@-1pc{
& \bullet \ar[ld]_-f \ar[rd] & & \ar[ld] \cdots \ar[rd] & & \bullet \ar[ld] \ar[rd]^-g \\
\bullet & & \bullet & & \bullet & & \bullet. 
}
\]
Note that there are $2 \times 2=4$ possibilities for the situation at the two ends of the sequence. 
One of them is drawn above. We call such an alternating sequence of arcs \textbf{lateral path}
\footnote{Similar notion called \textbf{alternating walk} is considered in \cite{Crofts2010}. However, 
it is defined on the set of nodes. In \cite{Crofts2010}, it is used as an auxiliary means 
to obtain a bipartition of directed networks. On the other hand, lateral path in this paper 
is a central stuff associated with the dynamic mode of biological networks.}. 
Since the lateral path is associated with gluing of functions, we introduce the term \textit{coherence} 
for its general functionality. On the other hand, the functionality of the directed path is considered 
as \textit{transport} on a network here. 

We claim that \textit{the notion of lateral path is a natural network structure associated with 
the dynamic mode}. This claim is precisely formulated and proved in Section 2. In this rough sketch, 
it may be enough to give the following explanation: in the network transformation $L$, function of 
a node is represented by a single arc. However, we can represent function by a more complicated graph. 
It does not even need to be a graph. We can show that the representation of function corresponding to $L$ occupies 
a special position among \textit{all} representations of function of a node (in the language of category theory, 
it satisfies a certain universality): any representation of function gives rise to a map on the set of arcs that 
materializes the idea ``interaction as interface between functions''. This map in turn induces 
an equivalence relation on the set of arcs by its fibers. The claim is that the equivalence relation for 
the above map $\varphi$ induced by $L$ gives the finest partition of the set of arcs among all 
such equivalence relations. 

This paper is organized as follows. 
In Section 2, we describe a mathematically precise formulation of the story sketched above. 
In Section 3, we introduce betweenness centralities of arcs for the static and dynamic modes based on 
directed path and lateral path, respectively. We see how these two modes are embedded in various real biological network data. 
As a result, we find that there is a trade-off relationship between the two centralities: if the value of one is large then 
the value of the other is small. 
In Section 4, we propose an optimization model of networks based on a quality function involving intensities 
of the static and dynamic modes. By numerical simulation, we see how networks with the above trade-off relationship can emerge. 
We show that the trade-off relationship emerges after the evolution only when the dynamic mode is dominant 
in the quality function. We also show that the evolved networks have features qualitatively similar to real biological networks 
by standard complex network analysis. 
In Section 5, we give conclusions and indicate some future directions. 
In Appendix, we present a generalization of the theory for directed networks described in Section 2 to more general presheaves. 

\section{Theory of interface in directed networks}
In this section, we develop category theoretic formulation of the story sketched in Introduction. 
For necessary category theoretic background, readers are referred to \cite{MacLane1998} or \cite{Borceux1994}. 

Although this section is the main part of this paper, the readers who do not want to go into 
category theory can safely skip this section and jump to Section 3 because the main concern in this section 
is a mathematical justification of the connection between the dynamic mode of directed networks and the notion of lateral path 
and how the relationship between the static and the dynamic modes is associated with an adjunction. 
The subject of Section 3 and 4 in turn is applications of the lateral path and the directed path. 

\subsection{Left Kan extension along the Yoneda embedding functor}
In this subsection, we review the left Kan extension along the Yoneda embedding functor which is a key material for our theory. 

Let ${\mathcal C}$ be a small category, ${\mathcal D}$ a cocomplete category and $M:{\mathcal C} \to {\mathcal D}$ 
a functor. Recall that a \textbf{presheaf} on ${\mathcal C}$ is a contravariant functor from $\mathcal{C}$ to $\mathbf{Set}$. 
The category of presheaves on ${\mathcal C}$ is a functor category $\mathbf{Set}^{\mathcal{C}^{\rm op}}$. We 
denote it by $\hat{{\mathcal C}}$. The \textbf{Yoneda embedding functor} ${\mathbf y}:{\mathcal C} \to \hat{{\mathcal C}}$ is 
defined by sending each object $c \in {\mathcal C}$ to the presheaf ${\mathbf y}(c)={\mathcal C}(-,c)$ and 
each morphism $f:c \to c'$ in ${\mathcal C}$ to the natural transformation 
${\mathbf y}(f)=f \circ (-) : {\mathbf y}(c) \to {\mathbf y}(c')$. 

The \textbf{left Kan extension} of functor $M$ along the Yoneda embedding functor ${\mathbf y}$ is the functor 
${\rm Lan}_{\mathbf y}M : \hat{{\mathcal C}} \to {\mathcal D}$ defined as follows: 
for any object $G \in {\hat{\mathcal C}}$, ${\rm Lan}_{\mathbf y}M(G)$ is given by the colimit 
\[
{\rm Lan}_{\mathbf y}M(G)={\rm colim} \left( 
\xymatrix{
{\rm Elts}(G) \ar[r]^-{\pi_G} & {\mathcal C} \ar[r]^-M & {\mathcal D}
}
\right). 
\]
Here, ${\rm Elts}(G)$ is the category of elements of $G$ whose objects are pairs $(c,x)$ of object $c$ in ${\mathcal C}$ and 
$x \in G(c)$ and morphisms from $(c,x)$ to $(c',x')$ are morphisms $f:c \to c'$ in ${\mathcal C}$ such that $G(f)(x')=x$. 
$\pi_G:{\rm Elts}(G) \to {\mathcal C}$ is the projection functor to ${\mathcal C}$. 

As we will review below, the functor ${\rm Lan}_{\mathbf y}M$ is the left adjoint functor to 
the ${\rm Hom}$ functor.  As usual, we denote ${\rm Lan}_{\mathbf y}M(G)$ by $G \otimes M$. 

For any morphism $\theta:G \to H$ in $\hat{\mathcal{C}}$, morphism 
$\theta \otimes M : G \otimes M \to H \otimes M$ in $\mathcal{D}$ is defined as follows: 
let $\{\mu_{(c,x)}^{M,G} : M \circ \pi_G(c,x)=M(c) \to G \otimes M\}_{(c,x) \in {\rm Elts}(G)}$ be the colimit on $M \circ \pi_G$ 
defining $G \otimes M$. For any $x \in G(c)$ and $c \in \mathcal{C}$, we have $\theta_c(x) \in H(c)$. Hence, 
$\{\mu_{(c,\theta_c(x))}^{M,H} : M \circ \pi_H(c,\theta_c(x))=M(c) \to H \otimes M\}_{(c,x) \in {\rm Elts}(G)}$
is a cocone on $M \circ \pi_G$, where $\mu^{M,H}$ is the colimit on $M \circ \pi_H$ defining $H \otimes M$. By the universality of $G \otimes M$, 
there is a unique morphism $\xi:G \otimes M \to H \otimes M$ such that 
$\mu_{(c,\theta_c(x))}^{M,H}=\xi \circ \mu_{(c,x)}^{M,G}$ for all $(c,x) \in {\rm Elts}(G)$. 
We define $\theta \otimes M:=\xi$. 

The following diagram commutes up to isomorphism: 
\[
\xymatrix{
{\mathcal C} \ar[r]^M \ar[rd]_{\mathbf y} & {\mathcal D} \\
& {\hat{\mathcal C}}. \ar[u]_{{\rm Lan}_{\mathbf y}M}
}
\]
Namely, $\mathbf{y} \otimes M \cong M$. 

The functor $(-) \otimes M : \hat{\mathcal{C}} \to \mathcal{D}$ is the left adjoint functor to 
the functor ${\mathcal D}(M(-_2),-_1) : \mathcal{D} \to \hat{\mathcal{C}}$. Namely, we have a 
natural bijection 
\[
\mathcal{D}(G \otimes M, F) \cong \hat{\mathcal{C}}(G,\mathcal{D}(M(-),F))
\]
for any $G \in \hat{\mathcal C}$ and $F \in \mathcal{D}$. Indeed, if we define the natural transformation 
$\eta^{M,G}:G \to \mathcal{D}(M(-),G \otimes M)$ by 
$\eta_c^{M,G} : G(c) \to \mathcal{D}(M(c),G \otimes M) : x \mapsto \mu_{(c,x)}^{M,G}$ for each object 
$c \in {\mathcal C}$, then for any morphism $\theta : G \to \mathcal{D}(M(-),F)$ there exists 
a unique morphism $\hat{\theta}:G \otimes M \to F$ such that 
$\theta=\mathcal{D}(M(-),\hat{\theta}) \circ \eta^{M,G}$: 
since $\{\theta_c(x) : M \circ \pi_G(c,x)=M(c) \to F\}_{(c,x) \in {\rm Elts}(G)}$ 
is a cocone on $M \circ \pi_G$, there exists a unique morphism $\hat{\theta} : G \otimes M \to F$ such that 
$\theta_c(x)=\hat{\theta} \circ \mu_{(c,x)}^{M,G}$ for all $(c,x) \in {\rm Elts}(G)$. The condition 
that $\theta_c(x)=\hat{\theta} \circ \mu_{(c,x)}^{M,G}$ for any $(c,x) \in {\rm Elts}(G)$ 
is equivalent to the equality $\theta=\mathcal{D}(M(-),\hat{\theta}) \circ \eta^{M,G}$. 
\[
\xymatrix{
G \ar[r]^-{\eta^{M,G}} \ar[rd]_{\theta} & {\mathcal D}(M(-),G \otimes M) \ar[d]^{{\mathcal D}(M(-),\hat{\theta})} & G \otimes M \ar[d]^{\hat{\theta}} \\
& {\mathcal{D}(M(-),F)} & F 
}
\]

\subsection{Interface of a representation}
\label{subsec2-1}
Let $\upuparrows$ be the category freely generated by the following directed graph: 
\[
\xymatrix{
0 \ar@<1ex>[r]^-{m_0} \ar@<-1ex>[r]_-{m_1} & 1.
}
\]
We can identify the category of directed networks consisting of directed networks and homomorphisms between 
directed networks with $\hat{\upuparrows}$, the category of presheaves on $\upuparrows$. 
Indeed, for any object $G \in \hat{\upuparrows}$, the quadruplet $(G(1),G(0),G(m_0),G(m_1))$ defines a 
directed network. On the other hand, any directed network $(A,N,\partial_0,\partial_1)$ gives rise to 
a presheaf $G$ on $\upuparrows$ by $G(0)=N, G(1)=A,G(m_0)=\partial_0$ and $G(m_1)=\partial_1$. 
It is obvious that directed network homomorphisms between directed networks are translated into natural transformations 
between corresponding presheaves, and vice versa. 

We can regard $\upuparrows$ as a kind of data type corresponding to directed networks and 
any functor $M: \upuparrows \to \mathcal{D}$ is a representation of it. Hence, hereafter, 
for any small category $\mathcal{C}$ and category $\mathcal{D}$, we call any functor $M: \mathcal{C} \to \mathcal{D}$ 
a \textbf{representation} of ${\mathcal C}$. 
In the following, we frequently assume that $\mathcal{D}$ is bicomplete for technical reasons. 
Even under this constraint, since any presheaf category is bicomplete, we still have 
a large freedom for representations. 

Let us define the representation $M_0 : \upuparrows \to \hat{\upuparrows}$ as follows: 
\[
M_0(0)=
\xymatrix{
p_0 \ar[r]^-a & p_1
}, \ \ 
M_0(1)=
\xymatrix{
q_0 \ar[r]^-{b_0} & q_1 \ar[r]^-{b_1} & q_2, 
}
\]
$m_0 \cdot a=b_0$ and $m_1 \cdot a=b_1$, where we write $m_i \cdot a:=M_0(m_i)(1)(a)$ for $i=0,1$. 
$M_0(0)$ is the representation of a node and $M_0(1)$ is the representation of interaction between nodes (arc). 
In the terminology of Introduction, we can regard function of a node as being represented by the arc $a$ in $M_0$ and 
an arc between two nodes (interaction) in a directed network is represented by the node $q_1$ connecting the two arcs 
(representations of function of a node) $b_0$ and $b_1$ (``interaction as interface between functions''). 

In general, we regard any representation $M:\mathcal{\upuparrows} \to \mathcal{D}$ as a representation of 
function of the unit of directed network, namely an arc together with its source and target nodes. For a directed network $G$, 
$G \otimes M$ can be seen as a representation of the dynamic mode -- \textit{a network is a pattern constructed 
by gluing functions of entities constituting the network} -- because $G \otimes M$ is constructed by gluing copies of $M$. 

One can see that the network transformation $L$ introduced in Introduction is isomorphic to 
the left Kan extension of $M_0$ along the Yoneda embedding functor and the network transformation $R$ is 
isomorphic to its right adjoint. Namely, $L \cong (-) \otimes M_0$ and $R \cong \hat{\upuparrows}(M_0(-_2),-_1)$. 
Thus, we obtain the functor associated with the static mode as the right adjoint functor to the functor 
associated with the dynamic mode. 

The node $q_1$ in $M_0(1)$ is a pullback of $M_0(m_0)$ and $M_0(m_1)$. 
Motivated by this fact, we define the notion of interface as follows: 

\begin{Def}
Let $\mathcal{C}$ be a small category and $\mathcal{D}$ a complete category. 
For any representation $M:\mathcal{C} \to \mathcal{D}$ of $\mathcal{C}$, its \textbf{interface} 
is the functor ${\rm Int}_M : \mathcal{C}^{\rm op} \to \mathcal{D}$ defined as follows: 
for any object $c \in \mathcal{C}$, we define 
\[
{\rm Int}_M(c)={\rm lim} \left( 
\xymatrix{
{\rm Elts}({\mathbf y}(c)) \ar[r]^-{\pi_c} & {\mathcal C} \ar[r]^-{M} & {\mathcal D}
}
\right). 
\]
Let $f:c' \to c$ be a morphism in ${\mathcal C}$ and 
$\{\nu_{(c'',g)}^{M,c'} : {\rm Int}_M(c') \to M \circ \pi_{c'}(c'',g)=M(c'')\}_{(c'',g) \in {\rm Elts}({\mathbf y}(c'))}$ 
the limit on $M \circ \pi_{c'}$ defining ${\rm Int}_M(c')$. 
Since $\{\nu_{(c'',f \circ g)}^{M,c} : {\rm Int}_M(c) \to M \circ \pi_{c}(c'',f \circ g)=M(c'')\}_{(c'',g) \in {\rm Elts}({\mathbf y}(c'))}$ 
is a cone on $M \circ \pi_{c'}$, there exists a unique morphism $\xi : {\rm Int}_M(c) \to {\rm Int}_M(c')$ such that 
$\nu_{(c'',g)}^{M,c'} \circ \xi=\nu_{(c'',f \circ g)}^{M,c}$ for all $(c'',g) \in {\rm Elts}({\mathbf y}(c'))$. 
We define ${\rm Int}_M(f)=\xi$. 
\label{def1}
\end{Def}

\subsection{Interface transformation}
In this subsection, we introduce what we call interface transformation which generalizes the map $\varphi$ 
defined in Introduction. 

\begin{Def}
Let $\mathcal{C}$ be a small category, $\mathcal{D}$ a bicomplete category and $M:\mathcal{C} \to \mathcal{D}$ 
a representation of $\mathcal{C}$. Fix an object $G \in \hat{\mathcal{C}}$. 
For any $x \in G(c)$ and $c \in \mathcal{C}$, we put 
$\varphi_c^{M,G}(x):=\mu_{(c,x)}^{M,G} \circ \nu_{(c,{\rm id}_c)}^{M,c}$. 
As we will see below in Proposition \ref{pro1}, $\varphi^{M,G}$ is a natural transformation 
from the contravariant functor $G$ on $\mathcal{C}$ to the covariant functor $\mathcal{D}({\rm Int}_M(-),G \otimes M)$
\footnote{
Natural transformations from a contravariant functor to a covariant functor can be seen as dinatural transformations \cite{MacLane1998}. 
Let $\widetilde{G}$ be the composition 
$\xymatrix{\mathcal{C}^{\rm op} \times \mathcal{C} \ar[r]^-{\pi_1} & \mathcal{C}^{\rm op} \ar[r]^-{G} & \mathbf{Set}}$ and 
$\widetilde{H}$ the composition 
$\xymatrix{\mathcal{C}^{\rm op} \times \mathcal{C} \ar[r]^-{\pi_2} & \mathcal{C} \ar[r]^-{H} & \mathbf{Set}}$, 
where $H=\mathcal{D}({\rm Int}_M(-),G \otimes M)$. 
Then, $\varphi^{M,G}$ is a dinatural transformation from $\widetilde{G}$ to $\widetilde{H}$. 
}. 
We call $\varphi^{M,G}$ \textbf{interface transformation} of $G$ with respect to $M$. 
\[
\xymatrix{
{\rm Int}_M(c) \ar@/_1pc/[rr]_{\varphi_c^{M,G}(x)} \ar[r]^-{\nu_{(c,{\rm id}_c)}^{M,c}} & M(c) \ar[r]^-{\mu_{(c,x)}^{M,G}} & G \otimes M. 
}
\]
\label{def2}
\end{Def}

Let us check that interface transformations defined in Definition \ref{def2} are actually 
natural transformations from contravariant functors to covariant functors. 
Hereafter, we write $x \cdot f:=G(f)(x)$ for any $x \in G(c)$ and $f:c' \to c$ in $\mathcal{C}$. 

\begin{Pro}
Let $\mathcal{C}$ be a small category and $\mathcal{D}$ a bicomplete category. 
For any representation $M:\mathcal{C} \to \mathcal{D}$ of $\mathcal{C}$ and presheaf $G \in \hat{\mathcal{C}}$, 
any morphism $f:c' \to c$ in ${\mathcal C}$ makes the following diagram commute: 
\[
\xymatrix{
G(c) \ar[d]_-{G(f)} \ar[r]^-{\varphi_c^{M,G}} & {\mathcal D}({\rm Int}_M(c),G \otimes M) \\
G(c') \ar[r]^-{\varphi_{c'}^{M,G}} & {\mathcal D}({\rm Int}_M(c'),G \otimes M). \ar[u]_-{(-) \circ {\rm Int}_M(f)}
}
\]
\label{pro1}
\end{Pro}
\textit{Proof. }
Let us show $\varphi_c^{M,G}(x)=\varphi_{c'}^{M,G}(x \cdot f) \circ {\rm Int}_M(f)$ for any $x \in G(c)$. 
However, the equality follows from a diagram chasing over the following commutative diagram: 
\[
\xymatrix{
{\rm Int}_M(c) \ar[r]^-{\nu_{(c,{\rm id}_c)}^{M,c}} \ar[d]_-{{\rm Int}_M(f)} \ar[rd]^-{\nu_{(c',f)}^{M,c}} & M(c) \ar[rd]^-{\mu_{(c,x)}^{M,G}} \\
{\rm Int}_M(c') \ar[r]_-{\nu_{(c',{\rm id}_{c'})}^{M,c'}} & M(c') \ar[u]_-{M(f)} \ar[r]_-{\mu_{(c',x \cdot f)}^{M,G}} & G \otimes M. 
}
\]

\hfill $\Box$ \\

\subsection{Category $\mathcal{IT}_G$}
Now, we define a category where interface transformations live. 

\begin{Def}
Let $\mathcal{C}$ be a small category. For any $G \in \hat{\mathcal{C}}$, we define a category 
$\mathcal{IT}_G$ as follows: 
objects in $\mathcal{IT}_G$ are pairs $(H,\varphi)$ consisting of 
a functor $H : \mathcal{C} \to \mathbf{Set}$ and a natural transformation $\varphi:G \to H$. 
Here, the natural transformation $\varphi$ from the contravariant functor $G$ to the covariant functor $H$ is 
defined as a family of maps $\{ \varphi_c : G(c) \to H(c)\}_{c \in C}$ that makes the following diagram 
commute for any morphism $f:c' \to c \in C$: 
\[
\xymatrix{
G(c) \ar[d]_-{G(f)} \ar[r]^-{\varphi_c} & H(c) \\
G(c') \ar[r]^-{\varphi_{c'}} & H(c'). \ar[u]_-{H(f)}
}
\]
Morphisms from $(H,\varphi)$ to $(H',\varphi')$ are natural transformations $n:H \to H'$ such that 
$\varphi'=n \circ \varphi$:
\[
\xymatrix@-1pc{
H(c) \ar[rr]^-{n_c} && H'(c) \\
& G(c) \ar[ul]^-{\varphi_c} \ar[ur]_-{\varphi'_c}.
}
\]
\label{def3}
\end{Def}

From Proposition \ref{pro1}, we have the following fact:

\begin{Pro}
Let $\mathcal{C}$ be a small category and $\mathcal{D}$ a bicomplete category. For any 
representation $M:\mathcal{C} \to \mathcal{D}$ of $\mathcal{C}$ and presheaf $G \in \hat{\mathcal{C}}$, 
the interface transformation $\varphi^{M,G}$ of $G$ with respect to $M$ gives rise to the object 
$({\mathcal D}({\rm Int}_M(-),G \otimes M), \varphi^{M,G})$ in $\mathcal{IT}_G$. 
\label{pro2}
\end{Pro}

\subsection{Initiality of $\mathcal{IT}_G$}

\begin{Thm}
For any directed network $G \in \hat{\upuparrows}$, $\mathcal{IT}_G$ has an initial object. 
In particular, the object $(\hat{\upuparrows}({\rm Int}_{M_0}(-),G \otimes M_0), \varphi^{M_0,G})$ 
induced by the interface transformation $\varphi^{M_0,G}$ of $G$ with respect to the representation 
$M_0 : \upuparrows \to \hat{\upuparrows}$ given in subsection \ref{subsec2-1} is an initial object 
in $\mathcal{IT}_G$. 
\label{thm1}
\end{Thm}
\textit{Proof. }
We show that for any object $(H,\psi)$ in $\mathcal{IT}_G$, there exists a unique morphism 
$\iota:\hat{\upuparrows}({\rm Int}_{M_0}(-),G \otimes M_0) \to H$ such that $\psi=\iota \circ \varphi^{M_0,G}$. 

By ${\rm Int}_{M_0}(0)=M_0(0) \cong \mathbf{y}(1)$, we have 
$\hat{\upuparrows}({\rm Int}_{M_0}(0),G \otimes M_0) \cong G \otimes M_0(1) \cong G(0)$. 
If we identify the leftmost-hand side with the rightmost-hand side of this sequence of isomorphisms, 
then one can see that $\varphi_0^{M_0,G}$ becomes the identity ${\rm id}_{G(0)}$ by chasing how elements move by 
each isomorphism. Hence, we should take $\iota_0$ as $\iota_0:=\psi_0$. 

Let us construct $\iota_1$. 
In general, for any representation $M:\mathcal{C} \to \hat{\mathcal{C}}$ of a small category ${\mathcal C}$, 
we can take 
\[
G \otimes M(c)=\left( \sum_{c' \in \mathcal{C}} G(c') \times M(c')(c) \right)/\sim 
\]
for any $G \in \hat{\mathcal{C}}$ and $c \in \mathcal{C}$. 
Here, $\sim$ is the equivalence relation generated by the following binary relation $r$: 
$(x,f \cdot p) r (x \cdot f,p)$ for $x \in G(c''), \ p \in M(c')(c)$ and 
$f : c' \to c'' \in \mathcal{C}$. 
Components of the colimit $\mu^{M,G}$ are given by 
$\mu_{(c,x)}^{M,G}(c'):M(c)(c') \to G \otimes M(c'): p \mapsto [(x,p)]$, 
where $[(x,p)]$ is the equivalence class containing $(x,p)$. 

Now, we have ${\rm Int}_{M_0}(1)=\{(p_1,p_0)\} \cong {\mathbf y}(0)$ and 
$\hat{\upuparrows}({\rm Int}_{M_0}(1),G \otimes M_0) \cong G \otimes M_0(0)$. 
We also have $\nu_{(1,{\rm id}_1)}^{M_0,1}(0)(p_1,p_0)=q_1 \in M_0(1)(0)$, 
$\nu_{(0,m_i)}^{M_0,1}(0)(p_1,p_0)=p_{i+1} \in M_0(0)(0)$ for $i=0,1$. 
Since ${\rm Int}_{M_0}(0)=M_0(0)$, ${\rm Int}_{M_0}(m_i)=\nu_{(0,m_i)}^{M_0,1}$ for $i=0,1$. 
Here, $i+1$ is considered modulo 2. 
Hence, we have $\varphi_1^{M_0,G}(x)=[(x,q_1)]$ for any $x \in G(1)$. 
$(-) \circ {\rm Int}_{M_0}(m_i) : \hat{\upuparrows}({\rm Int}_{M_0}(0),G \otimes M_0)=G(0) \to \hat{\upuparrows}({\rm Int}_{M_0}(1),G \otimes M_0)=G \otimes M_0(0)$
is given by $x \mapsto [(x,p_{i+1})]$. 

For $i=0,1$, $\iota_1$ should make the following diagram commute:
\[
\xymatrix{
G \otimes M_0(0) \ar[r]^-{\iota_1} & H(1) \\
G(0) \ar[u]^-{[(-,p_{i+1})]} \ar[r]^-{\iota_0} & H(0). \ar[u]_-{H(m_i)}
}
\]
Namely, for any $x \in G(0)$, we should have 
$\iota_1([(x,p_{i+1})])=H(m_i) \circ \iota_0(x), \ i=0,1$. 
Let us see that $\iota_1$ is uniquely determined by this condition. 

First, any element of $G \otimes M_0(0)=\left[ G(0) \times \{p_0,p_1\} + G(1) \times\{q_0,q_1,q_2\} \right]/\sim$ 
has a representative of the form $(x,p), \ x \in G(0), \ p \in \{p_0,p_1\}$. 
Indeed, for any $y \in G(1)$, we have 
$(y,q_0)=(y,m_0 \cdot p_0) \sim (y \cdot m_0,p_0)$, 
$(y,q_1)=(y,m_0 \cdot p_1) \sim (y \cdot m_0,p_1)$ and $(y,q_2)=(y,m_1 \cdot p_1) \sim (y \cdot m_1,p_1)$. 

Second, $\iota_1$ is well-defined. 
It is enough to show that $H(m_0)\circ\iota_0(x)=H(m_1)\circ\iota_0(y)$ for any $x, y \in G(0)$ 
such that $x=z \cdot m_0, \ y=z \cdot m_1$ for some $z \in G(1)$. 
However, since $\psi$ is a natural transformation from $G$ to $H$, 
we obtain 
$H(m_0)\circ\iota_0(x)=H(m_0) \circ \psi_0(x)=H(m_0) \circ \psi_0 (z \cdot m_0)=\psi_1(z)
=H(m_1) \circ \psi_0 (z \cdot m_1)=H(m_1) \circ \psi_0(y)=H(m_1)\circ\iota_0(y)$. 

Finally, $\psi_1 =\iota_1 \circ \varphi_1^{M_0,G}$ holds because 
$\iota_1 \circ \varphi_1^{M_0,G}(x)=\iota_1([(x,q_1)])=\iota_1([(x \cdot m_0,p_1)])
=H(m_0) \circ \iota_0(x \cdot m_0)=H(m_0) \circ \psi_0(x \cdot m_0)=\psi_1(x)$ 
for any $x \in G(1)$. 

\hfill $\Box$\\

\begin{Cor}
Let $G$ be a directed network. 
For any representation $M:\mathcal{\upuparrows} \to \mathcal{D}$ of $\upuparrows$ 
in any bicomplete category $\mathcal{D}$, 
let $\sim_M$ be the equivalence relation on $G(1)$ induced by the fibers of the map 
$\varphi_1^{M,G}:G(1) \to {\mathcal D}({\rm Int}_{M}(1),G \otimes M)$, 
the 1-component of the interface transformation of $G$ with respect to $M$. 
Namely, we define $x \sim_{M} y : \Tot \varphi_1^{M,G}(x)=\varphi_1^{M,G}(y)$. 
Then, $\sim_{M_0} \subseteq \sim_{M}$ for any $M$. 
\label{cor1}
\end{Cor}

The equivalence relation $\sim_{M_0}$ on $G(1)$ in Corollary \ref{cor1} can be given explicitly as follows: 
for $x,y \in G(1)$, since $\varphi_1^{M_0,G}(x)=\varphi_1^{M_0,G}(y) \Tot [(x,q_1)]=[(y,q_1)]$, 
we have $x \sim_{M_0} y$ if and only if there exists an alternating sequence of arcs between 
$x$ and $y$ in $G$ such as 
\[
\xymatrix@-1pc{
& \bullet \ar[ld]_-x \ar[rd] & & \ar[ld] \cdots \ar[rd] & & \bullet \ar[ld] \ar[rd]^-y \\
\bullet & & \bullet & & \bullet & & \bullet. 
}
\]
Note that there are four cases for the situation at both ends of the sequence. 
For $x,y \in G(1)$, we say that $x$ and $y$ are \textbf{laterally connected} if $\varphi_1^{M_0,G}(x)=\varphi_1^{M_0,G}(y)$. 
The alternating sequence of arcs connecting $x$ and $y$ is called \textbf{lateral path}. 

\subsection{Examples}
\begin{Ex}\normalfont
Let us consider the representation $\mathbf{y}: \upuparrows \to \hat{\upuparrows}$. 
For any directed network $G \in \hat{\upuparrows}$, $\varphi_0^{\mathbf{y},G}$ is 
an isomorphism by ${\rm Int}_{\mathbf{y}}(0)=\mathbf{y}(0)$. 
$\varphi_1^{\mathbf{y},G}$ is a map to a singleton set because ${\rm Int}_{\mathbf{y}}(1)=\emptyset$. 
Consequently, the equivalence relation $\sim_{\mathbf{y}}$ on $G(1)$ induced by the map $\varphi_1^{\mathbf{y},G}$ 
is $G(1) \times G(1)$. Namely, any $x,y \in G(1)$ are equivalent. 
\label{ex1}
\end{Ex}

\begin{Ex}\normalfont
Let us define the representation $M_u : \upuparrows \to \hat{\upuparrows}$ as follows: 
\begin{align*}
M_u(0) &=
\xymatrix{
p_0 \ar[r]^-a & p_1
}, \\
M_u(1) &=
\xymatrix{
q_0 \ar[r]^-{b_0} & q_1 & q_2 \ar[l]_-{b_1}, 
}
\end{align*}
$m_0 \cdot a=b_0$ and $m_1 \cdot a=b_1$. We have 
${\rm Int}_{M_u}(0)=M_u(0)$ and ${\rm Int}_{M_u}(1)=\{(p_1,p_1)\}$. 

For any directed network $G \in \hat{\upuparrows}$, $\varphi_0^{M_u,G}$ is an isomorphism. 
The equivalence relation $\sim_{M_u}$ on $G(1)$ defined by the fibers of $\varphi_1^{M_u,G}$ 
is generated by the equations 
$x \cdot m_0=y \cdot m_0, \ x \cdot m_0=y \cdot m_1, \ x \cdot m_1=y \cdot m_0$ and $x \cdot m_1=y \cdot m_1$ 
because $\varphi_1^{M_u,G}(x)=\varphi_1^{M_u,G}(y) \Tot [(x,q_1)]=[(y,q_1)]$ for any $x,y \in G(1)$ and 
$q_1=m_0 \cdot p_1=m_1 \cdot p_1$. 
Thus, for any $x,y \in G(1)$, $x \sim_{M_u} y$ if and only if there exists a sequence of arcs connecting 
$x$ and $y$ in $G$ without regard to the direction of arcs. We call this kind of sequence of arcs 
\textbf{undirected path}. 
\label{ex2}
\end{Ex}

\begin{Ex}\normalfont
In general, we can calculate the interface transformation $\varphi^{M,G}$ of a directed network $G$ 
with respect to a representation $M:\upuparrows \to \hat{\upuparrows}$ as follows: 
for any $x \in G(0)$, the 0-component 
\[
\varphi_0^{M,G}(x) : \xymatrix{ {\rm Int}_M(0)=M(0) \ar[rr]^-{{\nu_{(0,{\rm id}_0)}^{M,0}}={\rm id}_{M(0)}} && M(0) \ar[r]^-{\mu_{(0,x)}^{M,G}} & G \otimes M } 
\]
is given by $\varphi_0^{M,G}(x)(0)(p)=[(x,p)]$ for any $p \in M(0)(0)$ and 
$\varphi_0^{M,G}(x)(1)(a)=[(x,a)]$ for any $a \in M(0)(1)$. 

Next we consider the 1-component $\varphi_1^{M,G}$. 
We have 
${\rm Int}_M(1)(0)=\{ (p,p') \in M(0)(0) \times M(0)(0) | m_0 \cdot p=m_1 \cdot p' \}$ and 
${\rm Int}_M(1)(1)=\{ (a,a') \in M(0)(1) \times M(0)(1) | m_0 \cdot a=m_1 \cdot a' \}$. 
In the commutative diagram 
\[
\xymatrix{
{\rm Int}_M(1) \ar[r]^-{\nu_{(0,m_1)}^{M,1}} \ar[d]_-{\nu_{(0,m_0)}^{M,1}} \ar[rd]^-{\nu_{(1,{\rm id}_1)}^{M,1}} & M(0) \ar[d]^-{M(m_1)} \\
M(0) \ar[r]_-{M(m_0)} & M(1), 
}
\]
morphisms $\nu_{(0,m_0)}^{M,1}$ and $\nu_{(0,m_1)}^{M,1}$ consists of projection maps to the first component and 
the second component, respectively. 
Consequently, for any $x \in G(1)$, 
\[
\varphi_1^{M,G}(x) : \xymatrix{ {\rm Int}_M(1) \ar[r]^-{\nu_{(1,{\rm id}_1)}^{M,1}} & M(1) \ar[r]^-{\mu_{(1,x)}^{M,G}} & G \otimes M } 
\]
is given by 
$\varphi_1^{M,G}(x)(0)(p,p')=[(x,q)]$ for any $(p,p') \in {\rm Int}_M(1)(0)$, where $q:=m_0 \cdot p=m_1 \cdot  p'$, 
and $\varphi_1^{M,G}(x)(1)(a,a')=[(x,b)]$ for any $(a,a') \in {\rm Int}_M(1)(1)$, where $b:=m_0 \cdot a=m_1 \cdot  a'$. 
\label{ex3}
\end{Ex}

\begin{Ex}\normalfont
We define the representation $M_1 : \upuparrows \to \hat{\upuparrows}$ by 
\begin{align*}
M_1(0) &=
\xymatrix{
p_0 \ar[r]^-{a_0} & p_1 \ar[r]^-{a_1} & p_2
}, \\ 
M_1(1) &=
\xymatrix{
q_0 \ar[r]^-{b_0} & q_1 \ar[r]^-{b_1} & q_2 \ar[r]^-{b_2} & q_3, 
}
\end{align*}
$m_0 \cdot a_0=b_0, \ m_0 \cdot a_1=b_1, \ m_1 \cdot a_0=b_1$ and $m_1 \cdot a_1=b_2$. 
We have ${\rm Int}_{M_1}(1)(0)=\{(p_1,p_0), (p_2,p_1)\}$ and 
${\rm Int}_{M_1}(1)(1)=\{(a_1,a_0)\}$. 

Let $G$ be a directed network. 
Since $m_0 \cdot a_1=b_1$ and $m_1 \cdot a_0=b_1$, 
$\varphi_0^{M_1,G}(x)=\varphi_0^{M_1,G}(y) \Tot 
\textrm{ there exist } f_1,g_1,\cdots,f_k,g_k \in G(1) \ (k \geq 0) 
\textrm{ such that } x=f_1 \cdot m_1, f_i \cdot m_0=g_i \cdot m_0, g_i \cdot m_1=f_{i+1} \cdot m_1 \ (i=1,2,\cdots,k-1), y=g_k \cdot m_1
\textrm{ and }
\textrm{ there exist } d_1,e_1,\cdots,d_l,e_l \in G(1) \ (l \geq 0) 
\textrm{ such that } x=d_1 \cdot m_0, d_i \cdot m_1=e_i \cdot m_1, e_i \cdot m_0=d_{i+1} \cdot m_0 \ (i=1,2,\cdots,l-1), y=e_l \cdot m_0
$ for any $x,y \in G(0)$. 

\[
\xymatrix@-1pc{
& \bullet \ar[ld]_-{f_1} \ar[rd]_-{g_1} & & \ar[ld]_-{f_2} \cdots \ar[rd]^-{g_{k-1}} & & \bullet \ar[ld]^-{f_k} \ar[rd]^-{g_k} \\
x \ar[rd]_-{d_1} & & \bullet & & \bullet & & y \ar[ld]^-{e_l} \\
& \bullet & & \cdots & & \bullet & \\
& & \bullet \ar[lu]^-{e_1} \ar[ru]^-{d_2} & & \bullet \ar[lu]_-{e_{l-1}} \ar[ru]_-{d_l} & &
}
\]

For any $x,y \in G(1)$, we have $\varphi_1^{M_1,G}(x)=\varphi_1^{M_1,G}(y) \Tot [(x,b_1)]=[(y,b_1)]$. 
One can see that $[(x,b_1)]=[(y,b_1)]$ is equivalent to the existence of a lateral path 
between $x$ and $y$. 
\label{ex4}
\end{Ex}

\begin{Ex}\normalfont
The representation $M_2 : \upuparrows \to \hat{\upuparrows}$ is defined by 
\[
M_2(0) &=
\xymatrix{
p_0 \ar[r]^-{a_0} & p_1 \ar[r]^-{a_1} & p_2 \ar[r]^-{a_2} & p_3
}, \\ 
M_2(1) &=
\xymatrix{
q_0 \ar[r]^-{b_0} & q_1 \ar[r]^-{b_1} & q_2 \ar[r]^-{b_2} & q_3 \ar[r]^-{b_3} & q_4 
}
\]
and $m_0 \cdot a_i=b_i$ and $m_1 \cdot a_i=b_{i+1}$ for $i=0,1,2$. 

Let $G$ be a directed network. For any $x,y \in G(1)$, one can see that 
if $x$ and $y$ are connected by a lateral path in $G$, then 
$\varphi_1^{M_2,G}(x)=\varphi_1^{M_2,G}(y)$, and 
if $\varphi_1^{M_2,G}(x)=\varphi_1^{M_2,G}(y)$, then there exists an undirected path 
between $x$ and $y$ in $G$. However, the converses do not necessarily hold for both implications. 
For the first one, if $G$ is given by the following directed network, then 
we have $\varphi_1^{M_2,G}(x)=\varphi_1^{M_2,G}(y)$: 
\[
\xymatrix@-1pc{
& \bullet \ar[r]^-{x} & \bullet \ar[rd] & \\
\bullet \ar[ru] \ar[rd] & & & \bullet \\
& \bullet \ar[r]_-{y} & \bullet. \ar[ru] &
}
\]
For the second one, if $G$ is the following directed network, then we have 
$\varphi_1^{M_2,G}(x) \neq \varphi_1^{M_2,G}(y)$: 
\[
\xymatrix{
\bullet \ar[r]^-{x} & \bullet \ar[r]^-{y} & \bullet.
}
\]
\label{ex5}
\end{Ex}

\begin{Ex}\normalfont
In general, we define the representation 
$M_n : \upuparrows \to \hat{\upuparrows}$ by 
\begin{align*}
M_n(0) &=
\xymatrix{
p_0 \ar[r]^-{a_0} & p_1 \ar[r]^-{a_1} & \cdots \ar[r]^-{a_n} & p_{n+1}
}, \\ 
M_n(1) &=
\xymatrix{
q_0 \ar[r]^-{b_0} & q_1 \ar[r]^-{b_1} & \cdots \ar[r]^-{b_{n+1}} & q_{n+2} 
}
\end{align*}
and $m_0 \cdot a_i=b_i$ and $m_1 \cdot a_i=b_{i+1}$ for $i=0,1,\cdots,n$. 

Let $G_{n}$ be the following directed network: 
\[
\xymatrix@-1pc{
& \bullet \ar[r]^-{x_1} & \cdots \ar[r]^-{x_{n-1}} & \bullet \ar[rd]^-{x_n} & \\
\bullet \ar[ru]^-{x_0} \ar[rd]_-{y_0} & & & & \bullet \\
& \bullet \ar[r]^-{y_1} & \cdots \ar[r]^-{y_{n-1}} & \bullet. \ar[ru]_-{y_n} &
}
\]
For $n \geq 1$, we have $\varphi_1^{M_{2n},G_{2n}}(x_n)=\varphi_1^{M_{2n},G_{2n}}(y_n)$ and 
$\varphi_1^{M_{2(n-1)},G_{2n}}(x_n) \neq \varphi_1^{M_{2(n-1)},G_{2n}}(y_n)$. 

$M_n$ can be seen as a $n+1$-fold tensor product of $M_0$. 
In the next subsection, we discuss tensor product of representations. 
\label{ex6}
\end{Ex}
\subsection{Tensor product of representations and morphisms in $\mathcal{IT}_G$}
Let $\mathcal{C}$ be a small category and $\mathcal{D}$ a cocomplete category. 
The tensor product of representations $M:\mathcal{C} \to \hat{\mathcal{C}}$ and 
$N:\mathcal{C} \to \mathcal{D}$ is defined by $M \otimes N:=\left( {\rm Lan}_{\mathbf{y}} N \right) \circ M$. 
\[
\xymatrix{
{\mathcal C} \ar[r]^-{N} \ar[rd]_-{\mathbf y} & {\mathcal D} & & \\
& {\hat{\mathcal C}} \ar[u]_{{\rm Lan}_{\mathbf y}N} & & \hat{\mathcal{C}} \ar[ll]^-{{\rm Lan}_{\mathbf y} M} \ar[llu]_-{\textrm{\hspace{5mm}}{\rm Lan}_{\mathbf{y}}\left( ({\rm Lan}_{\mathbf{y}} N) \circ M \right)} \\
& \mathcal{C} \ar[u]^-{M} \ar[rru]_-{\mathbf{y}} & &
}
\]
The above diagram commutes up to isomorphism. 
Indeed, since left adjoints preserve colimits, 
we have a canonical isomorphism 
\[
\alpha_G : \xymatrix{ {\rm Lan}_{\mathbf{y}} N \left( {\rm Lan}_{\mathbf{y}}M (G) \right) \ar[r]^-{\cong} & {\rm Lan}_{\mathbf{y}}\left( ({\rm Lan}_{\mathbf{y}} N) \circ M \right)(G) }
\]
for any presheaf $G \in \hat{\mathcal{C}}$. 
If we use the symbol $\otimes$, then this becomes 
\[
\alpha_G : \xymatrix{ (G \otimes M) \otimes N \ar[r]^-{\cong} & G \otimes (M \otimes N) }. 
\]
Note that $\alpha_G$ is a unique morphism such that $\alpha_G^{-1} \circ \mu_{(c,x)}^{M \otimes N,G}=\mu_{(c,x)}^{M,G} \otimes N$ 
for all $(c,x) \in {\rm Elts}(G)$. 

\begin{Ex}\normalfont
For any representation $M : \mathcal{C} \to \hat{\mathcal{C}}$ of a small category $\mathcal{C}$, we define 
$M^{\otimes n}:=\underbrace{M \otimes \cdots \otimes M}_{n \textrm{ times}}$ up to canonical isomorphisms. 
We have $M_0^{\otimes n+1} \cong M_n$, where $M_n$ is the representation introduced in Example \ref{ex6}.
\label{ex7}
\end{Ex}

Let ${\mathcal C}$ be a small category, $\mathcal{D}$ a bicomplete category and 
$G$ a presheaf on ${\mathcal C}$. 
From Proposition \ref{pro2}, we can obtain an object of $\mathcal{IT}_G$ from any 
representation $M:\mathcal{C} \to \mathcal{D}$. 
Now, we describe a way to obtain a morphism of $\mathcal{IT}_G$. 

Let 
$\{\nu_{(c',f)}^{M,c} : {\rm Int}_{M}(c) \to M \circ \pi_{c}(c',f)=M(c')\}_{(c',f) \in {\rm Elts}({\mathbf y}(c))}$ 
be the limit on $M \circ \pi_c$ defining ${\rm Int}_{M}(c)$ and 
$\{\nu_{(c',f)}^{M \otimes N,c} : {\rm Int}_{M \otimes N}(c) \to (M \otimes N) \circ \pi_{c}(c',f)=M(c') \otimes N\}_{(c',f) \in {\rm Elts}({\mathbf y}(c))}$ 
the limit on $(M \otimes N) \circ \pi_c$ defining ${\rm Int}_{M \otimes N}(c)$. 
Since 
$\{\nu_{(c',f)}^{M,c}\otimes N : {\rm Int}_{M}(c) \otimes N \to M(c') \otimes N\}_{(c',f) \in {\rm Elts}({\mathbf y}(c))}$ 
is a cone on $(M \otimes N) \circ \pi_c$, there exists a unique morphism $\beta_c:{\rm Int}_{M}(c) \otimes N \to {\rm Int}_{M \otimes N}(c)$ 
such that $\nu_{(c',f)}^{M,c}\otimes N=\nu_{(c',f)}^{M \otimes N,c} \circ \beta_c$ for all $(c',f) \in {\rm Elts}({\mathbf y}(c))$ 
by the universality of ${\rm Int}_{M \otimes N}(c)$. 
$\beta_c$ is natural with respect to $c$. Indeed, we have the following commutative diagram 
for any morphisms $\xymatrix{c'' \ar[r]^-{f'} & c' \ar[r]^-{f} & c}$ in $\mathcal{C}$: 
\[
\xymatrix{
{\rm Int}_M(c) \otimes N \ar[rr]^-{\beta_c} \ar[dd]_-{{\rm Int}_M(f) \otimes N} \ar[rd]^-{\nu_{(c'',f \circ f')}^{M,c} \otimes N} && {\rm Int}_{M \otimes N}(c) \ar[ld]_-{\nu_{(c'',f \circ f')}^{M \otimes N,c}} \ar[dd]^-{{\rm Int}_{M \otimes N}(f)} \\
& M \otimes N(c'') & \\
{\rm Int}_M(c') \otimes N \ar[ru]^-{\nu_{(c'',f')}^{M,c'} \otimes N} \ar[rr]_-{\beta_{c'}} && {\rm Int}_{M \otimes N}(c'). \ar[lu]_-{\nu_{(c'',f')}^{M \otimes N,c'}} 
}
\]

\begin{Lem}
Let $\mathcal{A}$ and $\mathcal{B}$ be categories and $F:\mathcal{A} \to \mathcal{B}$ a functor. 
Assume that $F$ has a limit $\{\nu_A:L \to F(A)\}_{A \in \mathcal{A}}$. 
Given two morphisms $u,v:B \to L$ in $\mathcal{B}$, 
if $\nu_A \circ u=\nu_B \circ v$ for all $A \in \mathcal{A}$, then $u=v$. 
\label{lem1}
\end{Lem}
\textit{Proof. }
Put $\delta_A:=\nu_A \circ u=\nu_A \circ v$. 
Since $\{\delta_A : B \to F(A)\}_{A \in \mathcal{A}}$ is a cone on $F$, 
there exists a unique morphism $\eta$ in $\mathcal{B}$ such that $\nu_A \circ \eta=\delta_A$ for all $A \in \mathcal{A}$. 
However, both $u$ and $v$ satisfy this condition. Consequently, $u=\eta=v$. 

\hfill $\Box$ \\

By applying Lemma \ref{lem1} with $u={\rm Int}_{M \otimes N}(f) \circ \beta_c$ and $v=\beta_{c'} \circ {\rm Int}_M(f) \otimes N$, 
we get the equality ${\rm Int}_{M \otimes N}(f) \circ \beta_c=\beta_{c'} \circ {\rm Int}_M(f) \otimes N$. 
Consequently, we obtain a natural transformation $\beta:{\rm Int}_{M}(-) \otimes N \to {\rm Int}_{M \otimes N}$. 

\begin{Pro}
Let $\mathcal{C}$ be a small category and $\mathcal{D}$ a bicomplete category. 
For representations $M:\mathcal{C} \to \hat{\mathcal{C}}$ and $N:\mathcal{C} \to \mathcal{D}$, 
assume that the natural transformation $\beta:{\rm Int}_{M}(-) \otimes N \to {\rm Int}_{M \otimes N}$ is 
an isomorphism. 
Then, for any presheaf $G \in \hat{\mathcal{C}}$, a morphism 
$j: (\hat{{\mathcal C}}({\rm Int}_M(-),G \otimes M), \varphi^{M,G}) \to ({\mathcal D}({\rm Int}_{M \otimes N}(-),G \otimes (M \otimes N)), \varphi^{M \otimes N,G})$ 
in $\mathcal{IT}_G$ is obtained by defining 
$j_c(\theta)$ as the composition 
\[
\xymatrix{ {\rm Int}_{M \otimes N}(c) \ar[r]^-{\beta_c^{-1}} & {\rm Int}_M(c) \otimes N \ar[r]^-{\theta \otimes N} & (G \otimes M) \otimes N \ar[r]^-{\alpha_G} & G \otimes (M \otimes N)}
\]
for any $c \in \mathcal{C}$ and morphism $\theta : {\rm Int}_M(c) \to G \otimes M$. 
\label{pro3}
\end{Pro}
\textit{Proof. }
First, let us check that $j$ is a natural transformation. We have to show that the following diagram commutes 
for any morphism $f:c \to c'$ in $\mathcal{C}$: 
\[
\xymatrix{
\hat{{\mathcal C}}({\rm Int}_M(c),G \otimes M) \ar[r]^-{j_c} \ar[d]_-{(-) \circ {\rm Int}_M(f)} & {\mathcal D}({\rm Int}_{M \otimes N}(c),G \otimes (M \otimes N)) \ar[d]^-{(-)\circ {\rm Int}_{M \otimes N}(f)} \\
\hat{{\mathcal C}}({\rm Int}_M(c'),G \otimes M) \ar[r]^-{j_{c'}} & {\mathcal D}({\rm Int}_{M \otimes N}(c'),G \otimes (M \otimes N)). 
}
\]
However, for any morphism $\theta:{\rm Int}_M(c) \to G \otimes M$, we have 
\begin{align*}
j_{c'}(\theta \circ {\rm Int}_M(f)) &= \alpha_G \circ \left( \left( \theta \circ {\rm Int}_M(f) \right) \otimes N \right) \circ \beta_{c'}^{-1} \\
&= \alpha_G \circ \left( \theta \otimes N \right) \circ \left({\rm Int}_M(f) \otimes N \right) \circ \beta_{c'}^{-1} \\
&= \alpha_G \circ \left( \theta \otimes N \right) \circ \beta_{c}^{-1} \circ {\rm Int}_{M \otimes N}(f) \\
&= j_c(\theta) \circ {\rm Int}_{M \otimes N}(f), 
\end{align*}
where the second equality follows from the functoriality of $(-) \otimes N$ and the third equality follows from the 
naturality of $\beta$. 

Next, we show that the following diagram commutes for any $c \in \mathcal{C}$: 
\[
\xymatrix@-1pc{
\hat{{\mathcal C}}({\rm Int}_M(c),G \otimes M) \ar[rr]^-{j_c} && {\mathcal D}({\rm Int}_{M \otimes N}(c),G \otimes (M \otimes N)) \\
& G(c). \ar[lu]^-{\varphi_c^{M,G}} \ar[ru]_-{\varphi_c^{M \otimes N,G}} &
}
\]
However, for any $c \in \mathcal{C}$ and $x \in G(c)$, 
we have 
\begin{align*}
j_c \left( \varphi_c^{M,G}(x) \right) &= j_c \left( \mu_{(c,x)}^{M,G} \circ \nu_{(c,{\rm id}_c)}^{M,c} \right) \\
&= \alpha_G \circ \left( \left( \mu_{(c,x)}^{M,G} \circ \nu_{(c,{\rm id}_c)}^{M,c} \right) \otimes N \right) \circ \beta_c^{-1} \\
&= \alpha_G \circ \left( \mu_{(c,x)}^{M,G} \otimes N \right) \circ \left( \nu_{(c,{\rm id}_c)}^{M,c} \otimes N \right) \circ \beta_c^{-1} \\
&= \mu_{(c,x)}^{M \otimes N,G} \circ \nu_{(c,{\rm id}_c)}^{M \otimes N,c} \\
&= \varphi_c^{M \otimes N,G}(x), 
\end{align*}
where the fourth equality follows from the defining properties of $\alpha_G$ and $\beta_c$. 

\hfill $\Box$ \\

\begin{Ex}\normalfont
Let $\mathcal{C}$ be a small category. 
For any representation $M:\mathcal{C} \to \hat{\mathcal{C}}$ and $\mathbf{y}:\mathcal{C} \to \hat{\mathcal{C}}$, 
let us consider the tensor product $M \otimes \mathbf{y} \cong M$. 
By ${\rm Int}_{M}(-) \otimes \mathbf{y} \cong {\rm Int}_M \cong {\rm Int}_{M \otimes \mathbf{y}}$, 
the morphism $j:\varphi^{M,G} \to \varphi^{M \otimes \mathbf{y},G}$ induced in $\mathcal{IT}_G$ is 
an isomorphism for any $G \in \hat{\mathcal{C}}$. 
\label{ex8}
\end{Ex}

\begin{Ex}\normalfont
Let $\mathcal{D}$ be a bicomplete category. 
Let us consider the tensor product $\mathbf{y} \otimes N \cong N$ where $N:\upuparrows \to \mathcal{D}$ 
is a representation. 

We have ${\rm Int}_{\mathbf{y} \otimes N}(1) \cong {\rm Int}_N(1)$. 
Since ${\rm Int}_{\mathbf{y}}(1)=\emptyset$, 
${\rm Int}_{\mathbf{y}}(1) \otimes N=\emptyset \otimes N$ is an initial object of $\mathcal{D}$. 
Consequently, ${\rm Int}_{\mathbf{y} \otimes N}(1) \not\cong {\rm Int}_{\mathbf{y}}(1) \otimes N$ 
if ${\rm Int}_N(1)$ is not an initial object in $\mathcal{D}$. 
\label{ex9}
\end{Ex}

\begin{Ex}\normalfont
We define the representation $M_{1'} : \upuparrows \to \hat{\upuparrows}$ by 
\begin{align*}
M_{1'}(0)
&= \xymatrix{
p_0 \ar[r]^-{a_0} & p_1 \ar[r]^-{a_1} & p_2
}, \\ 
M_{1'}(1)
&= \xymatrix{
q_0 \ar[r]^-{b_0} & q_1 \ar@<1ex>[r]^-{b_{10}} \ar@<-1ex>[r]_-{b_{11}} & q_2 \ar[r]^-{b_2} & q_3 
}
\end{align*}
and $m_0 \cdot a_0=b_0, \ m_0 \cdot a_1=b_{10}, \ m_1 \cdot a_0=b_{11}$ and $m_1 \cdot a_1=b_2$. 

We have ${\rm Int}_{M_{1'}}(1) \cong \{\bullet \ \bullet \}$ and 
${\rm Int}_{M_{1'}}(1) \otimes M_0 \cong \{\xymatrix{ \bullet \ar[r] & \bullet \ \bullet \ar[r] & \bullet }\}$. 
On the other hand, 
${\rm Int}_{M_{1'} \otimes M_0}(1) \cong \{\xymatrix{ \bullet \ar[r] & \bullet \ar[r] & \bullet }\}$ 
because $M_{1'} \otimes M_0 \cong M_0^{\otimes 3}$. 
\label{ex10}
\end{Ex}

\begin{Ex}\normalfont
One can show $\beta:\xymatrix{ {\rm Int}_{M_0^{\otimes n}} (-) \otimes M_0 \ar[r]^-{\cong} & {\rm Int}_{M_0^{\otimes n+1}} }$
by direct calculation. Therefore, we have a sequence of morphisms 
\[
\xymatrix{ \varphi^{M_0,G} \ar[r] & \varphi^{M_0^{\otimes 2},G} \ar[r] & \cdots \ar[r] & \varphi^{M_0^{\otimes n},G} \ar[r] & \cdots }
\]
in $\mathcal{IT}_G$ for any directed network $G \in \hat{\upuparrows}$. 
In particular, by Example \ref{ex6}, we can take every morphism in the following sequence 
\[
\xymatrix{ \varphi^{M_0,G} \ar[r] & \varphi^{M_0^{\otimes 3},G} \ar[r] & \varphi^{M_0^{\otimes 5},G} \ar[r] & \cdots }
\]
as a non-isomorphism. 
\label{ex11}
\end{Ex}

\subsection{Stability with respect to representations and gluing condition for representations}
In this subsection, we study significance of $M_0$ from a different point of view. 

\begin{Def}
Let $\mathcal{C}$ be a small category, $\mathcal{D}$ a bicomplete category and 
$M:\mathcal{C} \to \mathcal{D}$ a representation. 
We say that a presheaf $G \in \hat{\mathcal{C}}$ is \textbf{stable} with respect to the representation $M$ 
if the unit $\eta^{M,G}:G \to \mathcal{D}(M(-),G \otimes M)$ is an isomorphism
\footnote{A stable directed network $G$ with respect to $M$ is an instance of free algebra of the monad 
induced by the adjunction $(-)\otimes M \dashv \mathcal{D}(M(-_2),-_1)$. }. 
\label{def4}
\end{Def}

Stability of $G \in \hat{\mathcal{C}}$ with respect to a representation $M$ concerns 
when we can recover $G$ from $G \otimes M$, a representation of its dynamic mode. 

\begin{Ex}\normalfont
Let $\mathcal{C}$ be a small category. 
Any presheaf $G \in \hat{\mathcal{C}}$ is stable with respect to the representation 
$\mathbf{y}:\mathcal{C} \to \hat{\mathcal{C}}$. 
\label{ex12}
\end{Ex}

\begin{Ex}\normalfont
A directed network $G \in \hat{\upuparrows}$ is stable with respect to the representation $M_0$ 
if and only if it satisfies the following two conditions \cite{Pultr1979,Haruna2007}:
(i) if $\xymatrix{ \bullet \ar@<1ex>[r]^-{f} \ar@<-1ex>[r]_-{g} & \bullet }$ in $G$ then $f=g$. \\
(ii) If 
\[
\xymatrix{
\bullet \ar[r]^-{f} & \bullet \\
\bullet \ar[ru]^-{g} \ar[r]^-{h} & \bullet 
}
\]
in $G$, then there exists an arc $\xymatrix{\partial_0(f) \ar[r]^-{k} & \partial_1(h)}$ in $G$. 
\label{ex13}
\end{Ex}

\begin{Ex}\normalfont
$G \in \hat{\upuparrows}$ is stable with respect to the representation $M_u$ 
if and only if $G$ is a graph of an equivalence relation. 
\label{ex14}
\end{Ex}

A trivial observation is that the representation $\mathbf{y}:\mathcal{C} \to \hat{\mathcal{C}}$ 
is a minimal representation with respect to stability in the sense that for any representation $M$, 
if $G \in \hat{\mathcal{C}}$ is stable with respect to $M$, then $G$ is stable with respect to 
$\mathbf{y}$ by Example \ref{ex12}. However, if we restrict representations into those with a 
`nice' property, then the representation $M_0$ in turn becomes a minimal one for $\mathcal{C}=\upuparrows$. 

\begin{Def}
Let $\mathcal{D}$ be a bicomplete category. We say that 
a representation $M:\upuparrows \to \mathcal{D}$ 
satisfies the \textbf{gluing condition} if the pullback diagram defining 
${\rm Int}_M(1)$ 
\[
\xymatrix{
{\rm Int}_M(1) \ar[r]^-{\nu_{(0,m_1)}^{M,1}} \ar[d]_-{\nu_{(0,m_0)}^{M,1}} & M(0) \ar[d]^-{M(m_1)} \\
M(0) \ar[r]_-{M(m_0)} & M(1)
}
\]
is also a pushout diagram.
\label{def5}
\end{Def}

Intuitively, in a representation $M$ satisfying the gluing condition, $M(1)$ precisely consists 
of two images of $M(0)$ glued on the interface. 

\begin{Thm}
Let $\mathcal{D}$ be a bicomplete category and 
a representation $M:\upuparrows \to \mathcal{D}$ satisfy the gluing condition. 
If a directed network $G \in \hat{\upuparrows}$ is stable with respect to $M$, 
then $G$ is also stable with respect to $M_0$. 
\label{thm2}
\end{Thm}
\textit{Proof. }
We show that the two conditions (i) and (ii) in Example \ref{ex13} hold for a directed network $G$ when 
$\eta^{M,G}:G \to \mathcal{D}(M(-),G \otimes M)$ is an isomorphism. 
First, we consider (ii). 
Let us assume 
\[
\xymatrix{
\bullet \ar[r]^-{f} & \bullet \\
\bullet \ar[ru]^-{g} \ar[r]^-{h} & \bullet 
}
\]
in $G$. 
Let us consider the diagram 
\[
\xymatrix{
{\rm Int}_M(1) \ar[r]^-{\nu_{(0,m_1)}^{M,1}} \ar[d]_-{\nu_{(0,m_0)}^{M,1}} & M(0) \ar[d]_-{M(m_1)} \ar@/^/[rdd]^-{\eta_0^{M,G}(\partial_1(h))} & \\
M(0) \ar[r]^-{M(m_0)} \ar@/_/[rrd]_-{\eta_0^{M,G}(\partial_0(f))} & M(1) \ar@{.>}[dr]|-{\delta} & \\
& & G \otimes M.
}
\]
Note that ${\rm Int}_M(m_i)=\nu_{(0,m_i)}^{M,1}$ for $i=0,1$. 

If the outside square is commutative, then we obtain a unique morphism $\delta:M(1) \to G \otimes M$ making 
the diagram commute because the inner square is pushout. 
Let us put $k=\left(\eta_1^{M,G}\right)^{-1}(\delta) \in G(1)$. 
Since $\delta \circ M(m_0)=\eta_0^{M,G}(\partial_0(f))$ and $\delta \circ M(m_1)=\eta_0^{M,G}(\partial_1(h))$, 
we have 
$\partial_0(k)=\partial_0 \left(\left(\eta_1^{M,G}\right)^{-1}(\delta)\right)=\left(\eta_0^{M,G}\right)^{-1}(\delta \circ M(m_0))=\partial_0(f)$ 
and 
$\partial_1(k)=\partial_1 \left(\left(\eta_1^{M,G}\right)^{-1}(\delta)\right)=\left(\eta_0^{M,G}\right)^{-1}(\delta \circ M(m_1))=\partial_1(h)$ 
by the naturality of $\left(\eta^{M,G}\right)^{-1}$. 
This means that we get an arc $\xymatrix{\partial_0(f) \ar[r]^-{k} & \partial_1(h)}$ in $G$. 

Thus, it is enough to show that $\eta_0^{M,G}(\partial_0(f)) \circ {\rm Int}_M(m_0)=\eta_0^{M,G}(\partial_1(h)) \circ {\rm Int}_M(m_1)$. 
However, since $\varphi_0^{M,G}=\eta_0^{M,G}$ and 
$\varphi_1^{M,G}=\eta_0^{M,G}(\partial_i(-)) \circ {\rm Int}_M(m_i)$ for $i=0,1$ by the naturality of $\varphi^{M,G}$, 
we have 
$\eta_0^{M,G}(\partial_0(f)) \circ {\rm Int}_M(m_0)=\varphi_1^{M,G}(f)
=\eta_0^{M,G}(\partial_1(f)) \circ {\rm Int}_M(m_1)=\eta_0^{M,G}(\partial_1(g)) \circ {\rm Int}_M(m_1)
=\varphi_1^{M,G}(g)=\eta_0^{M,G}(\partial_0(g)) \circ {\rm Int}_M(m_0)
=\eta_0^{M,G}(\partial_0(h)) \circ {\rm Int}_M(m_0)=\varphi_1^{M,G}(h)
=\eta_0^{M,G}(\partial_1(h)) \circ {\rm Int}_M(m_1)$. 

Next we consider (i). 
Assume that $\xymatrix{ \bullet \ar@<1ex>[r]^-{f} \ar@<-1ex>[r]_-{g} & \bullet }$ in $G$. 
Namely, $\partial_0(f)=\partial_0(g)$ and $\partial_1(f)=\partial_1(g)$ for $f,g \in G(1)$. 
We have $\eta_0^{M,G}(\partial_i(f))=\eta_0^{M,G}(\partial_i(g))$ for $i=0,1$. 
On the other hand, since $\eta_0^{M,G}(\partial_0(f)) \circ {\rm Int}_M(m_0)=\varphi_1^{M,G}(f)
=\eta_0^{M,G}(\partial_1(f)) \circ {\rm Int}_M(m_1)$, we obtain a unique morphism 
$\delta:M(1) \to G \otimes M$ such that $\delta \circ M(m_i)=\eta_0^{M,G}(\partial_i(f))$ for $i=0,1$ 
by the same way as in the proof of (ii) above. 
Now, $\eta_1^{M,G}(f)$ satisfies the conditions for $\delta$, it follows that 
$\delta=\eta_1^{M,G}(f)$ by the uniqueness. Since $\partial_i(f)=\partial_i(g)$ for $i=0,1$, 
$\eta_1^{M,G}(g)$ also satisfies the condition. Thus, we also have $\delta=\eta_1^{M,G}(g)$. 
Consequently, we have $\eta_1^{M,G}(f)=\eta_1^{M,G}(g)$. 
Since $\eta^{M,G}$ is an isomorphism, the equality $f=g$ holds. 

\hfill $\Box$\\

The next proposition provides a way to produce a new representation satisfying the gluing condition 
from an old one. 

\begin{Pro}
Let $\mathcal{D}$ be a bicomplete category. If a representation $M:\upuparrows \to \hat{\upuparrows}$ 
satisfies the gluing condition, then for any representation $N:\upuparrows \to \mathcal{D}$, 
the tensor product representation $M \otimes N$ also satisfies the gluing condition. 
\label{pro4}
\end{Pro}
\textit{Proof. }
Since $M$ satisfies the gluing condition, 
\[
\xymatrix{
{\rm Int}_M(1) \ar[r]^-{\nu_{(0,m_1)}^{M,1}} \ar[d]_-{\nu_{(0,m_0)}^{M,1}} & M(0) \ar[d]^-{M(m_1)} \\
M(0) \ar[r]_-{M(m_0)} & M(1)
}
\]
is a pushout diagram. Since $(-) \otimes N$ preserves colimits, the outside square of 
the following commutative diagram is also a pushout: 
\[
\xymatrix{
{\rm Int}_M(1) \otimes N \ar@/^/[rrd]^-{\nu_{(0,m_1)}^{M,1} \otimes N} \ar@/_/[rdd]_-{\nu_{(0,m_0)}^{M,1} \otimes N} & & \\
& {\rm Int}_{M \otimes N}(1) \ar[r]_-{\nu_{(0,m_1)}^{M \otimes N,1}} \ar[d]^-{\nu_{(0,m_0)}^{M \otimes N,1}} & M(0) \otimes N \ar[d]^-{M(m_1) \otimes N} \\
& M(0) \otimes N \ar[r]_-{M(m_0) \otimes N}  & M(1) \otimes N. 
}
\]
The inner square is a pullback by the definition of interface. 
Then, we can show that the inner square is automatically a pushout, which proves that 
$M \otimes N$ satisfies the gluing condition. 

Indeed, let us consider a diagram 
\[
\xymatrix{
E \ar@/^/[rrd] \ar@/_/[rdd] \ar@{.>}[dr]|-{\delta} & & & \\
& D \ar[r] \ar[d] & B \ar[d] \ar@/^/[rdd] & \\
& A \ar[r] \ar@/_/[rrd] & C \ar@{.>}[dr]|-{\delta'} & \\
& & & F
}
\]
in an arbitrary category $\mathcal{A}$. 
Let us assume that 
\[
\xymatrix{
D \ar[r] \ar[d] & B \ar[d] \\
A \ar[r] & C
}
\]
is a pullback and 
\[
\xymatrix{
E \ar[r] \ar[d] & B \ar[d] \\
A \ar[r] & C
}
\]
is a pushout. 
By the universality of the pullback, there exists a unique morphism 
$\delta:E \to D$ making the appropriate parts of the diagram commute. 
If 
\[
\xymatrix{
D \ar[r] \ar[d] & B \ar[d] \\
A \ar[r] & F
}
\]
commutes, then 
\[
\xymatrix{
E \ar[r] \ar[d] & B \ar[d] \\
A \ar[r] & F
}
\]
also commutes. By the universality of the pushout, there exists a unique morphism 
$\delta':C \to F$ making the appropriate parts of the diagram commute. 
This means that 
\[
\xymatrix{
D \ar[r] \ar[d] & B \ar[d] \\
A \ar[r] & C
}
\]
is a pushout diagram. 

\hfill $\Box$\\

One can check that $M_0$ satisfies the gluing condition easily. 
Consequently, by Proposition \ref{pro4}, we obtain a representation $M_0 \otimes N$ 
satisfying the gluing condition from $M_0$ and any representation 
$N:\upuparrows \to \mathcal{D}$ where $\mathcal{D}$ is a bicomplete category. 

\section{Static and Dynamic Modes in Real Biological Networks}
In Section 2, we have shown that the pair of two network transformations (adjoint functors) $R$ and $L$ considered in Introduction 
is a canonical one associated with the idea ``interaction is interface between functions''. 
$R$ is associated with the static mode and $L$ is with the dynamic mode. Natural path notions for the two modes 
have turned out to be the directed path and the lateral path, respectively. 
The correspondence among the two modes, functors, path notions and their functionality 
is summarized in the table below. 

\begin{center}
\begin{tabular}{cccc} \toprule
mode & functor & path & functionality \\ \midrule
dynamic & $L$ & lateral & coherence \\
static & $R$ & directed & transport \\ \bottomrule
\end{tabular}
\end{center}

In this section, we study how the 
static and dynamic modes are embedded in real biological networks by measuring importance of arcs with respect to 
each path notion. The measure we use is \textbf{betweenness centrality} \cite{Anthonisse1971,Freeman1977}. 
In general, betweenness centrality measures importance of a node or an arc in a given network 
by counting the number of geodesic paths that pass through the node or the arc. 
In this section and next section, we only consider simple directed networks, namely, 
directed networks without self-loops and multi-arcs. 

Let $G=(A,N,\partial_0,\partial_1)$ be a simple directed network and $f \in A$. 
The \textbf{lateral betweenness centrality (LBC)} of $f$ is defined by 
\[
{\rm LBC}_f=\frac{1}{\sum_{g,h \in A, \ l_{gh}>0} (s_{gh}+1)} \sum_{g,h \in A, \ l_{gh}>0} \frac{l_{gh}^f}{l_{gh}}, 
\]
where $l_{gh}$ is the number of geodesic lateral paths from $g$ to $h$, $l_{gh}^f$ is the number of 
geodesic lateral paths from $g$ to $h$ that pass through $f$ and $s_{gh}$ is the geodesic lateral distance 
from $g$ to $h$. 
Similarly, the \textbf{directed betweenness centrality (DBC)} of $f$ is defined by 
\[
{\rm DBC}_f=\frac{1}{\sum_{g,h \in A, \ d_{gh}>0} (t_{gh}+1)} \sum_{g,h \in A, \ d_{gh}>0} \frac{d_{gh}^f}{d_{gh}}. 
\]
where $d_{gh}$ is the number of geodesic directed paths from $g$ to $h$, $d_{gh}^f$ is the number of 
geodesic directed paths from $g$ to $h$ that pass through $f$ and $t_{gh}$ is the geodesic directed distance 
from $g$ to $h$. 
Here, we normalize both ${\rm LBC}_f$ and ${\rm DBC}_f$ so that 
$\sum_{f \in A}{\rm LBC}_f=1$ and $\sum_{f \in A}{\rm DBC}_f=1$. 
LBC can be seen as a measure of importance with respect to coherence of the network because 
the lateral path is associated with gluing of functions. On the other hand, DBC is a measure of 
importance with respect to transport on the network. 

\begin{figure}
\centering
\begin{tabular}{cc}
\includegraphics[width=5.5cm]{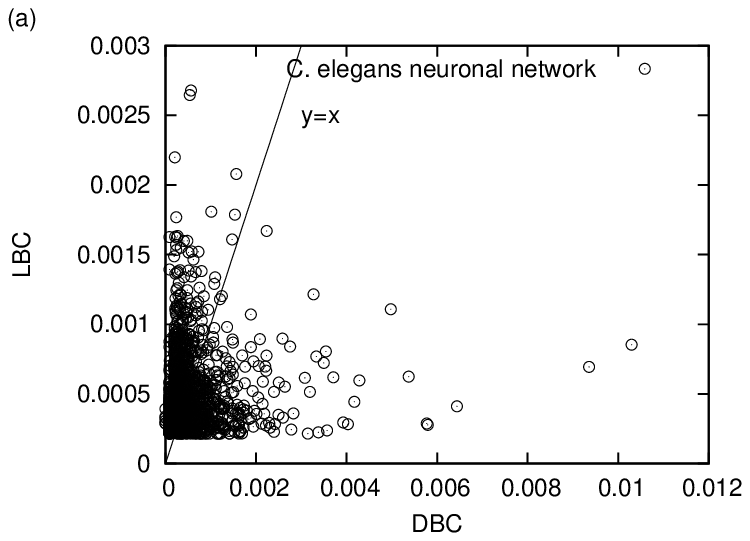} &
\includegraphics[width=5.5cm]{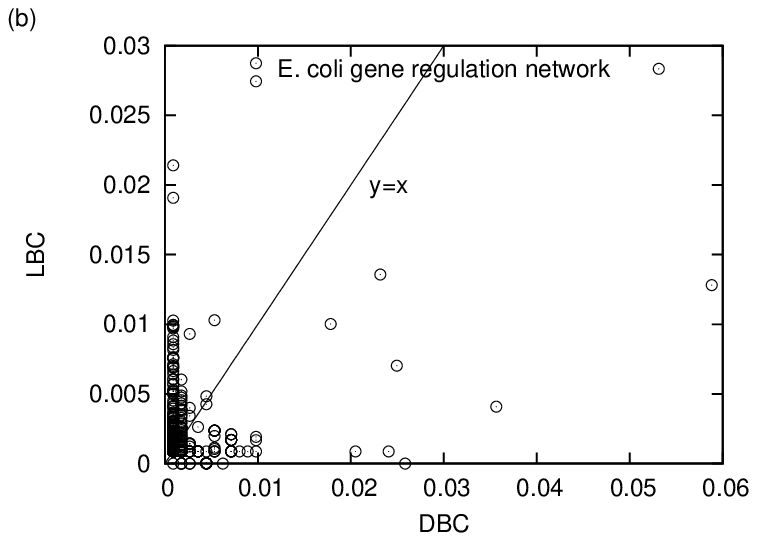} \\
\includegraphics[width=5.5cm]{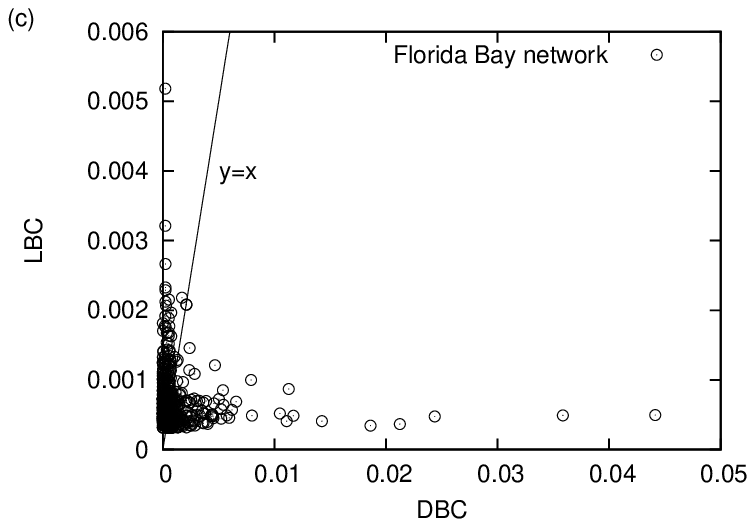} &
\includegraphics[width=5.5cm]{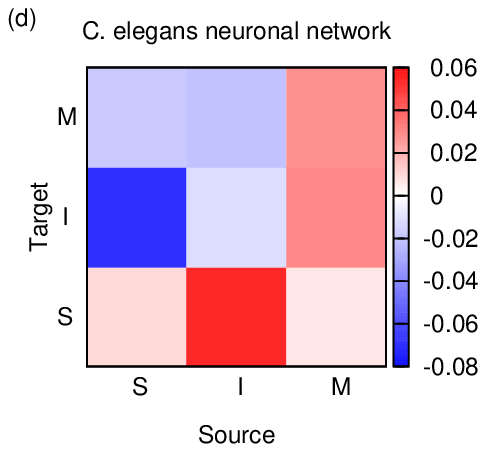}
\end{tabular}
\caption{
(a)-(c) The relationship between DBC value and LBC value for each arc in 
the three biological networks is shown. 
(d) The relationship between types of neurons and dominance of the two betweenness centralities. 
For each of nine blocks, a group of arcs is assigned based on types of their source neurons (horizontal axis, 
S: sensory neuron, I: inter neuron, M: motor neuron) and their target neurons (vertical axis). 
The color of each block indicates the value of the sum of DBC value minus LBC value over 
all arcs in the corresponding group. 
}
\label{fig1}
\end{figure}

In Figure \ref{fig1} (a)-(c), we show the relationship between DBC and LBC for three directed biological networks: 
a neuronal network of \textit{C. elegans} \cite{White1986,Varshney2011} (nodes: neurons, arcs: chemical synapses, $|N|=279$, $|A|=2194$), 
a gene transcription regulation network of \textit{E. coli} \cite{Shen-Orr2002}
 (nodes: genes (operons), arcs: regulation relations, $|N|=328$, $|A|=456$) and 
an ecological flow network of Florida Bay \cite{Ulanowicz1998} (nodes: taxa, arcs: carbon flows, $|N|=121$, $|A|=1767$)
\footnote{
Three nodes classified as detritus, one node corresponding to roots and arcs whose source or target are 
one of these four nodes are removed from the original data to focus on prey-predator interactions. 
}. 
Each point in each figure corresponds to an arc in the network. 
In any of the three networks, we can see a trade-off relationship between the values of LBC and DBC. 
Namely, given an arc, if its LBC is large, then its DBC is small and vice versa. Thus, 
any arc in these networks cannot have large values of LBC and DBC at the same time. 

Figure \ref{fig1} (d) shows the relationship between types of neurons 
(sensory, inter and motor) and average dominance of LBC or DBC in the neuronal network of \textit{C. elegans}. 
Each arc is classified into one of the nine groups depending on types of its source and target neurons. For each group, 
the sum of DBC value minus LBC value over all arcs in the group is shown as a color. 
We can see that LBC is dominant over DBC on average in the forward direction such as connections from sensory to inter or motor and 
from inter to motor neurons. On the other hand, DBC is dominant over LBC on average in the backward directions such as 
connections from motor to sensory or inter and from inter to sensory neurons. 
Thus, in this case, there is a clear relationship between types of synaptic connections and dominance of the two 
betweenness centralities. 

\begin{figure}
\centering
\begin{tabular}{cc}
\includegraphics[width=5.5cm]{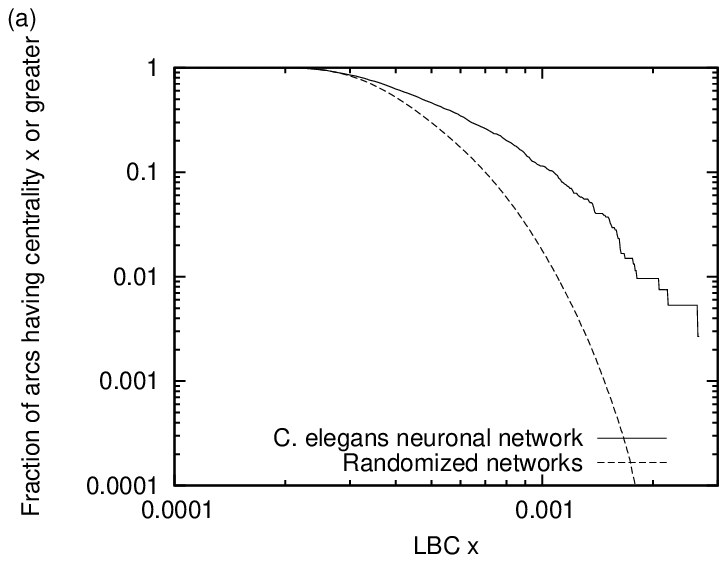} &
\includegraphics[width=5.5cm]{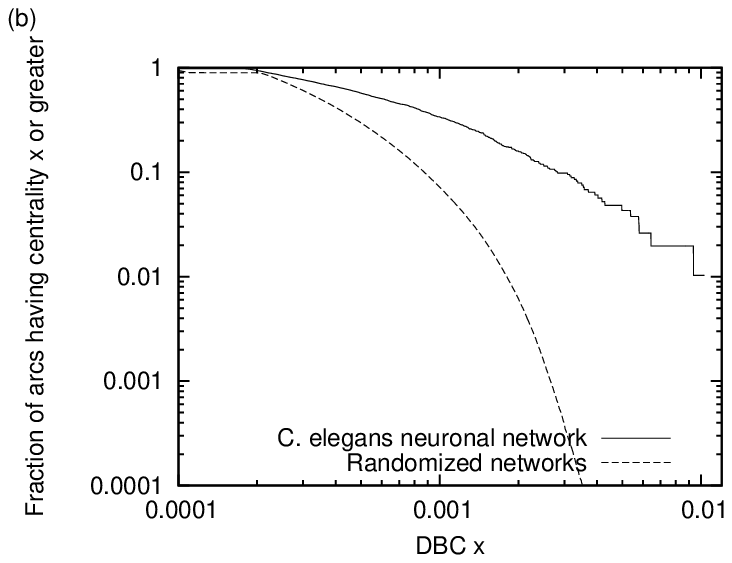} \\
\includegraphics[width=5.5cm]{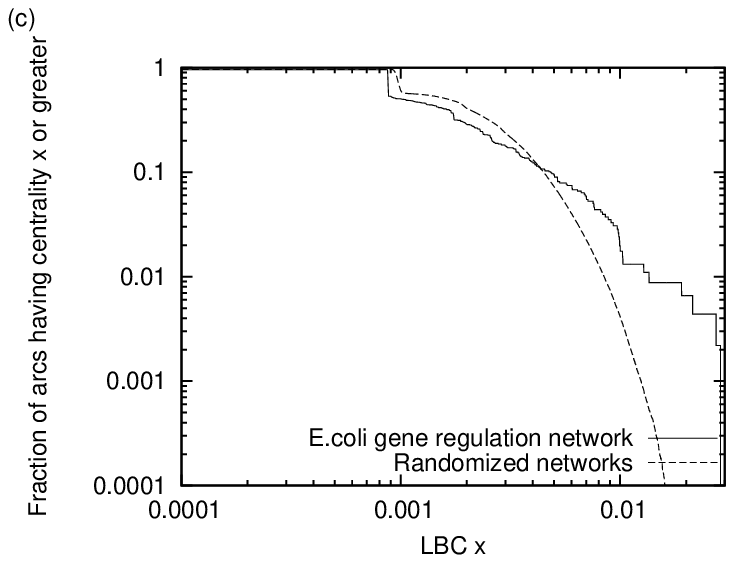} &
\includegraphics[width=5.5cm]{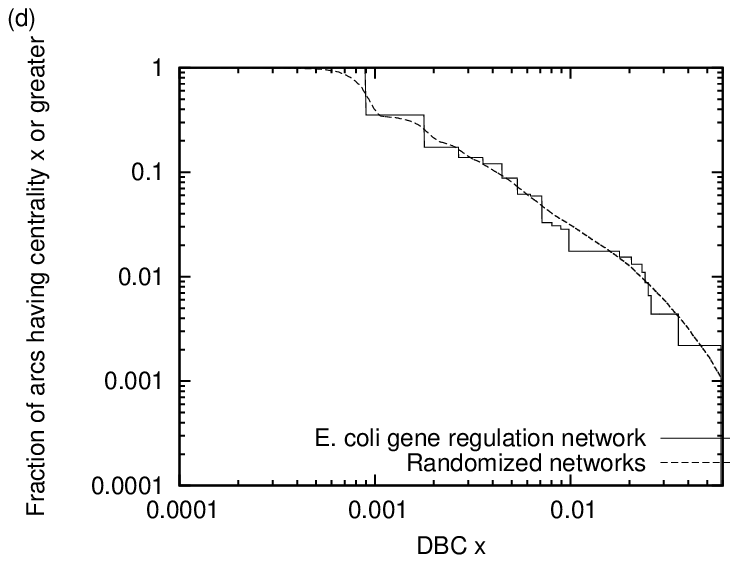} \\
\includegraphics[width=5.5cm]{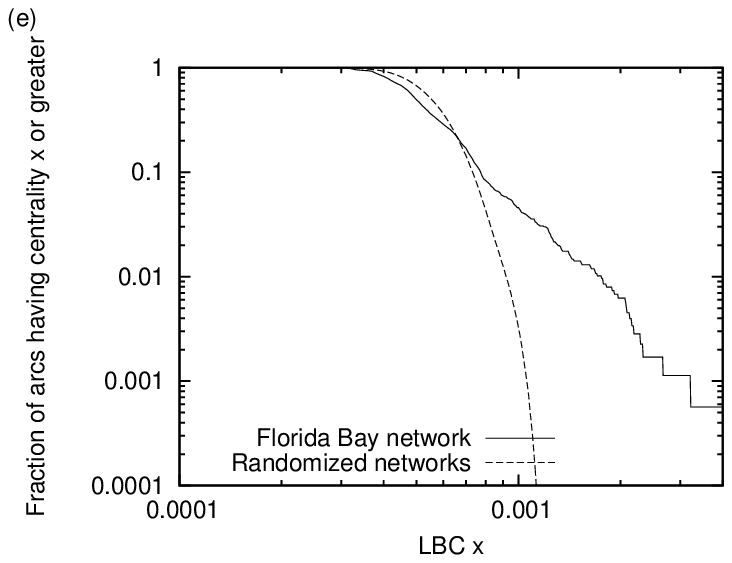} &
\includegraphics[width=5.5cm]{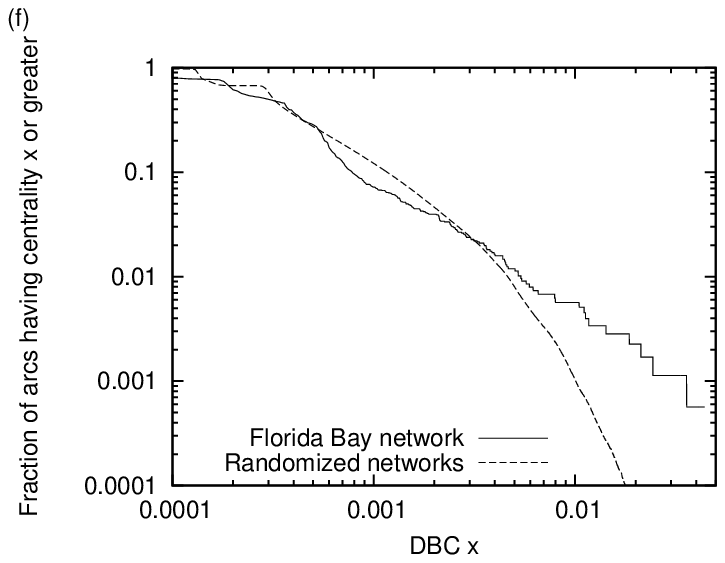} 
\end{tabular}
\caption{
For each of the three biological networks, the cumulative frequency distribution functions of 
LBC and DBC values and the averages of those for randomized networks preserving 
the degree of each node. The average is taken over 1000 randomized networks. 
(a) LBC for the \textit{C. elegans} neuronal network. 
(b) DBC for the \textit{C. elegans} neuronal network. 
(c) LBC for the \textit{E. coli} gene regulation network. 
(d) DBC for the \textit{E. coli} gene regulation network. 
(e) LBC for the Florida Bay ecological flow network. 
(f) DBC for the Florida Bay ecological flow network. 
}
\label{fig2}
\end{figure}

To see whether this kind of trade-off relationship is significant, we compare cumulative frequency distribution functions 
of LBC and DBC values for real networks with those for randomized networks (Figure \ref{fig2}). 
In preparing a randomized network for each real network, we repeatedly swap the targets of two randomly chosen arcs. 
Thus, the degree of each node is preserved in the randomized network. 
Clearly, the cumulative frequency distribution functions for real networks have longer tails than the averages of those for randomized networks 
except for that of DBC for the gene transcription regulation network. This suggests that the trade-off relationship 
between LBC and DBC seen in Figure \ref{fig1} is not trivial. 

The trade-off relationship between LBC and DBC values indicates that a ``division of labor'' between 
coherence of the network and transport on the network have emerged through the evolution of 
the biological networks. In the next section, we discuss how the trade-off relationship emerges 
through an optimization model of complex networks. 

\section{An Optimization Model of Directed Networks}
In the previous section, we found the trade-off relationship between LBC and DBC values in three 
biological networks such that if an arc has a large value in one of the two betweenness centralities, 
then the arc has a small value in the other. The comparison with the null-model suggests that 
this trade-off relationship is not an accidental one. In this section, we propose a simple optimization model 
of complex networks to study how the trade-off relationship emerges through evolution. 

In general, an optimization model of complex networks sets a quality function to be optimized. 
In the course of evolution, the topology of a network is changed gradually to the direction 
in which the quality function increases or decreases. Typically, quality functions involve 
two quantities that trade-off each other. 
For example, \cite{Cancho2003} took a quality function as a linear combination of 
average geodesic distance and a normalized number of edges to express the trade-off between 
transport efficiency and cost to make edges in a network. \cite{Brede2009} considered 
a linear combination of average geodesic distance and the largest eigenvalue of adjacency matrix. 
In this case, the trade-off between transport efficiency and robustness of a network is 
the target of the optimization. In our model, we take a quality function incorporating two 
efficiencies \cite{Latora2001} with respect to the lateral path and the directed path. 
This choice is intended to see what kind of networks evolve 
under a specified balance between intensity of the dynamic mode and that of the static mode in a network. 
Namely, our target of optimization is the hypothesized trade-off between coherence and transport efficiency 
of networks. 

For a simple directed network $G=(A,N,\partial_0,\partial_1)$, we define its 
\textbf{lateral efficiency} by 
\[
{\rm Eff}_L=\frac{1}{|A|(|A|-1)}\sum_{f,g \in A, \ f \neq g} \frac{1}{s_{fg}}, 
\]
where $|A|$ is the number of arcs and $s_{fg}$ is the geodesic lateral distance between two arcs $f$ and $g$ 
as in the previous section. 
${\rm Eff}_L$ is a measure of the intensity of the dynamic mode or coherence of the network. 
Similarly, \textbf{directed efficiency} of $G$ is 
defined by 
\[
{\rm Eff}_D=\frac{1}{|A|(|A|-1)}\sum_{f,g \in A, \ f \neq g} \frac{1}{t_{fg}}, 
\]
where $t_{fg}$ is the geodesic directed distance between two arcs $f$ and $g$. 
${\rm Eff}_D$ is a measure of the intensity of the static mode or transport efficiency on the network. 
The following quality function is used to evolve networks: 
\[
Q(\lambda)=\lambda{\rm Eff}_D+(1-\lambda){\rm Eff}_L, \ \ 0 \leq \lambda \leq 1. 
\]
$\lambda$ is the parameter which controls the balance of contributions to the quality function 
from the two efficiencies. The evolutionary algorithm is given as follows: 
\begin{description}
\item[(1)]
Initial networks are given by random networks with $n$ nodes and the probability of arcs 
between nodes $p$. 
In the following evolutionary process, the number of arcs $a$ is preserved. 
\item[(2)]
Remove a randomly chosen arc and add an arc randomly. 
\item[(3)]
If the quality function $Q$ increases in (2), then accept the result, otherwise repeat (2). 
\item[(4)]
Stop when $Q$ does not increase $a$ (the number of arcs) times successively in (3). 
\end{description}
Although we can incorporate other factors such as growing process into the evolutionary algorithm \cite{Brede2011}, 
here we consider one of the simplest model to study the basic performance of the evolutionary algorithm 
using the quality function $Q$. 
The results in what follows are obtained from numerical simulations with $n=100$ and $p=0.05$. 
$\lambda$ is incremented from 0 to 1 by 0.05. 
For each value of $\lambda$, we generated 100 optimized networks evolved from different initial networks. 

\begin{figure}
\centering
\begin{tabular}{cc}
\includegraphics[width=5.5cm]{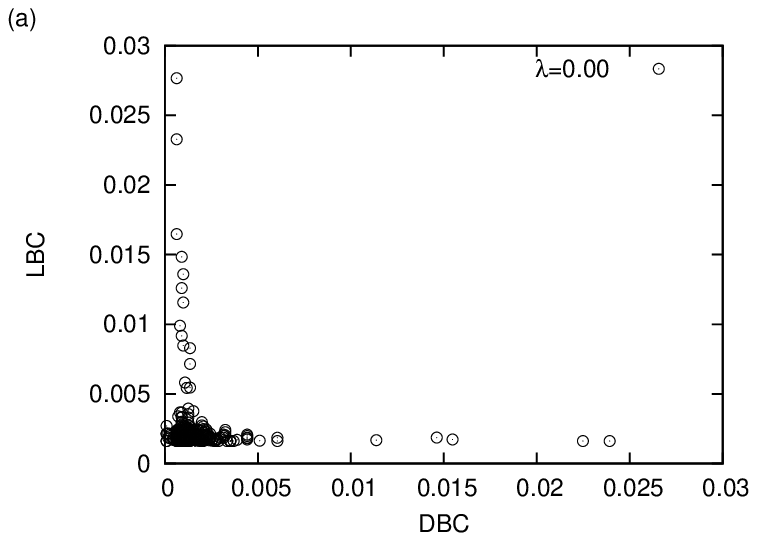} &
\includegraphics[width=5.5cm]{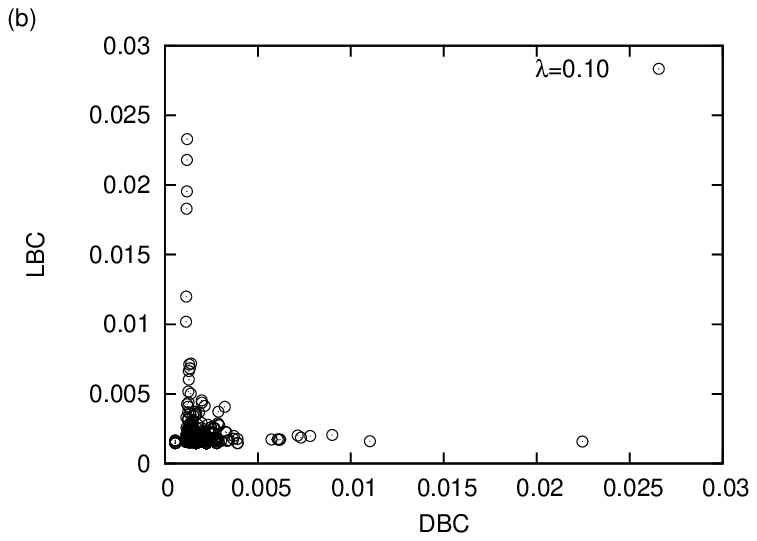} \\
\includegraphics[width=5.5cm]{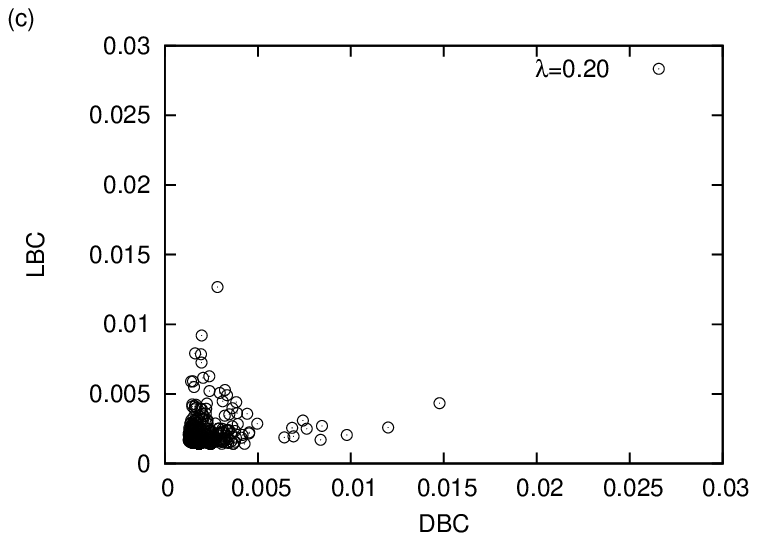} &
\includegraphics[width=5.5cm]{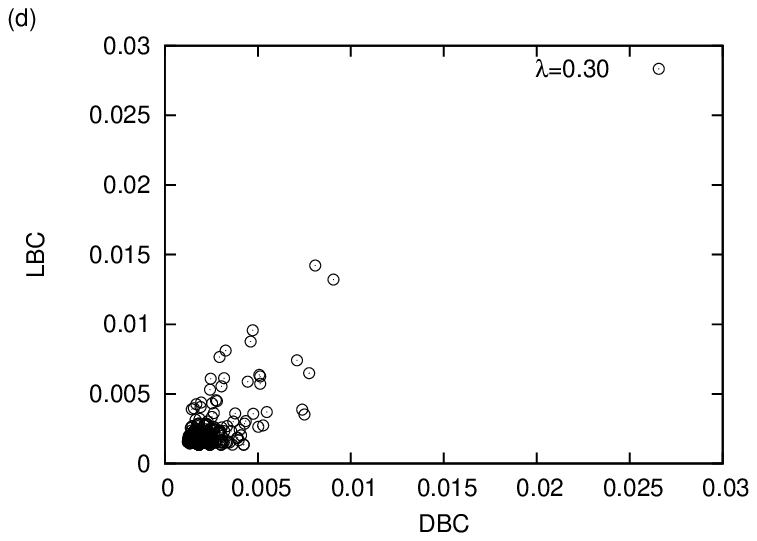} \\
\includegraphics[width=5.5cm]{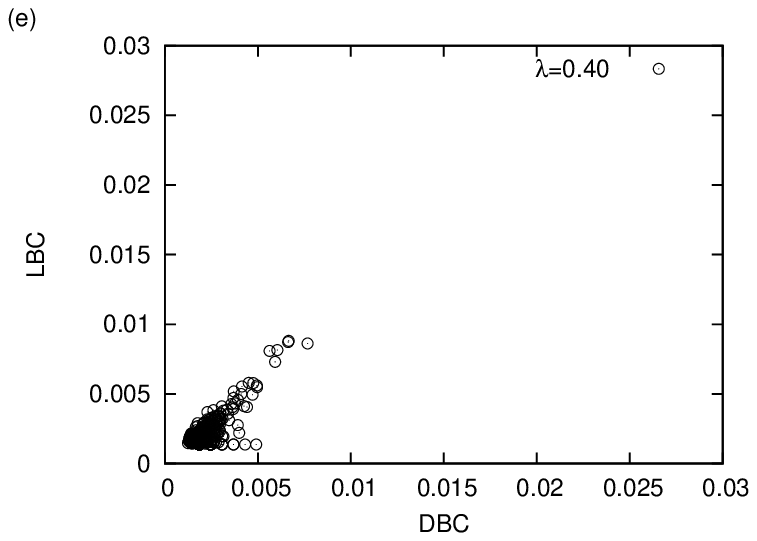} &
\includegraphics[width=5.5cm]{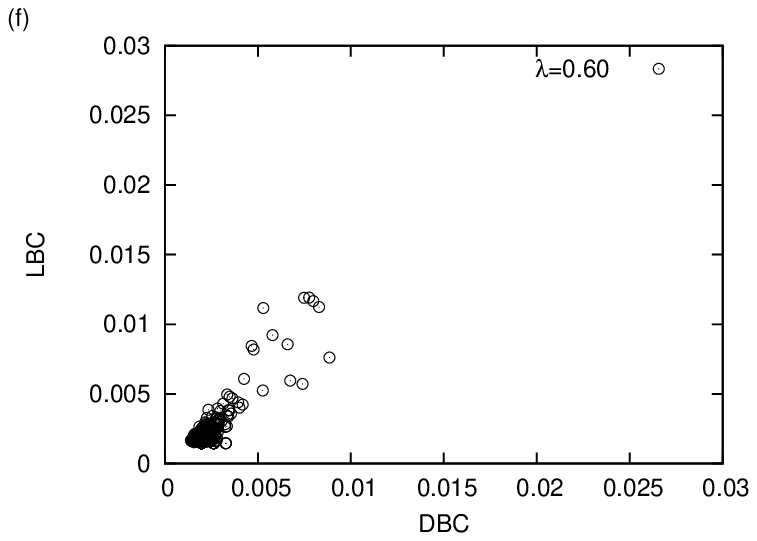} \\
\includegraphics[width=5.5cm]{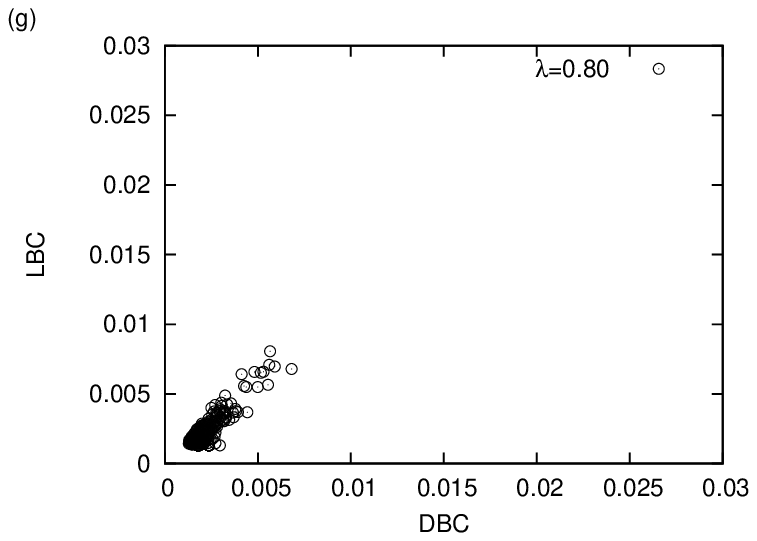} &
\includegraphics[width=5.5cm]{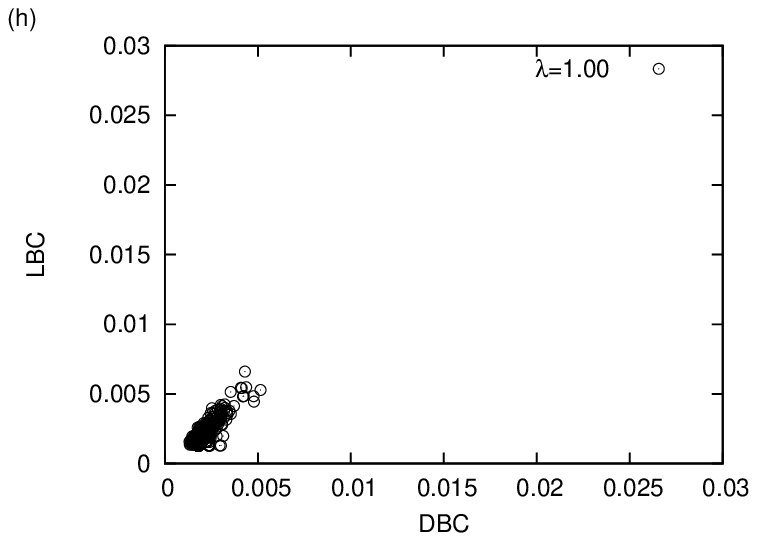} 
\end{tabular}
\caption{
The relationship between LBC and DBC for some values of $\lambda$ in an optimized network. 
(a) $\lambda=0.00$. 
(b) $\lambda=0.10$. 
(c) $\lambda=0.20$. 
(d) $\lambda=0.30$. 
(e) $\lambda=0.40$. 
(f) $\lambda=0.60$. 
(g) $\lambda=0.80$. 
(h) $\lambda=1.00$. 
}
\label{fig3}
\end{figure}

\begin{figure}
\centering
\begin{tabular}{c}
\includegraphics[width=7cm]{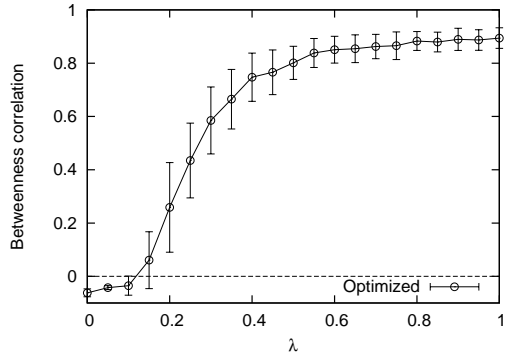} 
\end{tabular}
\caption{
The averaged value of the correlation coefficient between LBC and DBC for each $\lambda$. 
The averages are taken over the 100 optimized networks. Bars are standard deviations. 
}
\label{fig4}
\end{figure}

In Figure \ref{fig3}, the relationship between LBC and DBC for an optimized network for each value of 
$\lambda$ is shown. We can see that the trade-off relationship between LBC and DBC holds for 
the values of $\lambda$ up to around 0.2. In Figure \ref{fig4}, we plot the correlation coefficient 
between LBC and DBC averaged over the 100 optimized networks for each $\lambda$. 
It is negative for $0 \leq \lambda \leq 0.10$ and positive for $\lambda \geq 0.15$. 
A negative correlation coefficient between LBC and DBC indicates the trade-off relationship 
between them. From visual inspection into Figure \ref{fig3}, 
one might suppose that the trade-off relationship between LBC and DBC would hold 
for $\lambda=0.15, 0.20$. However, there are many arcs that have small values for both LBC and DBC. 
Correlation coefficients for these values of $\lambda$ become non-negative due to the existence of these arcs. 
Thus, the correlation coefficient may not necessarily a good indicator for the trade-off relationship. 
However, at least we can conclude that dominance of ${\rm Eff}_L$ in the quality function $Q$ results in 
the trade-off relationship between LBC and DBC. 

\begin{figure}
\centering
\begin{tabular}{cc}
\includegraphics[width=5.5cm]{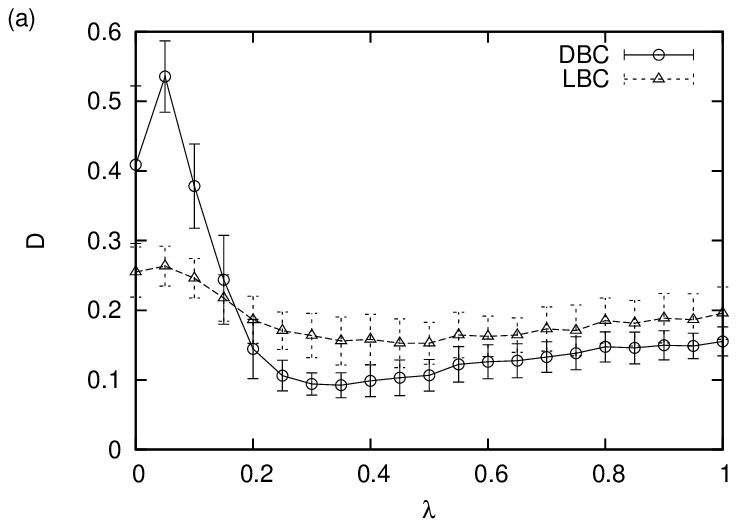} &
\includegraphics[width=5.5cm]{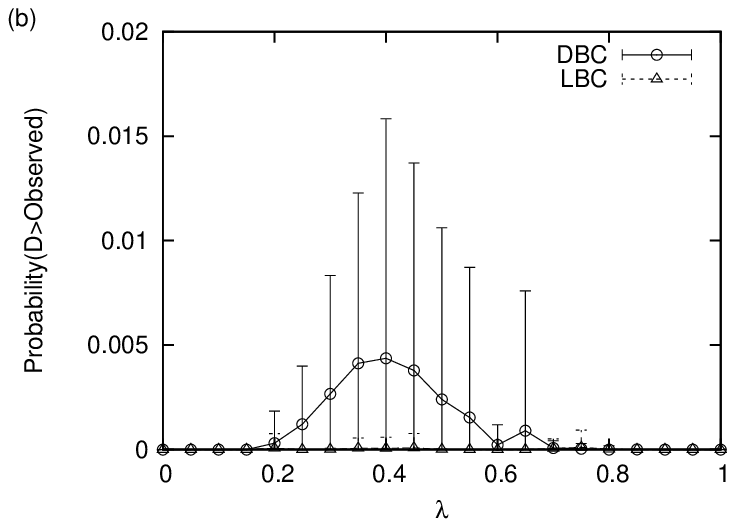} 
\end{tabular}
\caption{
The result of the Kolmogorov-Smirnov test for distribution functions of LBC and DBC 
in optimized networks. (a) $D$ statics and (b) $p$-values. See the main text for details. 
}
\label{fig5}
\end{figure}

We apply the Kolmogorov-Smirnov test for the values of LBC and DBC for each optimized network to see 
whether their distributions are statistically significant. Here, we take reference cumulative distribution functions for 
LBC and DBC as the averaged cumulative frequency distribution functions for them over 100 randomly rewired networks preserving 
the degree of each node, respectively. The averaged value of $D$ static and $p$-values are shown in 
Figure \ref{fig5}. In the worst case, the averaged $p$-value for DBC takes values up to 0.005 around $\lambda=0.40$ 
and here the distributions of $p$-value have large standard deviations. However, on average, we can reject the null-hypothesis 
at level 0.005 for any $\lambda$ and for both LBC and DBC. These results indicate that the distributions of LBC and DBC 
in optimized networks are significantly different from those in the randomized networks. 

We apply standard complex network analysis \cite{Newman2010} to the undirected networks obtained by 
ignoring the direction of arcs in the optimized networks. By abuse of language, we also refer to 
these undirected networks as the optimized networks. 

\begin{figure}
\centering
\begin{tabular}{cc}
\includegraphics[width=5.5cm]{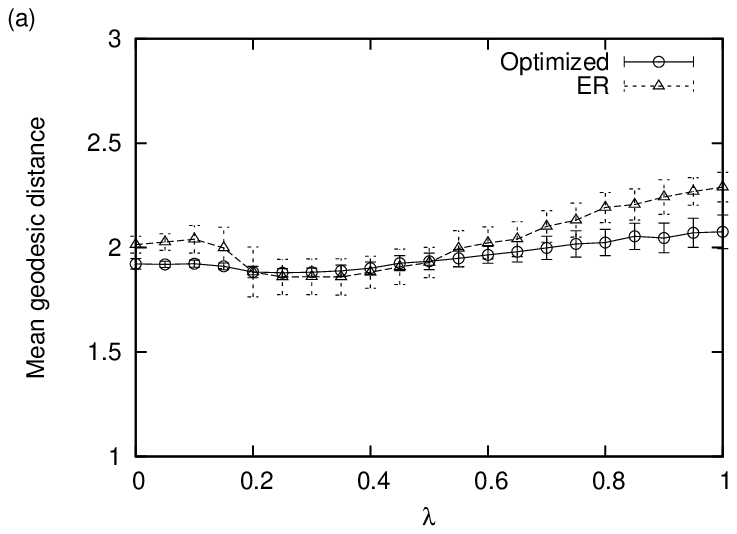} &
\includegraphics[width=5.5cm]{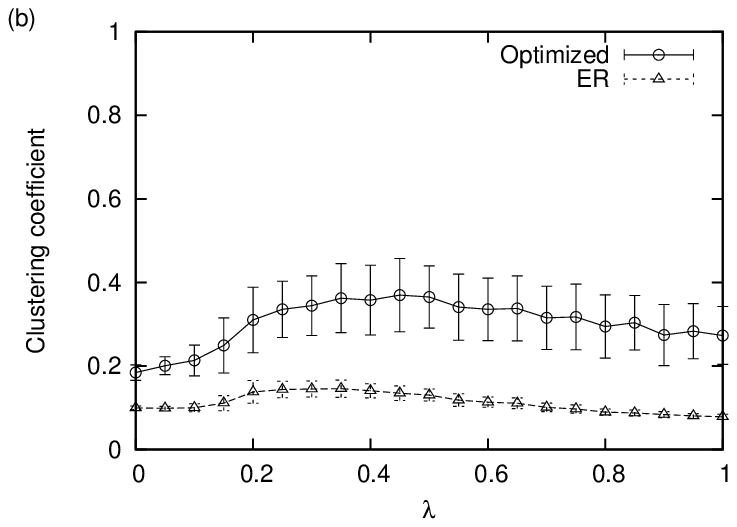} 
\end{tabular}
\caption{
(a) The mean geodesic distance between nodes and (b) the clustering coefficient for 
each $\lambda$ averaged over the 100 optimized networks. Bars are standard deviations. 
Legends with subscript ER are values for corresponding Erd\"{o}s-R\'{e}nyi random networks. 
See the main text for details. 
}
\label{fig6}
\end{figure}

\begin{figure}
\centering
\begin{tabular}{c}
\includegraphics[width=7cm]{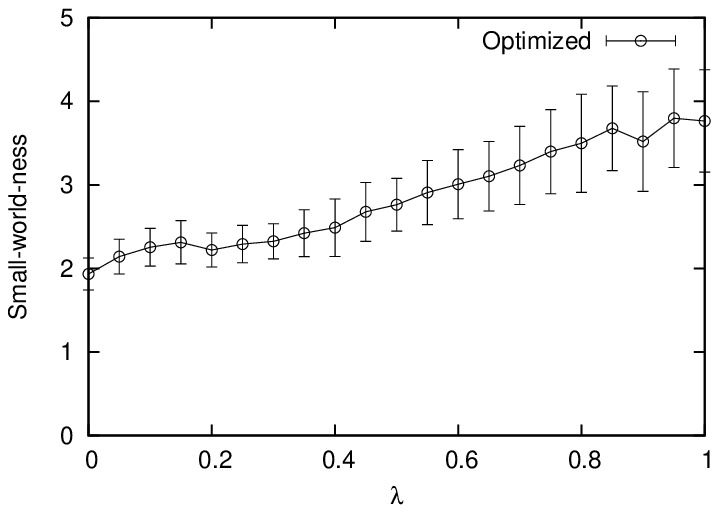} 
\end{tabular}
\caption{
The Small-world-ness for each $\lambda$ averaged over the 100 optimized networks. 
Bars are standard deviations. 
}
\label{fig7}
\end{figure}

Figure \ref{fig6} shows the mean geodesic distance between nodes $l$ and the clustering coefficient $C$ 
(which is degree of transitivity in a network defined by the number of triangles times 6 divided by the number 
of paths of length 2) averaged over the 100 optimized networks for each $\lambda$. 
As references, we also plot the clustering coefficient $C_{ER}=\langle k \rangle/n$ and the mean geodesic distance 
$l_{ER}=\ln(n)/\ln(\langle k \rangle)$ averaged over 100 Erd\"{o}s-R\'{e}nyi random networks with the same 
number of nodes and edges as the optimized networks, where $n$ is the number of nodes and $\langle k \rangle$ is 
mean degree. From Figure \ref{fig6}, the mean geodesic distances for the optimized networks are as small as $l_{ER}$ 
and the clustering coefficients for them are at least about two times larger than $C_{ER}$ for any value of $\lambda$. 
These results suggest that the optimized networks are small-world networks \cite{Watts1998}. 
In Figure \ref{fig7}, we show the average value of the Small-world-ness \cite{Humphries2008} which is a 
quantitative measure of the degree of small-world for each $\lambda$. The Small-world-ness for a network 
with clustering coefficient $C$ and mean geodesic distance $l$ is defined by 
\[
S=\frac{l_{ER}}{l}\frac{C}{C_{ER}}, 
\]
where $l_{ER}$ and $C_{ER}$ are as above. It is proposed that if $S>$1, then we say that the network is 
small-world. Since the value of $S$ for any $\lambda$ exceeds about 2 in Figure \ref{fig7}, 
it is supported that the optimized networks have the small-world property also from this point of view. 

\begin{figure}
\centering
\begin{tabular}{cc}
\includegraphics[width=5.5cm]{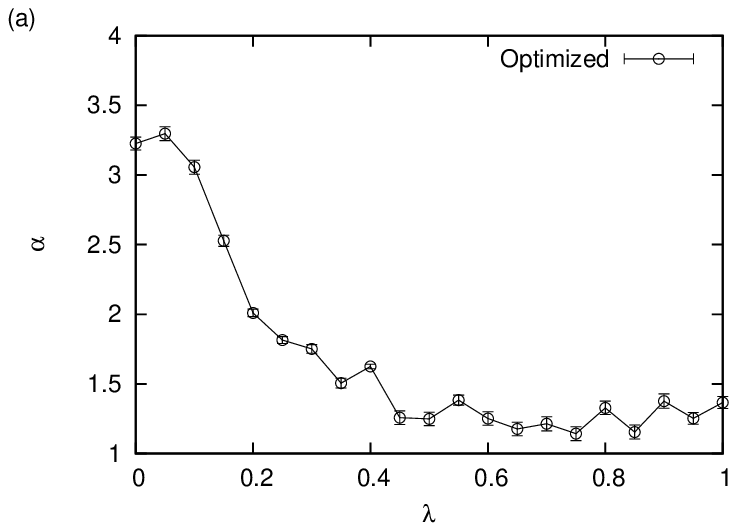} &
\includegraphics[width=5.5cm]{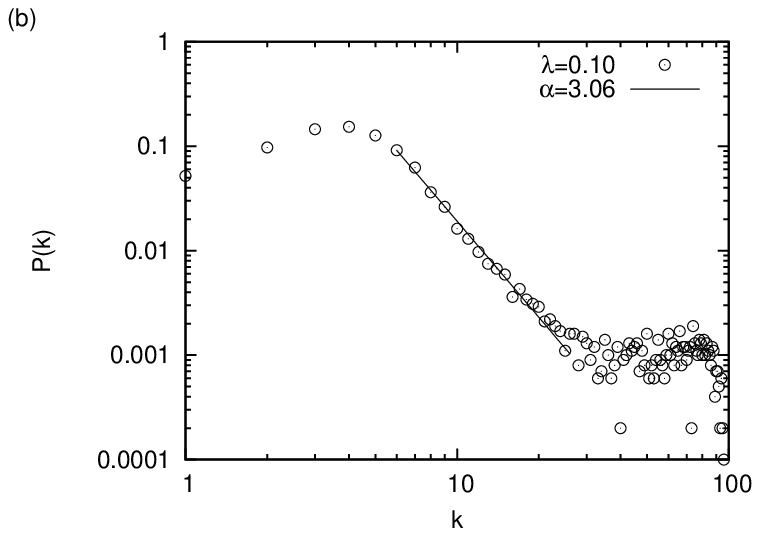} 
\end{tabular}
\caption{
(a) The value of exponent for the power law approximation of the average degree distribution 
of the optimized networks for each $\lambda$. 
(b) The average degree distribution of the optimized networks for $\lambda=0.10$. 
}
\label{fig8}
\end{figure}

Next, we focus on degree distributions of the optimized networks. 
In Figure \ref{fig8} (a), we show the average value of exponent for the power law function which best approximates 
the average degree distribution $P(k)$ of the 100 optimized networks for each $\lambda$. 
We estimated the exponents by applying the maximum likelihood method \cite{Clauset2009} to 
the power law function 
\[
P(k)=k^{-\alpha}/\zeta(\alpha,k_{min},k_{max}), \ \ \zeta(\alpha,k_{min},k_{max})=\sum_{k=k_{min}}^{k_{max}} k^{-\alpha}, 
\]
where $[k_{min},k_{max}]$ is the interval on which we assumed the power law. 
The interval $[k_{min},k_{max}]$ was also estimated by the method in \cite{Clauset2009}. 

From Figure \ref{fig8} (a), we can see that $\alpha$ has a decreasing trend as $\lambda$ increases. 
It is remarkable that $\alpha$ takes the values about 2 to about 3 which are frequently observed 
in real world networks whose degree distributions obey power laws when $\lambda \leq 0.2$. 

\begin{figure}
\centering
\begin{tabular}{cc}
\includegraphics[width=5.5cm]{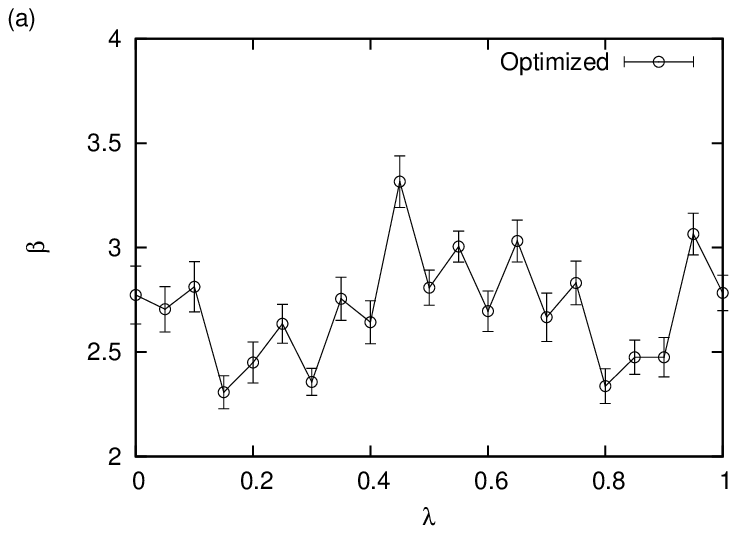} &
\includegraphics[width=5.5cm]{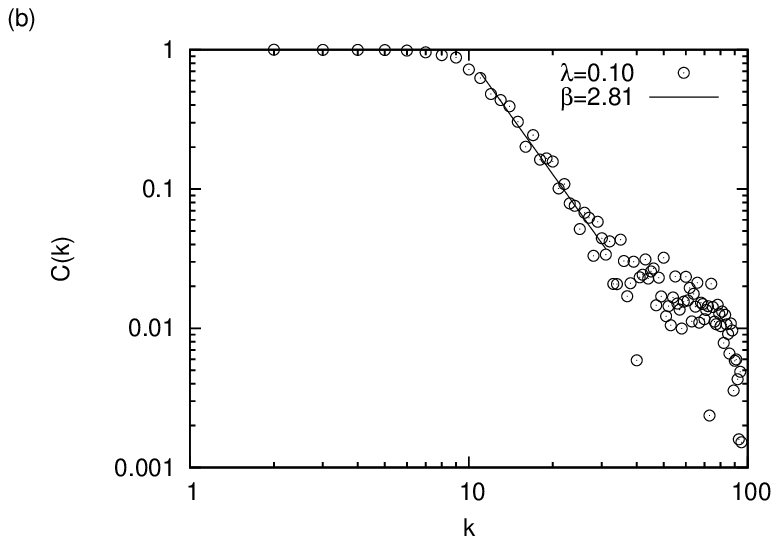} 
\end{tabular}
\caption{
Degree dependency of local clustering coefficients $C(k)$ averaged over the 100 optimized networks 
for each $\lambda$. (a) The exponents $\beta$ of the approximation $C(k) \propto k^{-\beta}$ on 
estimated ranges are shown. (b) $C(k)$ for $\lambda=0.10$. See the main text for details. 
}
\label{fig9}
\end{figure}

The result of modularity scaling analysis based on degree dependency of local clustering coefficients 
is shown in Figure \ref{fig9}. 
The local clustering coefficient $C_i$ for a node $i$ is defined as the number of edges between 
neighbors of $i$ divided by the number of unordered pairs of neighbors of $i$. 
$C(k)$ is the average local clustering coefficients over all nodes with degree $k$. 
Figure \ref{fig9} (a) shows that the exponent $\beta$ of the power law approximation $k^{-\beta}$ 
for $C(k)$ averaged over the 100 optimized networks for each $\lambda$. 
The exponents were estimated by the least square method. The intervals on which 
the power laws hold were determined by the method in \cite{Clauset2009} using $\chi^2$ statistic. 
The estimated values of $\beta$ are within the range between 2 and 3.5. However, the values themselves 
may not be important. Here, what is notable is that $C(k)$ are decreasing functions of $k$ on certain ranges. 
This property of $C(k)$ can be associated with modularity of networks \cite{Ravasz2003,Newman2010}, 
that is a common feature of biological networks. 

\begin{figure}
\centering
\begin{tabular}{c}
\includegraphics[width=7cm]{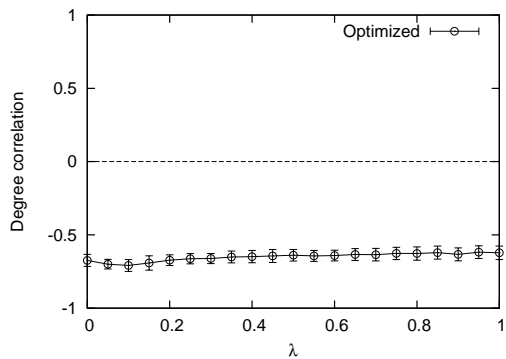} 
\end{tabular}
\caption{
The degree correlation coefficient averaged over the 100 optimized networks for each $\lambda$. 
Bars are standard deviations. 
}
\label{fig10}
\end{figure}

Finally, we analyze degree correlation \cite{Newman2002} in the optimized networks. 
The degree correlation coefficient $r$ for a network is defined by 
\[
r=\frac{\sum_{i,j} (A_{ij} - k_i k_j /2a)k_i k_j}{\sum_{i,j} (k_i \delta_{ij} - k_i k_j /2a)k_i k_j}, 
\]
where $A$ is the adjacency matrix of the network, $k_i$ is the degree of node $i$, $a$ is the number of edges and 
$\delta_{ij}$ is the Kronecker's delta. 
It is known that the value of $r$ has a certain trend depending on types of networks 
such as social, technological, biological and so on. 
In particular, biological networks typically have negative degree correlations \cite{Newman2002}. 
Figure \ref{fig10} shows the degree correlation coefficient $r$ averaged over the 100 optimized networks 
for each $\lambda$. We can see that $r$ takes negative values for each $\lambda$. 

By combining the result for the relationship between LBC and DBC and that for standard complex networks measures, 
we conclude that when the lateral efficiency term is dominant ($\lambda<0.2$) in the quality function $Q$, 
networks that have similar qualitative features found typically in biological networks emerge through the 
optimization process. This suggests that the dynamic mode plays a significant role in the evolution of 
biological networks. 

\section{Conclusions and Future Work}
Let us summarize the results achieved in this paper. 

In Introduction, we introduced the two modes of directed networks inspired by biological networks. 
In the static mode, a network is a passage on which something flows. Thus, the network itself 
is regarded as a static structure. On the other hand, in the dynamic mode, a network is seen as 
a pattern constructed by gluing functions of entities constituting the network. In the dynamic mode, 
interaction between entities is seen as interface between functions. 

In Section 2, first we formalized the dynamic mode in the framework of category theory. Then, 
a canonical representation of the dynamic mode was deduced in terms of the category theoretic universality 
with respect to the idea ``interaction as interface between functions''. A path notion called 
lateral path is naturally associated with it. The directed path, the path notion corresponding to the static mode, 
is obtained from the dual (right adjoint) functor to the functor induced by the canonical representation. 
The results in Section 2 are generalized in Appendix. 

In Section 3, we examined how the dynamic and static modes are embedded in real biological networks by using 
centralities of arcs called lateral betweenness centrality (LBC) and directed betweenness centrality (DBC). 
DBC is a measure of importance of arcs with respect to the static mode or transport on a network. 
On the other hand, LBC is a measure of importance of arcs with respect to the dynamic mode or coherence of a network. We found that in several types of real biological networks, 
there is a trade-off relationship between LBC and DBC, namely, if one has large value, then the value of 
the other is small, that indicates the existence of ``division of labor'' with respect to the two modes. 

In Section 4, we proposed an evolutionary model of complex networks to see how the trade-off relationship 
between LBC and DBC can emerge. Evolution of networks is based on an optimization process toward a quality function. 
We introduced the quality function as the linear combination of the lateral efficiency and the directed efficiency. 
These quantities are indices of intensities of the dynamic and static mode in a network, respectively. 
We found that the trade-off relationship between LBC and DBC emerged only when the lateral efficiency was 
dominant in the quality function. We also found that the optimized networks have features qualitatively similar 
to typical real biological networks. It is suggested that the dynamic mode is an important factor 
in the evolution of biological networks. 

Finally, we indicate two applications of the lateral path in complex directed networks 
to be investigated in future work. 
One is on robustness of networks and the other is on community detection. 
For the configuration model of directed networks (a model of random graph with a specified degree sequence), 
we can compute the size of giant component with respect to the lateral connectedness by applying 
generating function formalism \cite{Newman2001}. Let us consider a percolation problem when arcs are 
removed randomly with probability $p$. The size of the giant component $S_{GC}$ will decrease as $p$ 
approaches 1. The area $A$ under the graph of $S_{GC}$ as a function of $p$ can be a measure of robustness 
of networks \cite{Schneider2011}. By comparing $A$ for a real world network and one for the configuration model 
with the same degree sequence, we can quantify how network structures other than the degree sequence contribute 
to the robustness of the network and discuss the difference from the cases of strong connectedness or 
weak connectedness. 

For community detection \cite{Fortunato2010}, we can consider the lateral path version of the 
Girvan-Newman algorithm \cite{Girvan2002}. In this algorithm, arcs with large LBC will be regarded as 
sitting between communities. An arc that has the largest value of LBC is removed successively as 
LBC of each arc is recalculated. In this process, if a laterally connected component divides into 
several components, then we interpret that it is divided into small communities. To obtain a 
meaningful division into communities, we need a quality function that measures how good the division is 
to stop the process. A candidate is the density of arcs in the lateral direction which is the 
lateral path version of the link density \cite{Ahn2010}. We can define it by means of 
the stability condition with respect to the standard representation $M_0$ obtained in Section 2. 
For another community detection algorithm, one can define a new modularity by considering random walk 
on the set of arcs \cite{Evans2009} with respect to the lateral path. 
These applications will shed a new right on the meaning of the lateral path which is somewhat abstract 
in this work. 

\appendix
\section{General theory of interface}
In this appendix, we generalize the results of Section 2. 

\subsection{Universality of interface transformation induced by a standard representation}

\begin{Def}
Let $\mathcal{C}$ be a small category and $|\mathcal{C}|$ the set of all objects in $\mathcal{C}$. 
The category $\mathcal{C}$ is called \textbf{well-founded} 
if for any nonempty subset $X \subseteq |\mathcal{C}|$ 
there exists a minimal element in $X$ with respect to the binary relation $<$ on $|\mathcal{C}|$ 
defined by $c < c' :\Tot \textrm{ there exists a morphism } f : c \to c' \textrm{ such that }f\textrm{ is not an identity}$, 
where $x \in X$ is called minimal if $y \not< x$ for all $y \in X$. 
\label{def6}
\end{Def}

\begin{Def}
We say that a small category $\mathcal{C}$ has a \textbf{dual structure} 
if there exists an isomorphism $\sigma: \xymatrix{ \mathcal{C} \ar[r]^-{\cong} & \mathcal{C}^{\rm op} }$. 
\label{def7}
\end{Def}

The isomorphism $\sigma$ in Definition \ref{def7} gives rise to an isomorphism that goes backward: 
$\xymatrix{ \mathcal{C}^{\rm op} \ar[r]^-{\cong} & \mathcal{C} }$. However, note that 
$\sigma$ is not necessarily the inverse of itself in general. 

\begin{Def}
Let $\mathcal{C}$ be a small category, $\mathcal{D}$ a bicomplete category and 
$M:\mathcal{C} \to \mathcal{D}$ a representation. 
A functor ${\rm coInt}_M : \mathcal{C}^{\rm op} \to \mathcal{D}$ defined below 
is called \textbf{cointerface} of the representation $M$: 
for any object $c \in \mathcal{C}$, we define 
\[
{\rm coInt}_M(c)={\rm colim} \left( 
\xymatrix{
{\rm Elts}({\mathbf y}^o(c)) \ar[r]^-{\pi_c} & {\mathcal C} \ar[r]^-{M} & {\mathcal D}
}
\right), 
\]
where $\mathbf{y}^o:\mathcal{C}^{\rm op} \to \mathbf{Set}^{\mathcal C}$ is the Yoneda embedding functor 
defined by $\mathbf{y}^o(c)=\mathcal{C}(c,-)$ and the category of elements ${\rm Elts}(F)$ for any functor 
$F \in \mathbf{Set}^{\mathcal C}$ is defined as the category whose objects are pairs 
$(c,x)$ such that $c \in \mathcal{C}$ and $x \in F(c)$ and morphisms $f:(c,x) \to (c',x')$ are morphisms $f:c \to c'$ 
in $\mathcal{C}$ such that $F(f)(x)=x'$. 

For a morphism $f:c \to c' \in \mathcal{C}$, 
let 
$\{ \mu_{(c'',g)}^{M,c'} :  M \circ \pi_{c'}(c'',g)=M(c'') \to {\rm coInt}_M(c') \}_{(c'',g) \in {\rm Elts}({\mathbf y}^o(c'))}$
be the colimit on $M \circ \pi_{c'}$ defining ${\rm coInt}_M(c')$. 
Since 
$\{ \mu_{(c'',g \circ f)}^{M,c} :  M \circ \pi_{c}(c'',g \circ f)=M(c'') \to {\rm coInt}_M(c) \}_{(c'',g) \in {\rm Elts}({\mathbf y}^o(c'))}$
is a cocone on $M \circ \pi_{c'}$, there exists a unique morphism ${\rm coInt}_M(f) : {\rm coInt}_M(c') \to {\rm coInt}_M(c)$ 
such that $\mu_{(c'',g \circ f)}^{M,c}= {\rm coInt}_M(f) \circ \mu_{(c'',g)}^{M,c'}$ for all 
$(c'',g) \in {\rm Elts}({\mathbf y}^o(c'))$. 
\label{def8}
\end{Def}

\begin{Def}
Let $\mathcal{C}$ be a small category equipped with a dual structure 
$\sigma: \xymatrix{ \mathcal{C} \ar[r]^-{\cong} & \mathcal{C}^{\rm op} }$. 
The representation $M_0 : \mathcal{C} \to \hat{\mathcal{C}}$ defined by the composition 
$M_0:={\rm coInt}_{\mathbf{y}} \circ \sigma$ is called the \textbf{standard representation} of $\mathcal{C}$ 
with respect to the dual structure $\sigma$. 
\label{def9}
\end{Def}

For a small category $\mathcal{C}$ with a dual structure $\sigma$, 
there is a natural transformation $\gamma : \mathbf{y} \circ \sigma \to {\rm Int}_{M_0}$ defined as follows
\footnote{
We regard $\sigma$ in the domain of $\gamma$ as the isomorphism $\xymatrix{ \mathcal{C}^{\rm op} \ar[r]^-{\cong} & \mathcal{C} }$. 
}: 
For an object $c \in \mathcal{C}$, 
let 
$\{\nu_{(c',f)}^{M_0,c} :  {\rm Int}_{M_0}(c) \to M_0 \circ \pi_{c}(c',f)=M_0(c')\}_{(c',f) \in {\rm Elts}({\mathbf y}(c))}$ 
be the limit on $M_0 \circ \pi_c$ defining ${\rm Int}_{M_0}(c)$. 
The family of morphisms 
\[
\{
\xymatrix{
\mathbf{y}(\sigma(c)) \ar[rr]^-{\mathbf{y}(\sigma(f))} && \mathbf{y}(\sigma(c')) \ar[rr]^-{ \mu_{(\sigma(c'),{\rm id}_{\sigma(c')})}^{\mathbf{y},\sigma(c')} } && M_0(c')
}
\}_{(c',f) \in {\rm Elts}({\mathbf y}(c))}
\]
is a cone on $M_0 \circ \pi_c$, where 
$\mu_{(\sigma(c'),{\rm id}_{\sigma(c')})}^{\mathbf{y},\sigma(c')}: 
\mathbf{y}(\sigma(c'))=\mathbf{y} \circ \pi_{\sigma(c')}(\sigma(c'),{\rm id}_{\sigma(c')}) \to {\rm coInt}_{\mathbf{y}}(\sigma(c'))=M_0(c')$ 
is a component of the colimit defining $M_0(c')$. 
Indeed, the following diagram commutes for any morphism $h: (c',f) \to (c'',g)$ in ${\rm Elts}(\mathbf{y}(c))$: 
\[
\xymatrix{
& \mathbf{y}(\sigma(c)) \ar[ld]_-{\mu_{(\sigma(c'),{\rm id}_{\sigma(c')})}^{\mathbf{y},\sigma(c')} \circ \mathbf{y}(\sigma(f))} \ar[rd]^-{\mu_{(\sigma(c''),{\rm id}_{\sigma(c'')})}^{\mathbf{y},\sigma(c'')} \circ \mathbf{y}(\sigma(g))} & \\
M_0(c') \ar[rr]_-{M_0(h)} && M_0(c''). 
}
\]
First, note that $h$ satisfies $g \circ h=f$ which is equivalent to $\sigma(h) \circ \sigma(g)=\sigma(f)$, 
which is also equivalent to $\mathbf{y}(\sigma(h)) \circ \mathbf{y}(\sigma(g))=\mathbf{y}(\sigma(f))$. 
Second, by the defining property of $M_0(h)$, we have 
$M_0(h) \circ \mu_{(\sigma(c'),{\rm id}_{\sigma(c')})}^{\mathbf{y},\sigma(c')} = \mu_{(\sigma(c'),h)}^{\mathbf{y},\sigma(c'')}$. 
Third, since $\sigma(h) : (\sigma(c''),{\rm id}_{\sigma(c'')}) \to (\sigma(c'),\sigma(h))$ is a morphism in 
${\rm Elts}(\mathbf{y}^o(\sigma(c'')))$, we have 
$\mu_{(\sigma(c'),\sigma(h))}^{\mathbf{y},\sigma(c'')} \circ \mathbf{y}(\sigma(h)) = \mu_{(\sigma(c''),{\rm id}_{\sigma(c'')})}^{\mathbf{y},\sigma(c'')}$. 
By combining these three equalities, we obtain 
\begin{align*}
M_0(h) \circ \mu_{(\sigma(c'),{\rm id}_{\sigma(c')})}^{\mathbf{y},\sigma(c')} \circ \mathbf{y}(\sigma(f)) 
&= \mu_{(\sigma(c'),\sigma(h))}^{\mathbf{y},\sigma(c'')} \circ \mathbf{y}(\sigma(h)) \circ \mathbf{y}(\sigma(g)) \\
&= \mu_{(\sigma(c''),{\rm id}_{\sigma(c'')})}^{\mathbf{y},\sigma(c'')} \circ \mathbf{y}(\sigma(g)) 
\end{align*}
as required. 

Consequently, by the universality of the limit ${\rm Int}_{M_0}(c)$, 
there exists a unique morphism $\gamma_c:\mathbf{y}(\sigma(c)) \to {\rm Int}_{M_0}(c)$ such that 
$\nu_{(c',f)}^{M_0,c} \circ \gamma_c=\mu_{(\sigma(c'),{\rm id}_{\sigma(c')})}^{\mathbf{y},\sigma(c')} \circ \mathbf{y}(\sigma(f))$
for all $(c',f) \in {\rm Elts}({\mathbf y}(c))$. 
$\gamma_c$ turns out to be natural in $c$. 
Indeed, we have the following commutative diagram for morphisms 
$\xymatrix{ c'' \ar[r]^-{f'} & c' \ar[r]^-{f} & c}$ in $\mathcal{C}$: 
\[
\xymatrix{
\mathbf{y}(\sigma(c)) \ar[rr]^-{\gamma_c} \ar[dd]_-{\mathbf{y}(\sigma(f))} \ar[rd]^-{\delta} && {\rm Int}_{M_0}(c) \ar[ld]_-{\nu_{(c'',f \circ f')}^{M_0 ,c}} \ar[dd]^-{{\rm Int}_{M_0}(f)} \\
& M_0(c'') & \\
\mathbf{y}(\sigma(c')) \ar[ru]^-{\delta'} \ar[rr]_-{\gamma_{c'}} && {\rm Int}_{M_0}(c'), \ar[lu]_-{\nu_{(c'',f')}^{M_0,c'}} 
}
\]
where we put 
$\delta=\mu_{(\sigma(c''),{\rm id}_{\sigma(c'')})}^{\mathbf{y},\sigma(c'')} \circ \mathbf{y}(\sigma(f') \circ \sigma(f))$ and 
$\delta'=\mu_{(\sigma(c''),{\rm id}_{\sigma(c'')})}^{\mathbf{y},\sigma(c'')} \circ \mathbf{y}(\sigma(f'))$. 
By applying Lemma \ref{lem1}, we obtain the desired result. 

Now, let us obtain a concrete expression of the natural transformation $\gamma : \mathbf{y} \circ \sigma \to {\rm Int}_{M_0}$. 
For any objects $c,c' \in \mathcal{C}$, we take ${\rm Int}_{M_0}(c)(c')$ as the set of tuples 
$\left( a_{(d,f)} \right) \in \prod_{(d,f) \in {\rm Elts}(\mathbf{y}(c))} M_0(d)(c')$ 
such that $h \cdot a_{(d,f)}=a_{(e,g)}$ for any morphism $h:d \to e$ in $\mathcal{C}$ such that $g \circ h =f$, 
where $h \cdot a_{(d,f)}=M_0(h)(c')\left( a_{(d,f)} \right)$. 
Maps $\nu_{(d,f)}^{M_0,c}(c') : {\rm Int}_{M_0}(c)(c') \to \left( M_0 \circ \pi_c(d,f) \right) (c')=M_0(d)(c')$ 
are projections. 

Putting $\overline{c'}=\sigma^{-1}(c')$, we take 
\[
M_0(d)(c')= \left( \sum_{c'' \in \mathcal{C}} \mathbf{y}^o(\sigma(d))(\sigma(c'')) \times \mathbf{y}(\sigma(c''))(\sigma(\overline{c'})) \right)/\sim, 
\]
where $\sim$ is the equivalence relation generated by the relation $r$ defined by 
$(\sigma(i),\sigma(j)) r (\sigma(k) \circ \sigma(i),\sigma(k) \circ \sigma(j))$ 
for morphisms 
\[
\xymatrix@-1pc{
\sigma(d) \ar[rd]^-{\sigma(i)} & & \\
& \sigma(c'') \ar[r]^-{\sigma(k)} & \sigma(c''') \\
\sigma(\overline{c'}) \ar[ru]_-{\sigma(j)} & &
}
\]
in $\mathcal{C}$. For any morphism $\sigma(f'): \sigma(d) \to \sigma(c'')$ in $\mathcal{C}$, 
the map 
\[
\mu_{(\sigma(c''),\sigma(f'))}^{\mathbf{y},\sigma(d)}(\sigma(\overline{c'})) : 
\mathbf{y}(\sigma(c''))(\sigma(\overline{c'})) \to M_0(d)(\sigma(\overline{c'}))
\]
is given by $\sigma(g') \mapsto [(\sigma(f'),\sigma(g'))]$. 

By this expression of $M_0$, we have 
$h \cdot a_{(d,f)}=[(\sigma(i) \circ \sigma(h),\sigma(j))]$ for $a_{(d,f)}=[(\sigma(i),\sigma(j))]$. 

The map $\gamma_c(c'): \mathbf{y}(\sigma(c))(c') \to {\rm Int}_{M_0}(c)(c')$ is uniquely determined by 
the following commutative diagram: 
\[
\xymatrix{
\mathbf{y}(\sigma(c))(c') \ar[r]^-{\gamma_c(c')} \ar[d]_-{\mathbf{y}(\sigma(f))(c')} & {\rm Int}_{M_0}(c)(c') \ar[d]^-{\nu_{(d,f)}^{M_0,c}(c')} \\
\mathbf{y}(\sigma(d))(c') \ar[r]^-{\mu_{(\sigma(d),{\rm id}_{\sigma(d)})}^{\mathbf{y},\sigma(d)}(c')} & M_0(d)(c'). 
}
\]
Hence, for a morphism $\sigma(s):\sigma(\overline{c'}) \to \sigma(c)$ in $\mathcal{C}$, if we put 
$\gamma_c(c')(\sigma(s))=(a_{(d,f)})$, 
then $a_{(d,f)}=[({\rm id}_{\sigma(d)},\sigma(f) \circ \sigma(s))]$. 

\begin{Pro}
Let $\mathcal{C}$ be a small category equipped with a dual structure 
$\sigma: \xymatrix{ \mathcal{C} \ar[r]^-{\cong} & \mathcal{C}^{\rm op} }$. 
For any presheaf $G \in \hat{\mathcal{C}}$ and $c \in \mathcal{C}$, 
if $\gamma_c: \mathbf{y}(\sigma(c)) \to {\rm Int}_{M_0}(c)$ is an isomorphism, 
then we have the following expression for the component of the natural transformation 
$\varphi_c^{M_0,G}:G(c) \to G \otimes M_0(\sigma(c))$: 
for any $x \in G(c)$, $\varphi_c^{M_0,G}(x)=[(x,[({\rm id}_{\sigma(c)},{\rm id}_{\sigma(c)})])]$, 
where we identify $\hat{\mathcal{C}}({\rm Int}_{M_0}(c),G \otimes M_0)$ with $G \otimes M_0(\sigma(c))$ by 
\[
\xymatrix{
\hat{\mathcal{C}}({\rm Int}_{M_0}(c),G \otimes M_0) \ar[r]_-{(-) \circ \gamma_c}^-{\cong} & \hat{\mathcal{C}}(\mathbf{y}(\sigma(c)),G \otimes M_0) 
\ar[r]_-{\theta_c}^-{\cong} & G \otimes M_0 (\sigma(c))
}
\]
and $\theta_c$ is a natural isomorphism given by the Yoneda Lemma. 
\label{pro5}
\end{Pro}
\textit{Proof. }
By definition, $\varphi_c^{M_0,G}:G(c) \to G \otimes M_0(\sigma(c))$ is the map given by 
\begin{align*}
\varphi_c^{M_0,G}(x)
&= \theta_c \left( \mu_{(c,x)}^{M_0,G} \circ \nu_{(c,{\rm id}_c)}^{M_0,c} \circ \gamma_c \right) \\
&= \left( \mu_{(c,x)}^{M_0,G} \circ \nu_{(c,{\rm id}_c)}^{M_0,c} \circ \gamma_c \right) (\sigma(c)) ({\rm id}_{\sigma(c)}) 
\end{align*}
for each $x \in G(c)$. However, the composition 
\[
\xymatrix{
\mathbf{y}(\sigma(c))(\sigma(c)) \ar[r]^-{\gamma_c(\sigma(c))} & {\rm Int}_{M_0}(c)(\sigma(c)) \ar[r]^-{\nu_{(c,{\rm id}_c)}^{M_0,c}} & M_0(c)(\sigma(c)) \ar[r]^-{\mu_{(c,x)}^{M_0,G}} & G \otimes M_0(\sigma(c)) 
}
\]
gives 
\[
\xymatrix@-1pc{
{\rm id}_{\sigma(c)} \ar@{|->}[r] & \left( [({\rm id}_{\sigma(d)},\sigma(f))] \right)_{(d,f) \in {\rm Elts}(\mathbf{y}(c))} \ar@{|->}[r] & [({\rm id}_{\sigma(c)},\sigma({\rm id}_c))] \ar@{|->}[r] & [(x,[({\rm id}_{\sigma(c)},{\rm id}_{\sigma(c)})])].
}
\]

\hfill $\Box$ \\

\begin{Thm}
Let $\mathcal{C}$ be a well-founded small category equipped with a dual structure 
$\sigma: \xymatrix{ \mathcal{C} \ar[r]^-{\cong} & \mathcal{C}^{\rm op} }$ and 
$M_0 : \mathcal{C} \to \hat{\mathcal{C}}$ the standard representation with respect to $\sigma$. 
If the natural transformation $\gamma : \mathbf{y} \circ \sigma \to {\rm Int}_{M_0}$ is an natural isomorphism, 
then $(\hat{\mathcal{C}}({\rm Int}_{M_0}(-),G \otimes M_0),\varphi^{M_0,G})$ is an initial object of $\mathcal{IT}_G$ 
for any presheaf $G \in \hat{\mathcal{C}}$. 
\label{thm3}
\end{Thm}
\textit{Proof. }
Let $\psi:G \to H$ be an object in $\mathcal{IT}_G$. We shall show that 
there exists a unique morphism $\iota:(\hat{\mathcal{C}}({\rm Int}_{M_0}(-),G \otimes M_0),\varphi^{M_0,G}) \to (H,\psi)$ 
such that $\psi=\iota \circ \varphi^{M_0,G}$. 

We appeal to a well-founded induction. 
Recall that a well-founded induction on a set $X$ with a well-founded binary relation $<$ 
(a binary relation $<$ on $X$ is called well-founded if any nonempty subset of $X$ has 
a minimal element with respect to $<$) asserts that for a condition $P(x)$ on $X$ 
\[
\forall x \in X \left( \left( \forall y \in X \left( y < x \To P(y) \right) \right) \To P(x) \right)
\]
implies $P(x)$ holds for any $x \in X$. 

We apply a well-founded induction on the set $|\mathcal{C}|$ with the 
binary relation $<$ on it defined in Definition \ref{def6} to the following condition $P(c)$: 
$P(c')$ is true for any $c'<c$ and 
there exists a unique map $\iota_c : G \otimes M_0(\sigma(c)) \to H(c)$ such that 
$\iota_c \circ G \otimes M_0(\sigma(f))=H(f) \circ \iota_{c'}$ and 
$\psi_c=\iota_c \circ \varphi_c^{M_0,G}$ for any morphism $f:c' \to c$, 
namely, the following two diagrams commute for any morphism $f:c' \to c$: 
\[
\xymatrix{
G \otimes M_0(\sigma(c)) \ar[r]^-{\iota_c} & H(c) \\
G \otimes M_0(\sigma(c')) \ar[u]^-{G \otimes M_0(\sigma(f))} \ar[r]^-{\iota_{c'}} & H(c'). \ar[u]_-{H(f)}
}
\]
and 
\[
\xymatrix@-1pc{
G \otimes M_0(\sigma(c)) \ar[rr]^-{\iota_c} && H(c) \\
& G(c) \ar[ul]^-{\varphi_c^{M_0,G}} \ar[ur]_-{\psi_c}.
}
\]

If $c \in |\mathcal{C}|$ is a minimal element, then $\sigma(c)$ is maximal and 
thus any morphism $\sigma(f):\sigma(c) \to \sigma(c')$ must be ${\rm id}_{\sigma(c)}$. 
Consequently, we have 
\begin{align*}
G \otimes M_0(\sigma(c)) 
&= \left( \sum_{c' \in \mathcal{C}} G(c') \times M_0(c')(\sigma(c)) \right) / \sim \\
&\cong \left( \sum_{c' \in \mathcal{C}} G(c') \times \mathbf{y}(\sigma(c'))(\sigma(c)) \right) / \sim \\
&= G(c) \times \{ {\rm id}_{\sigma(c)} \} \cong G(c). 
\end{align*}
By the concrete expression of $\varphi_c^{M_0,G}$ given in Proposition \ref{pro5}, 
one can see that $\varphi_c^{M_0,G}$ is an isomorphism. Hence, we obtain a unique map 
that satisfies $P(c)$ by $\iota_c:=\psi_c \circ (\varphi_c^{M_0,G})^{-1}$. 

Now, let $c \in |\mathcal{C}|$ is not minimal and assume that $P(c')$ is true 
for any $c'<c$. We would like to prove that $P(c)$ is also true. 

For any $[(x,\omega)] \in G \otimes M_0(\sigma(c))$, if $\omega=[(\sigma(g),\sigma(f))]$ 
with 
\[
\xymatrix@-1pc{
\sigma(c'') \ar[rd]^-{\sigma(g)} & \\
& \sigma(c') \\
\sigma(c) \ar[ru]_-{\sigma(f)} &
}
\]
in $\mathcal{C}$, then for $\omega':=[(\sigma(g),{\rm id}_{\sigma(c')})] \in M_0(c'')(\sigma(c'))$ 
we have 
\[
G \otimes M_0(\sigma(f)) ([(x,\omega')])=[(x,\omega)], 
\]
where $G \otimes M_0 (\sigma(f))$ sends $[(x,[(\sigma(i),\sigma(j))])] \in G \otimes M_0(\sigma(c'))$ to 
$[(x,[(\sigma(i),\sigma(j) \circ \sigma(f))])] \in G \otimes M_0(\sigma(c))$ in general. 
\[
\xymatrix@-1pc{
& \sigma(c'') \ar[rd]^-{\sigma(i)} & \\
& & \sigma(c''') \\
\sigma(c) \ar[r]_-{\sigma(f)} & \sigma(c') \ar[ru]_-{\sigma(j)} &
}
\]
Note that we can always take $f$ as a non-identity because $c$ is not minimal by the assumption of 
the well-founded induction and hence we can take ${\rm dom}(f)=c'<c$. 

Consequently, we define the map $\iota_c : G \otimes M_0(\sigma(c)) \to H(c)$ as follows: for each 
$[(x,\omega)] \in G \otimes M_0(\sigma(c))$, we take $c'$, $f$ and $\omega'$ as above and define 
\[
\iota_c([(x,\omega)])=H(f) \left( \iota_{c'}([(x,\omega')]) \right). 
\]

Now, let us check $\iota_c$ is a well-defined map. To do so, we take two steps. 
First, we show that $\iota_c([(x,\omega)])$ is independent of the choice of representative of 
$[(x,\omega)]$ and second it does not depend on the choice of representative of $\omega$. 

Let $[(x,\omega)]=[(y,\omega')] \in G \otimes M_0(\sigma(c))$. 
It is sufficient to consider the case when 
$y=x \cdot f$ and $\omega=f \cdot \omega'$ for some morphism $f:c_2 \to c_1$ in $\mathcal{C}$, 
where $x \in G(c_1)$, $y \in G(c_2)$, $\omega \in M_0(c_1)(\sigma(c))$ and 
$\omega' \in M_0(c_2)(\sigma(c))$. 

Let $(\sigma(i),\sigma(j)) \in \mathbf{y}^o(\sigma(c_2))(\sigma(c')) \times \mathbf{y}(\sigma(c'))(\sigma(c))$ 
be a representative of $\omega'$ such that $j$ is a non-identity. 
\[
\xymatrix@-1pc{
\sigma(c_1) \ar[r]^-{\sigma(f)} & \sigma(c_2) \ar[rd]^-{\sigma(i)} & \\
& & \sigma(c') \\
& \sigma(c) \ar[ru]_-{\sigma(j)} & 
}
\]
Since $\omega=f \cdot \omega'=f \cdot [(\sigma(i),\sigma(j))]=[(\sigma(i) \circ \sigma(f),\sigma(j))]$, 
$(\sigma(i) \circ \sigma(f),\sigma(j))$ is a representative of $\omega$. 
Putting $\chi:=[(\sigma(i),{\rm id}_{\sigma(c')})]$, 
we have $[(y,\chi)]=[(x \cdot f,\chi)]=[(x, f \cdot \chi)]$. 
This completes the first step. 

For the second step, let $[(x,\omega)] \in G \otimes M_0(\sigma(c))$ and 
$(x,\omega) \in G(c') \times M_0(c')(\sigma(c))$. 
It is enough to show that if $(\sigma(i),\sigma(j)) \in \omega$ and $j$ is a non-identity, then 
\[
H(j \circ k) \circ \iota_{c'''}([(x,[(\sigma(k) \circ \sigma(i),{\rm id}_{\sigma(c''')})])])
=H(j) \circ \iota_{c''}([(x,[(\sigma(i),{\rm id}_{\sigma(c'')})])])
\]
for any morphism $k$ in $\mathcal{C}$ such that 
\[
\xymatrix@-1pc{
\sigma(c') \ar[rd]^-{\sigma(i)} & & \\
& \sigma(c'') \ar[r]^-{\sigma(k)} & \sigma(c'''). \\
\sigma(c) \ar[ru]_-{\sigma(j)} & &
}
\]
However, by the assumption of the well-founded induction, we have 
\begin{align*}
&\quad H(j \circ k) \circ \iota_{c'''}([(x,[(\sigma(k) \circ \sigma(i),{\rm id}_{\sigma(c''')})])]) \\
&= H(j) \circ H(k) \circ \iota_{c'''}([(x,[(\sigma(k) \circ \sigma(i),{\rm id}_{\sigma(c''')})])]) \\
&= H(j) \circ \iota_{c''} \circ G \otimes M_0(\sigma(k))([(x,[(\sigma(k) \circ \sigma(i),{\rm id}_{\sigma(c''')})])]) \\
&= H(j) \circ \iota_{c''}([(x,[(\sigma(k) \circ \sigma(i),\sigma(k))])]) \\
&= H(j) \circ \iota_{c''}([(x,[(\sigma(i),{\rm id}_{\sigma(c'')})])]). 
\end{align*}

By the construction of $\iota_c$, it is clear that we have 
$\iota_c \circ G \otimes M_0(\sigma(f))=H(f) \circ \iota_{c'}$ for any morphism $f:c' \to c$ 
and $\iota_c$ is a unique map satisfying this condition. 

Finally, let us show that $\psi_c=\iota_c \circ \varphi_c^{M_0,G}$. 
For any $x \in G(c)$, we have 
$\varphi_c^{M_0,G}(x)=[(x,[({\rm id}_{\sigma(c)},{\rm id}_{\sigma(c)})])]$ by Proposition \ref{pro5}. 
Since $c$ is not minimal, we can take a morphism $f:c' \to c$ in $\mathcal{C}$ such that $c' \neq c$. 
Consequently, we have 
\begin{align*}
\iota_c \circ \varphi_c^{M_0,G}(x)
&= \iota_c([(x,[({\rm id}_{\sigma(c)},{\rm id}_{\sigma(c)})])]) \\
&= \iota_c([(x,[(\sigma(f),\sigma(f))])]) \\
&= \iota_c \circ G \otimes M_0(\sigma(f))([(x,[(\sigma(f),{\rm id}_{\sigma(c')})])]) \\
&= H(f) \circ \iota_{c'}([(x \cdot f,[({\rm id}_{\sigma(c')},{\rm id}_{\sigma(c')})])]) \\
&= H(f) \circ \iota_{c'} \circ \varphi_{c'}^{M_0,G}(x \cdot f) \\
&= H(f) \circ \psi_{c'}(x \cdot f) \\
&= \psi_c(x), 
\end{align*}
where we use the naturality $\iota$ in the fourth equality, Proposition \ref{pro5} in the 
fifth equality, the assumption of the well-founded induction in the sixth equality 
and the naturality of $\psi$ in the last equality. 

\hfill $\Box$ \\

One can check that the category $\upuparrows$ satisfies the assumption of Theorem \ref{thm3}. 
Hence, Theorem \ref{thm3} is a generalization of Theorem \ref{thm1}. However, the condition 
that the natural transformation $\gamma : \mathbf{y} \circ \sigma \to {\rm Int}_{M_0}$ is a 
natural isomorphism seems to be too abstract. In the following, we shall translate this condition 
into a condition on the structure of the small category $\mathcal{C}$. However, translation in 
general situation is not so informative. Here, we focus on the case $\mathcal{C}$ is a category 
freely generated by a finite acyclic directed graph $\Gamma$ with a directed graph isomorphism 
$\sigma:\xymatrix{ \Gamma \ar[r]^-{\cong} & \Gamma^{\rm op}}$, where $\Gamma^{\rm op}$ is 
a directed graph obtained by reversing the direction of all the arcs in $\Gamma$. $\sigma$ 
induces an isomorphism of categories $\xymatrix{ \mathcal{C} \ar[r]^-{\cong} & \mathcal{C}^{\rm op}}$ 
which is also denoted by $\sigma$. Thus, $\mathcal{C}$ is a well-founded category with 
a dual structure $\sigma$. 

\begin{Def}
Let $\mathcal{C}$ be a category. A \textbf{span} in $\mathcal{C}$ is a pair of morphisms 
$(i,j)$ in $\mathcal{C}$ with the common domain. 
\[
\xymatrix@-1pc{
& c' \\
c \ar[ru]^-{i} \ar[rd]_-{j} & \\
& c''
}
\]

A span $(i,j)$ is called \textbf{maximal} if there is no 
non-identity morphism $k$ and span $(i',j')$ such that 
$i=i' \circ k$ and $j=j' \circ k$, namely, if we cannot find 
morphisms $k,i',j'$ with $k$ non-identity that make the following diagram commute: 
\[
\xymatrix{
& & c' \\
c \ar@/^/[rru]^-{i} \ar@/_/[rrd]_-{j} \ar[r]^-{k} & c''' \ar[ru]_-{i'} \ar[rd]^-{j'} & \\
& & c''
}
\]

A \textbf{maximal span for a span $(i,j)$} is a maximal span $(\tilde{i},\tilde{j})$ 
such that there exists a morphism $k$ such that $i=\tilde{i} \circ k$ and $j=\tilde{j} \circ k$. 
\label{def10}
\end{Def}

Let us consider the following condition (A) for a category $\mathcal{C}$: 
\begin{description}
\item[(A)]
for any span $(i,j)$, if there is a maximal span for it, then it is unique. 
\end{description}

\begin{Lem}
Let $\mathcal{C}$ be a well-founded small category equipped with a dual structure 
$\sigma: \xymatrix{ \mathcal{C} \ar[r]^-{\cong} & \mathcal{C}^{\rm op} }$. 
For any span $(i,j)$ there is a maximal span for it. 
\label{lem2}
\end{Lem}
\textit{Proof. }
Let $(i,j)$ be a span and put 
\[
X=\{ {\rm cod}(k) | \exists (i',j') \textrm{: span such that } i=i' \circ k \textrm{ and } j=j' \circ k \}. 
\]
Since $c:={\rm dom}(i)={\rm dom}(j) \in X$ with $i'=i$, $j'=j$ and $k={\rm id}_c$, 
we have $X \neq \emptyset$. 
$X$ has a maximal element $c_m$ with respect to the binary relation $<$ defined in Definition \ref{def6} 
because $\mathcal{C}$ is well-founded and $\mathcal{C} \cong \mathcal{C}^{\rm op}$. 
Consequently, there exist a morphism $k_m:c \to c_m$ and span $(i_m,j_m)$ with 
${\rm dom}(i_m)={\rm dom}(j_m)=c_m$ in $\mathcal{C}$ such that $i=i_m \circ k_m$ and $j=j_m \circ k_m$. 

Suppose that $(i_m,j_m)$ is not a maximal span. Then, there exist a non-identity morphism $k$ 
and a span $(i',j')$ such that $i_m=i' \circ k$ and $j_m=j' \circ k$. 
We have $i=i_m \circ k_m=i' \circ (k \circ k_m)$ and $j=j_m \circ k_m=j' \circ (k \circ k_m)$. 
Hence, ${\rm cod}(k \circ k_m) \in X$. 
On the other hand, we have $c_m < {\rm cod}(k)={\rm cod}(k \circ k_m)$ because $k$ is not an identity. 
However, this contradicts to the maximality of $c_m$ in $X$. 
Thus, $(i_m,j_m)$ is a maximal span for $(i,j)$. 
\[
\xymatrix{
& & c' \\
c \ar@/^/[rru]^-{i} \ar[r]^-{k_m} \ar@/_/[rrd]_-{j} & c_m \ar[ru]_-{i_m} \ar[r]^-{k} \ar[rd]^-{j_m} & c''' \ar[u]_-{i'} \ar[d]^-{j'} \\
& & c''
}
\]

\hfill $\Box$ \\

\begin{Lem}
Let $\mathcal{C}$ be a well-founded small category equipped with a dual structure 
$\sigma: \xymatrix{ \mathcal{C} \ar[r]^-{\cong} & \mathcal{C}^{\rm op} }$. 
Assume that the condition (A) holds for $\mathcal{C}$. 
If maximal spans $(i,j)$ and $(i',j')$ satisfies 
$[(\sigma(i),\sigma(j))]=[(\sigma(i'),\sigma(j'))] \in M_0(c)(\sigma(c'))$ 
where $c:={\rm cod}(i)={\rm cod}(i')$ and $c':={\rm cod}(j)={\rm cod}(j')$, 
then we have $(i,j)=(i',j')$. 
\label{lem3}
\end{Lem}
\textit{Proof. }
If $[(\sigma(i),\sigma(j))]=[(\sigma(i'),\sigma(j'))]$, then 
there exist a sequence of spans 
\[
(i,j)=(i_0,j_0), (i_1,j_1),\cdots,(i_n,j_n)=(i',j')
\]
and a sequence of morphisms $k_1,k_2,\cdots,k_n$ such that 
for each $p=1,2,\cdots,n$ either of the following two conditions 
holds: 
(1) $i_p=i_{p-1} \circ k_p$ and $j_p=j_{p-1} \circ k_p$ 
or (2) $i_{p-1}=i_{p} \circ k_p$ and $j_{p-1}=j_{p} \circ k_p$. 
Since $(i_0,j_0)$ and $(i_n,j_n)$ are maximal spans, 
we have (1) for $p=1$ and (2) for $p=n$. We can assume that 
(1) and (2) occur alternately without loss of generality. 
Thus, $n$ is an even number. 

We can also assume that $(i_p,j_p)$ is maximal if $p$ is an even number. 
Indeed, for example, if $(\tilde{i}_2,\tilde{j}_2)$ is a maximal span for $(i_2,j_2)$ 
with a morphism $\tilde{k}_2$ such that $i_2=\tilde{i}_2 \circ \tilde{k}_2$ 
and $j_2=\tilde{j}_2 \circ \tilde{k}_2$, then we have 
$i_1=\tilde{i}_2 \circ (\tilde{k}_2 \circ k_2)$ and $j_1=\tilde{j}_2 \circ (\tilde{k}_2 \circ k_2)$ 
for $p=2$ and 
$i_3=\tilde{i}_2 \circ (\tilde{k}_2 \circ k_3)$ and $j_3=\tilde{j}_2 \circ (\tilde{k}_2 \circ k_3)$ 
for $p=3$. 

Now, both $(i_0,j_0)$ and $(i_2,j_2)$ are maximal spans for $(i_1,j_1)$. 
It follows that $(i_0,j_0)=(i_2,j_2)$ by (A). 
By the same manner, we have 
$(i_0,j_0)=(i_2,j_2)=(i_4,j_4)= \cdots =(i_n,j_n)$. 

\hfill $\Box$ \\

\begin{Pro}
Let $\mathcal{C}$ be a well-founded small category equipped with a dual structure 
$\sigma: \xymatrix{ \mathcal{C} \ar[r]^-{\cong} & \mathcal{C}^{\rm op} }$. 
If the condition (A) holds for $\mathcal{C}$, then the map 
$\gamma_c(\sigma(c')):\mathbf{y}(\sigma(c))(\sigma(c')) \to {\rm Int}_{M_0}(c)(\sigma(c'))$ 
is injective for any $c,c' \in \mathcal{C}$. 
\label{pro6}
\end{Pro}
\textit{Proof. }
Assume that $\gamma_c(\sigma(c'))(\sigma(s))=\gamma_c(\sigma(c'))(\sigma(s'))$ for 
morphisms $s,s':c \to c'$ in $\mathcal{C}$. Then, we have 
\[
[({\rm id}_{\sigma(c'')},\sigma(f) \circ \sigma(s))]=[({\rm id}_{\sigma(c'')},\sigma(f) \circ \sigma(s'))]
\]
for all $(c'',f) \in {\rm Elts}(\mathbf{y}(c))$. 
Both $({\rm id}_{\sigma(c'')},\sigma(f) \circ \sigma(s))$ and 
$({\rm id}_{\sigma(c'')},\sigma(f) \circ \sigma(s'))$ are maximal spans because of 
the well-foundedness of $\mathcal{C}$. By Lemma \ref{lem3}, we have 
$\sigma(f) \circ \sigma(s)=\sigma(f) \circ \sigma(s')$, namely, $s \circ f=s' \circ f$ 
for each morphism $f:c'' \to c$. If we take $f={\rm id}_c$, then we have $s=s'$. 

\hfill $\Box$ \\

\begin{Pro}
Let $\mathcal{C}$ be a well-founded small category equipped with a dual structure 
$\sigma: \xymatrix{ \mathcal{C} \ar[r]^-{\cong} & \mathcal{C}^{\rm op} }$. 
Assume that the condition (A) holds for $\mathcal{C}$. 
For any $c,c' \in \mathcal{C}$, the map 
$\gamma_c(\sigma(c')):\mathbf{y}(\sigma(c))(\sigma(c')) \to {\rm Int}_{M_0}(c)(\sigma(c'))$ 
is surjective if and only if the following condition (B) holds: 
\begin{description}
\item[(B)]
For any family of maximal spans $\{(i_f,j_f)\}_{(c'',f) \in {\rm Elts}(\mathbf{y}(c))}$ such that 
${\rm cod}(i_f)=c''$ and ${\rm cod}(j_f)=c'$ for each $(c'',f) \in {\rm Elts}(\mathbf{y}(c))$, 
assume that if $f=g \circ h$, then 
$(i_g,j_g)$ is a unique maximal span for the span $(h \circ i_f,j_f)$. 
Then, $i_f={\rm id}_{c''}$ for any $(c'',f) \in {\rm Elts}(\mathbf{y}(c))$ and 
there exists a morphism $s : c \to c'$ such that $j_f=s \circ f$ for any $(c'',f) \in {\rm Elts}(\mathbf{y}(c))$. 
\[
\xymatrix{
c''' \ar[r]^-{i_f} \ar[d]_-{j_f} & c'' \ar[d]^-{f} \ar[rd]^-{h} & \\
c' & c \ar[l]^-{s} & d \ar[l]^-{g}
}
\]
\end{description}
\label{pro7}
\end{Pro}
\textit{Proof. }
Suppose that the condition (B) holds for $c,c' \in \mathcal{C}$. 
For an element $\left( [(\sigma(i_f),\sigma(j_f))] \right)_{(c'',f) \in {\rm Elts}(\mathbf{y}(c))}$, 
we can assume that each representative $(i_f,j_f)$ of $[(\sigma(i_f),\sigma(j_f))]$ is a maximal span 
by taking its unique maximal span if necessary because a span and its maximal span are contained in the 
same equivalence class. For a morphism $f:c'' \to c$, if 
$f=g \circ h$, then 
$[(\sigma(i_f) \circ \sigma(h),\sigma(j_f))]=h \cdot [(\sigma(i_f),\sigma(j_f))]=[(\sigma(i_g),\sigma(j_g))]$. 
If $(\tilde{i}_f,\tilde{j}_f)$ is a unique maximal span for the span $(h \circ i_f,j_f)$, 
then we have $\tilde{i}_f=i_g$ and $\tilde{j}_f=j_g$ by Lemma \ref{lem3}. 
By the condition (B), $i_f={\rm id}_{c''}$ for any morphism $f:c'' \to c$ and 
there exists $s:c \to c'$ such that $j_f=s \circ f$ for any morphism $f:c'' \to c$. 
Consequently, we have $[(\sigma(i_f),\sigma(j_f))]=[({\rm id}_{\sigma(c'')},\sigma(f) \circ \sigma(s))]$ 
for any $(c'',f) \in {\rm Elts}(\mathbf{y}(c))$ which means that 
$\gamma_c(\sigma(c'))$ is surjective.

For the reverse direction, suppose that $\gamma_c(\sigma(c'))$ is surjective. 
Take any family of maximal spans $\{(i_f,j_f)\}_{(c'',f) \in {\rm Elts}(\mathbf{y}(c))}$ such that 
${\rm cod}(i_f)=c''$ and ${\rm cod}(j_f)=c'$ for each $(c'',f) \in {\rm Elts}(\mathbf{y}(c))$ and 
assume that if $f=g \circ h$, then 
$(i_g,j_g)$ is a unique maximal span for the span $(h \circ i_f,j_f)$. 
Clearly, $p:=\left( [(\sigma(i_f),\sigma(j_f))] \right)_{(c'',f) \in {\rm Elts}(\mathbf{y(c)})}$ 
is an element of ${\rm Int}_{M_0}(c)(\sigma(c'))$. 
Since $\gamma_c(\sigma(c'))$ is a surjection, there exists a morphism $s:c \to c'$ such that 
$\gamma_c(\sigma(c'))(\sigma(s))=p$. 
Then, we have 
$[(\sigma(i_f),\sigma(j_f))]=[({\rm id}_{\sigma(c'')},\sigma(f) \circ \sigma(s))]$ 
for any $(c'',f) \in {\rm Elts}(\mathbf{y(c)})$. By Lemma \ref{lem3}, 
we obtain $i_f={\rm id}_{c''}$ and $j_f=s \circ f$ for any 
$(c'',f) \in {\rm Elts}(\mathbf{y}(c))$. 

\hfill $\Box$ \\

\begin{Lem}
If a category $\mathcal{C}$ is freely generated by a finite acyclic directed graph, 
then (A) holds for $\mathcal{C}$. 
\label{lem4}
\end{Lem}
\textit{Proof. }
For a span $(i,j)$, let $(\tilde{i}_1,\tilde{j}_1)$ and $(\tilde{i}_2,\tilde{j}_2)$ be 
its two maximal spans. Then, there exist morphisms $k_1$ and $k_2$ such that 
$i=\tilde{i}_1 \circ k_1=\tilde{i}_2 \circ k_2$ and $j=\tilde{j}_1 \circ k_1=\tilde{j}_2 \circ k_2$. 
Assume that $i$ can be written as 
$i=f_n \circ f_{n-1} \circ \cdots f_1$, where $f_l, \ l=1,2,\cdots,n$ are generating morphisms for $\mathcal{C}$, namely, 
morphisms corresponding to arcs in the base graph. 
If $(i,j)$ is maximal, then we have nothing to do. Hence, we assume that 
$(i,j)$ is not a maximal span. Consequently, morphisms $i$, $k_1$ and $k_2$ are 
not identities. 

We can write $k_1=f_{n_1} \circ \cdots \circ f_1$ and $k_2=f_{n_2} \circ \cdots \circ f_1$. 
Suppose that $n_1 \neq n_2$. 
Without loss of generality, we can assume that $n_1>n_2$. 
Putting $g=f_{n_1} \circ \cdots \circ f_{n_2+1}$, we have $k_1=g \circ k_2$. 
It follows that $\tilde{i}_2=\tilde{i}_1 \circ g$ and $\tilde{j}_2=\tilde{j}_1 \circ g$. 
Since $g$ is not an identity, this contradicts to the maximality of the span $(\tilde{i}_2,\tilde{j}_2)$. 
Consequently, we must have $k_1=k_2$ which implies that 
$(\tilde{i}_1,\tilde{j}_1)=(\tilde{i}_2,\tilde{j}_2)$. 

\hfill $\Box$ \\

\begin{Lem}
Let $\mathcal{C}$ be a category 
freely generated by a finite acyclic directed graph $\Gamma$ with a directed graph isomorphism 
$\sigma:\xymatrix{ \Gamma \ar[r]^-{\cong} & \Gamma^{\rm op}}$. 
For any $c,c' \in \mathcal{C}$, the condition (B) is equivalent to the following condition (B'): 
\begin{description}
\item[(B')]
For any family of maximal spans $\{(i_f,j_f)\}_{(c'',f) \in {\rm Elts}(\mathbf{y}(c))}$ such that 
${\rm cod}(i_f)=c''$ and ${\rm cod}(j_f)=c'$ for each $(c'',f) \in {\rm Elts}(\mathbf{y}(c))$, 
assume that if $f=g \circ h$, then 
$(i_g,j_g)$ is a unique maximal span for the span $(h \circ i_f,j_f)$. 
Then, $i_{{\rm id}_c}={\rm id}_{c}$. 
\end{description}
\label{lem5}
\end{Lem}
\textit{Proof. }
Assume that (B') holds for $c,c' \in \mathcal{C}$. 
Take any family of maximal spans $\{(i_f,j_f)\}_{(c'',f) \in {\rm Elts}(\mathbf{y}(c))}$ such that 
${\rm cod}(i_f)=c''$ and ${\rm cod}(j_f)=c'$ for each $(c'',f) \in {\rm Elts}(\mathbf{y}(c))$ and 
assume that if $f=g \circ h$, then 
$(i_g,j_g)$ is a unique maximal span for the span $(h \circ i_f,j_f)$. 
Then, for any morphism $f:c'' \to c$, $(i_{{\rm id}_c},j_{{\rm id}_c})$ is a maximal span for $(f \circ i_f,j_f)$ 
and there exists a morphism $k_f$ such that $f \circ i_f=i_{{\rm id}_c} \circ k_f$ and $j_f=j_{{\rm id}_c} \circ k_f$ 
because $f={\rm id}_c \circ f$. 
\[
\xymatrix{
& c'' \ar[r]^-{f} & c \\
d \ar[ru]^-{i_f} \ar[r]^-{k_f} \ar[rd]_-{j_f} & d' \ar[ru]_-{i_{{\rm id}_c}} \ar[d]^-{j_{{\rm id}_c}} & \\
& c' &
}
\]

Now, suppose that $i_f \neq {\rm id}_{c''}$. 
Let us write $i_f=a_n \circ \cdots \circ a_1$ and $j_f=b_m \circ \cdots \circ b_1$, where 
$a_p$ and $b_q$ are generating morphisms for $p=1,\cdots,n$ and $q=1,\cdots,m$. 
Then, we have $n \geq 1$. If $m=0$, then $j_f={\rm id}_{c'}$, which in turn implies that 
$k_f$ is an identity. If $m \geq 1$, then we must have $a_1 \neq b_1$ because 
$(i_f,j_f)$ is a maximal span. If $k_f$ is not an identity, then we can write 
$k_f=c_l \circ \cdots \circ c_1$ with $l \geq 1$, where $c_1, \cdots, c_l$ are generating morphisms. However, by the equalities 
$f \circ i_f=i_{{\rm id}_c} \circ k_f$ and $j_f=j_{{\rm id}_c} \circ k_f$, it follows 
that $a_1=c_1=b_1$. This is a contradiction. Therefore, $k_f$ is an identity also in this case. 
On the other hand, we have $i_{{\rm id}_c}={\rm id}_c$ by the condition (B'). 
Thus, we obtain $f \circ i_f={\rm id}_c$. However, this contradicts to the well-foundedness of 
$\mathcal{C}$. It follows that $i_f = {\rm id}_{c''}$, $f=k_f$ and $j_f=j_{{\rm id}_c} \circ f$. 
The condition (B) holds with $s=j_{{\rm id}_c}$. 

The converse is obvious. 

\hfill $\Box$ \\

\begin{Pro}
Let $\mathcal{C}$ be a category 
freely generated by a finite acyclic directed graph $\Gamma$ with a directed graph isomorphism 
$\sigma:\xymatrix{ \Gamma \ar[r]^-{\cong} & \Gamma^{\rm op}}$. 
For an object $c \in \mathcal{C}$, 
assume that there exist two different generating morphisms $f$ and $g$ 
such that ${\rm cod}(f)=c={\rm cod}(g)$ if $c \in |\mathcal{C}|$ is not minimal with respect to $<$. 
Then, (B') holds for $c$ and any $c' \in \mathcal{C}$. 
\label{pro8}
\end{Pro}
\textit{Proof. }
Take any family of maximal spans $\{(i_f,j_f)\}_{(c'',f) \in {\rm Elts}(\mathbf{y}(c))}$ such that 
${\rm cod}(i_f)=c''$ and ${\rm cod}(j_f)=c'$ for each $(c'',f) \in {\rm Elts}(\mathbf{y}(c))$ and 
assume that if $f=g \circ h$, then 
$(i_g,j_g)$ is a unique maximal span for the span $(h \circ i_f,j_f)$. 
If $c$ is a minimal element of $|\mathcal{C}|$, then we have no choice but to put $i_{{\rm id}_c}={\rm id}_c$. 
Otherwise, we have two different generating morphisms $f$ and $g$ such that ${\rm cod}(f)=c={\rm cod}(g)$ 
by the assumption of the claim. As in the proof of Lemma \ref{lem5}, we have $i_{{\rm id}_c} \circ k_f=f \circ i_f$. 
If $i_{{\rm id}_c} \neq {\rm id}_c$, then we can write $i_{{\rm id}_c}=a_n \circ \cdots \circ a_1$ for 
generating morphisms $a_1,\cdots,a_n$ and $n \geq 1$. Since $f$ is a generating morphism, it follows that 
$f=a_n$. By the same manner, we obtain $g=a_n$. Thus, we obtain $f=g$ but this is impossible because of the 
assumption of the claim. Consequently,  we have $i_{{\rm id}_c} = {\rm id}_c$. 

\hfill $\Box$ \\

\begin{Thm}
Let $\mathcal{C}$ be a category 
freely generated by a finite acyclic directed graph $\Gamma$ with a directed graph isomorphism 
$\sigma:\xymatrix{ \Gamma \ar[r]^-{\cong} & \Gamma^{\rm op}}$. 
Let $M_0 : \mathcal{C} \to \hat{\mathcal{C}}$ be the standard representation with respect to $\sigma$. 
Then, $(\hat{\mathcal{C}}({\rm Int}_{M_0}(-),G \otimes M_0),\varphi^{M_0,G})$ is an initial object of $\mathcal{IT}_G$ 
for any presheaf $G \in \hat{\mathcal{C}}$. 
\label{thm4}
\end{Thm}
\textit{Proof. }
For any $c,c' \in \mathcal{C}$, the map 
$\gamma_c(\sigma(c')):\mathbf{y}(\sigma(c))(\sigma(c')) \to {\rm Int}_{M_0}(c)(\sigma(c'))$ is injective 
by Lemma \ref{lem4} and Proposition \ref{pro6}. 
If $c$ is minimal or there exist two different generating morphisms whose codomains are $c$, then 
$\gamma_c(\sigma(c'))$ is surjective by Proposition \ref{pro7}, Lemma \ref{lem4}, Lemma \ref{lem5} and Proposition \ref{pro8}. 
Thus, if there exist two different generating morphisms whose codomains are $c$ for any non-minimal element 
$c \in |\mathcal{C}|$, then the natural transformation $\gamma : \mathbf{y} \circ \sigma \to {\rm Int}_{M_0}$ 
is a natural isomorphism. Consequently, the claim follows by Theorem \ref{thm3}. 

If there is an object $c$ such that the number of generating morphism $f$ such that ${\rm cod}(f)=c$ is 1, 
then $\gamma$ is not necessarily a natural isomorphism. Indeed, this happens in the simplest case 
$\Gamma=\{\bullet \to \bullet\}$. However, if $f:c' \to c$ is a unique generating morphism such that 
${\rm cod}(f)=c$, then we have $M_0(c') \cong M_0(c)$ and ${\rm Int}_{M_0}(c') \cong {\rm Int}_{M_0}(c)$. 
Consequently, $\iota_c$ is uniquely determined by $\iota_{c'}$, where $\iota_c$ and $\iota_{c'}$ are as in the proof of 
Theorem \ref{thm3}. Thus, the claim follows 
by slightly modifying the well-founded induction used in the proof of Theorem \ref{thm3}. 

\hfill $\Box$ \\

The representation $M_0$ defined in Section \ref{subsec2-1} is isomorphic to the standard representation 
for the isomorphism $\sigma:\xymatrix{ \upuparrows \ar[r]^-{\cong} & \upuparrows^{\rm op}}$ 
defined by $\sigma(0)=1,\ \sigma(1)=0,\ \sigma(m_0)=m_1,\ \sigma(m_1)=m_0$. 
In general, different dual structure on the same category gives rise to non-isomorphic standard representations. 
However, if $\mathcal{C}$ satisfies the assumption of Theorem \ref{thm3} or Theorem \ref{thm4}, 
then all standard representations induce the isomorphic objects in $\mathcal{IT}_G$ and 
they all are initial objects. For example, 
the dual structure $\sigma'$ on $\upuparrows$ given by $\sigma(m_0)=m_0,\ \sigma(m_1)=m_1$ induces 
a non-isomorphic standard representation to that defined in Section \ref{subsec2-1}. 

Let us see some examples to which Theorem \ref{thm4} are applied. 
We denote the category freely generated by the graph 
\[
\Gamma=
\xymatrix{
0 \ar[r]^-{m} & 1
}
\]
by $\uparrow$. 
The presheaf category $\hat{\uparrow}$ consists of all functions between two sets and 
pairs of functions between them satisfying an obvious commutative diagram. Let $\sigma$ be 
the unique dual structure on $\uparrow$. It is obvious that the standard representation $M_0$ 
with respect to $\sigma$ is given by maps from a singleton to a singleton for both $M_0(0)$ and 
$M_0(1)$. Hence, for any $G \in \hat{\uparrow}$, 
$\varphi_0^{M_0,G} : G(0) \to \hat{\uparrow}({\rm Int}_{M_0}(0),G \otimes M_0) \cong G(0)$ 
is a bijection and 
$\varphi_1^{M_0,G} : G(1) \to \hat{\uparrow}({\rm Int}_{M_0}(1),G \otimes M_0) \cong G(0)$ is 
isomorphic to the map $G(m):G(1) \to G(0)$. Theorem \ref{thm4} implies that 
the partition of the set $G(1)$ by the fibers of the map $G(m)$ is the finest one 
among those induced by all representations into bicomplete categories. 

Potentially important categories in applications to which Theorem \ref{thm4} are applied include 
the categories freely generated by the following graphs: 
\[
\xymatrix{
\bullet \ar[r] \ar[dr] & \bullet \\
\bullet \ar[r] \ar[ur] & \bullet 
}
\]
and 
\[
\xymatrix{
\bullet \ar@<1ex>[r] \ar@<-1ex>[r] & \bullet \ar@<1ex>[r] \ar@<-1ex>[r] & \cdots \ar@<1ex>[r] \ar@<-1ex>[r] & \bullet. 
}
\]
In the former case, its presheaf category is the category of bipartite directed networks. 
Bipartite networks are frequently used to represent metabolic networks. 
In the latter case, we obtain the category consisting of ``hierarchical directed networks''. 
This could be used to represent hierarchical systems. 
However, here, we do not pursue these examples further. 

\subsection{Stability with respect to representations and gluing condition for representations}
Now, we proceed to generalizations of Theorem \ref{thm2} and Proposition \ref{pro4}. 
Throughout this subsection, $\mathcal{C}$ is a well-founded small category equipped with 
a dual structure $\sigma$. 

For an object $c \in \mathcal{C}$, 
let ${\rm Elts}(\mathbf{y}(c))^*$ be a category whose objects are the same as 
those in ${\rm Elts}(\mathbf{y}(c))$ and morphisms are defined as follows: 
any morphism $h:(c',f) \to (c'',g)$ in ${\rm Elts}(\mathbf{y}(c))$ such that 
$(c'',g) \neq (c,{\rm id}_c)$ is a morphism in ${\rm Elts}(\mathbf{y}(c))^*$. 
We denote it by $h^*$ to indicate it is a morphism in ${\rm Elts}(\mathbf{y}(c))^*$. 
Note that if $(c'',g) \neq (c,{\rm id}_c)$, then $(c',f)$ cannot be $(c,{\rm id}_c)$ 
because of the well-foundedness of $\mathcal{C}$ for a morphism $h:(c',f) \to (c'',g)$ in ${\rm Elts}(\mathbf{y}(c))$. 
For any morphism $f:(c',f) \to (c,{\rm id}_c)$ in ${\rm Elts}(\mathbf{y}(c))$, 
we add a morphism $f^*$ which goes in the reverse way to ${\rm Elts}(\mathbf{y}(c))^*$: 
$f^*: (c,{\rm id}_c) \to (c',f)$. The composition of morphisms 
$f^*: (c,{\rm id}_c) \to (c',f)$ and $h^*:(c',f) \to (c'',g)$ in ${\rm Elts}(\mathbf{y}(c))^*$ is defined as 
$h^* \circ f^*:=g^*$. 

Let $\mathcal{D}$ be a complete category and $M:\mathcal{C} \to \mathcal{D}$ 
a representation. For each $c \in \mathcal{C}$, we define a functor $I_{M,c}:{\rm Elts}(\mathbf{y}(c))^* \to \mathcal{D}$ 
as follows: for objects in ${\rm Elts}(\mathbf{y}(c))^*$, 
$I_{M,c}(c,{\rm id}_c)={\rm Int}_M(c)$ and $I_{M,c}(c',f)=M(c')$ if $(c',f) \neq (c,{\rm id}_c)$. 
For morphisms, ${\rm id}_c:(c,{\rm id}_c) \to (c,{\rm id}_c)$ is sent to ${\rm id}_{{\rm Int}_M(c)}$, 
$f^*:(c,{\rm id}_c) \to (c',f)$ is sent to $\nu_{(c',f)}^{M,c}:{\rm Int}_M(c) \to M(c')$ if 
$(c',f) \neq (c,{\rm id}_c)$ and $h^*:(c',f) \to (c'',g)$ is sent to $M(h):M(c') \to M(c'')$ if 
$(c',f) \neq (c,{\rm id}_c)$ and $(c'',g) \neq (c,{\rm id}_c)$. 

One can check that $I_{M,c}$ is indeed a functor by using the facts that $M$ is a functor and 
$\nu^{M,c}$ is a cone on $M \circ \pi_c$. 

The following definition generalizes the gluing condition given in Definition \ref{def5}. 

\begin{Def}
Let $\mathcal{C}$ be a well-founded small category with a dual structure $\sigma$ and 
$\mathcal{D}$ a bicomplete category. We say that a representation $M:\mathcal{C} \to \mathcal{D}$ 
satisfies the \textbf{gluing condition} if 
the cocone 
\begin{align*}
& \{ \nu_{(c,{\rm id}_c)}^{M,c}:I_{M,c}(c,{\rm id}_c)={\rm Int}_M(c) \to M(c) \} \\
&\cup \{ M(f): I_{M,c}(c',f)=M(c') \to M(c) \}_{(c',f) \in {\rm Elts}(\mathbf{y}(c))^*, (c',f) \neq (c,{\rm id}_c)}
\end{align*}
is a colimit on $I_{M,c}$. 
\label{def11}
\end{Def}

\begin{Pro}
Let $\mathcal{C}$ be a well-founded small category with a dual structure $\sigma$. 
Then, the standard representation $M_0$ with respect to $\sigma$ satisfies the gluing condition. 
\label{pro9}
\end{Pro}
\textit{Proof. }
Let 
\begin{align*}
\{ \theta_{(c,{\rm id}_c)}:{\rm Int}_{M_0}(c) \to D \} 
\cup \{ \theta_{(c',f)}: M_0(c') \to D \}_{(c',f) \in {\rm Elts}(\mathbf{y}(c))^*, (c',f) \neq (c,{\rm id}_c)}
\end{align*}
be a cocone on $I_{M_0,c}$. Since colimits in $\hat{\mathcal{C}}$ can be computed pointwise, 
it is sufficient to show the universality of $M_0(c)(d)$ for each $d \in \mathcal{C}$. 
In what follows, we put $d=\sigma(\overline{d})$. 

Let us define a map $\rho:M_0(c)(d) \to D(d)$ in the following way: 
when $c$ is non-minimal element of $|\mathcal{C}|$ with respect to $<$, then 
for any element of $M_0(c)(d)$ we can take a representative $(\sigma(i),\sigma(j))$ of it such that 
$i \neq {\rm id}_c$. Then, we define $\rho([(\sigma(i),\sigma(j))])=\theta_{(e,i)}(d)([({\rm id}_{\sigma(e)},\sigma(j))])$, 
where $e={\rm dom}(i)={\rm dom}(j)$. 
\[
\xymatrix{
& M_0(e)(d) \ar[ld]_-{\theta_{(e,i)}(d)} \ar[d]^-{M_0(i)(d)} \\
D(d) & M_0(c)(d) \ar[l]^-{\rho} 
}
\]
If $c$ is minimal, then any element of $M_0(c)(d)$ takes a form 
$[({\rm id}_{\sigma(c)},\sigma(j))]=\{ ({\rm id}_{\sigma(c)},\sigma(j)) \}$ for some $j:c \to \overline{d}$ and 
${\rm Int}_{M_0}(c)=M_0(c)$. Thus, we define 
$\rho([({\rm id}_{\sigma(c)},\sigma(j))])=\theta_{(c,{\rm id}_c)}(d)([({\rm id}_{\sigma(c)},\sigma(j))])$. 
Note that in either case, we have the same form for $\rho$. 
It is straightforward to check that $\rho$ is a well-defined map. 

Now, we show that 
$\rho \circ M_0(f)(d)=\theta_{(c',f)}(d)$ for any $(c',f) \in {\rm Elts}(\mathbf{y}(c))^*$ such that $(c',f) \neq (c,{\rm id}_c)$ 
and $\rho \circ \nu_{(c,{\rm id}_c)}^{M_0,c}(d)=\theta_{(c,{\rm id}_c)}(d)$. 
For any morphism $f:c' \to c$ such that $c'\neq c$, we have 
\begin{align*}
\rho \circ M_0(f)(d)([(\sigma(i),\sigma(j))])
&= \rho([(\sigma(i) \circ \sigma(f),\sigma(j))]) \\
&= \theta_{(e,f \circ i)}(d)([({\rm id}_{\sigma(e)},\sigma(j))]) \\
&= \theta_{(c',f)}(d) \circ M_0(i)(d)([({\rm id}_{\sigma(e)},\sigma(j))]) \\
&= \theta_{(c',f)}(d) ([(\sigma(i),\sigma(j))]) 
\end{align*}
for any $[(\sigma(i),\sigma(j))] \in M_0(c')(d)$, where the third equality follows from the assumption 
that $\theta$ is a cocone on $I_{M_0,c}$. 
\[
\xymatrix@-1pc{
\sigma(c) \ar[r]^-{\sigma(f)} & \sigma(c') \ar[r]^-{\sigma(i)} & \sigma(e) \\
& \sigma(\overline{d}) \ar[ru]_-{\sigma(j)}
}
\]
For ${\rm id}_c:c \to c$, if $c$ is minimal, then $\nu_{(c,{\rm id}_c)}^{M_0,c}(d)$ is an identity. 
Hence, we have nothing to do. If $c$ is not minimal, then there exists a morphism $f:c' \to c$ such that 
$c' \neq c$. The desired equality follows because the following diagram commutes: 
\[
\xymatrix{
& {\rm Int}_{M_0}(c)(d) \ar[ldd]_-{\theta_{(c,{\rm id}_c)}(d)} \ar[d]^-{\nu_{(c',f)}^{M_0,c}(d)} \ar@/^2pc/[rdd]^-{\nu_{(c,{\rm id}_c)}^{M_0,c}(d)} & \\
& M_0(c')(d) \ar[ld]^-{\theta_{(c',f)}(d)} \ar[rd]^-{M_0(f)(d)} & \\
D(d) & & M_0(c)(d). \ar[ll]^-{\rho}
}
\]

Finally, we show the uniqueness of $\rho$. 
Suppose that there exists a map $\rho' : M_0(c)(d) \to D(d)$ such that 
$\rho' \circ M_0(f)(d)=\theta_{(c',f)}(d)$ for any $(c',f) \in {\rm Elts}(\mathbf{y}(c))^*$ such that $(c',f) \neq (c,{\rm id}_c)$ 
and $\rho' \circ \nu_{(c,{\rm id}_c)}^{M_0,c}(d)=\theta_{(c,{\rm id}_c)}(d)$. 
Then, for any $[(\sigma(i),\sigma(j))] \in M_0(c)(d)$, we have 
\begin{align*}
\rho'([(\sigma(i),\sigma(j))])
&= \rho' \circ M_0(i)(d)([({\rm id}_{\sigma(e)},\sigma(j))]) \\
&= \theta_{(e,i)}(d)([({\rm id}_{\sigma(e)},\sigma(j))]) \\
&= \rho([(\sigma(i),\sigma(j))]), 
\end{align*}
where we take $i$ as non-identity when $c$ is not minimal. 

\hfill $\Box$ \\

The next theorem is a generalization of Theorem \ref{thm2}. 

\begin{Thm}
Let $\mathcal{C}$ be a well-founded small category with a dual structure $\sigma$. 
For a presheaf $G \in \hat{\mathcal{C}}$, assume that the standard representation $M_0$ with respect 
to $\sigma$ induces an initial object in the category $\mathcal{IT}_G$. 
Let $M:\mathcal{C} \to \mathcal{D}$ be a representation into a bicomplete category $\mathcal{D}$ 
satisfying the gluing condition. 
If $G$ is stable with respect to $M$, then $G$ is also stable with respect to $M_0$.  
\label{thm5}
\end{Thm}
\textit{Proof. }
For a presheaf $G \in \hat{\mathcal{C}}$, assume that $\eta_c^{M,G}$ is an isomorphism 
for any $c \in \mathcal{C}$. We show that $\eta_c^{M_0,G}$ is also an isomorphism 
for any $c \in \mathcal{C}$ by appealing to a well-founded induction with respect to $<$ on 
$|\mathcal{C}|$. 

If $c$ is a minimal element of $|\mathcal{C}|$, then we have 
$M_0(c)={\rm Int}_{M_0}(c) \cong \mathbf{y}(\sigma(c))$. Hence, 
we can identify $\eta_c^{M_0,G}$ with $\varphi_c^{M_0,G}$ which has been shown to be an isomorphism 
in the proof of Theorem \ref{thm3}. 

Let $c$ be a non-minimal element of $|\mathcal{C}|$. 
Assume that $\eta_{c'}^{M_0,G}$ is an isomorphism for any $c' \in \mathcal{C}$ such that $c'<c$. 
We first show that $\eta_c^{M_0,G}$ is surjective. 
For any morphism $\delta:M_0(c) \to G \otimes M_0$, we shall construct a morphism 
$\overline{\delta}:M(c) \to G \otimes M$ and then show that 
$\eta_c^{M_0,G}\left( \left(\eta_c^{M,G}\right)^{-1}(\overline{\delta})\right)=\delta$. 

Let us construct a cocone 
\begin{align*}
\{ \theta_{(c,{\rm id}_c)}:{\rm Int}_M(c) \to G \otimes M \} 
\cup \{ \theta_{(c',f)}: M(c') \to G \otimes M \}_{(c',f) \in {\rm Elts}(\mathbf{y}(c))^*, (c',f) \neq (c,{\rm id}_c)}
\end{align*}
on $I_{M,c}$ by defining 
$\theta_{(c',f)}=\eta_{c'}^{M,G}\left( \left(\eta_{c'}^{M_0,G}\right)^{-1}(\delta \circ M_0(f))\right)$ 
for $(c',f) \neq (c,{\rm id}_c)$ and 
$\theta_{(c,{\rm id}_c)}=\iota_c \left( \delta \circ \nu_{(c,{\rm id}_c)}^{M_0,c} \right)$, 
where $\iota$ is a unique morphism from $\varphi^{M_0,G}$ to $\varphi^{M,G}$ in $\mathcal{IT}_G$. 
To show that $\theta$ is indeed a cocone on $I_{M,c}$, we have to show that 
$\theta_{(c',f)}=\theta_{(c'',g)} \circ M(h)$ for any morphism $h^*:(c',f) \to (c'',g)$ in ${\rm Elts}(\mathbf{y}(c))^*$ 
such that $(c'',g) \neq (c,{\rm id}_c)$ and 
$\theta_{(c,{\rm id}_c)}=\theta_{(c',f)} \circ \nu_{(c',f)}^{M,c}$ for any morphism $f^*:(c,{\rm id}_c) \to (c',f)$ 
in ${\rm Elts}(\mathbf{y}(c))^*$ such that $(c,{\rm id}_c) \neq (c',f)$. 

For the former equality, we have 
\begin{align*}
\theta_{(c'',g)} \circ M(h)
&= \eta_{c''}^{M,G}\left( \left(\eta_{c''}^{M_0,G} \right)^{-1}(\delta \circ M_0(g)) \right) \circ M(h) \\
&= \eta_{c'}^{M,G} \left( G(h) \left( \left( \eta_{c''}^{M_0,G} \right)^{-1}(\delta \circ M_0(g)) \right) \right) \\
&= \eta_{c'}^{M,G} \left( \left( \eta_{c''}^{M_0,G} \right)^{-1}(\delta \circ M_0(g) \circ M_0(h)) \right) \\
&= \eta_{c'}^{M,G} \left( \left( \eta_{c''}^{M_0,G} \right)^{-1}(\delta \circ M_0(f)) \right) \\
&= \theta_{(c',f)}, 
\end{align*}
where we use the naturality of $\eta^{M,G}$ in the second equality and that of $\eta^{M_0,G}$ in the 
third equality. For the latter equality, 
\begin{align*}
\theta_{(c,{\rm id}_c)}
&= \iota_c \left( \delta \circ \nu_{(c,{\rm id}_c)}^{M_0,c} \right) \\
&= \iota_c \left( \delta \circ M_0(f) \circ \nu_{(c',{\rm id}_{c'})}^{M_0,c'} \circ {\rm Int}_{M_0}(f) \right) \\
&= \iota_{c'} \left( \delta \circ M_0(f) \circ \nu_{(c',{\rm id}_{c'})}^{M_0,c'} \right) \circ {\rm Int}_M(f) \\
&= \iota_{c'} \left( \varphi_{c'}^{M_0,G}(x) \right) \circ {\rm Int}_M(f) \\
&= \varphi_{c'}^{M,G}(x) \circ {\rm Int}_M(f) \\
&= \mu_{(c',x)}^{M,G} \circ \nu_{(c',{\rm id}_{c'})}^{M,c'} \circ {\rm Int}_M(f) \\
&= \mu_{(c',x)}^{M,G} \circ \nu_{(c',f)}^{M,c'} \\
&= \theta_{(c',f)} \circ \nu_{(c',f)}^{M,c'}, 
\end{align*}
where the second equality follows because the diagram 
\begin{equation}
\begin{split}
\xymatrix{
{\rm Int}_{M_0}(c) \ar[r]^-{{\rm Int}_{M_0}(f)} \ar[d]_-{\nu_{(c,{\rm id}_c)}^{M_0,c}} \ar[rd]^-{\nu_{(c',f)}^{M_0,c}} & {\rm Int}_{M_0}(c') \ar[d]^-{\nu_{(c',{\rm id}_{c'})}^{M_0,c'}} \\
M_0(c) & M_0(c') \ar[l]_-{M_0(f)} 
}
\end{split}
\label{diag1}
\end{equation}
commutes and the third equality is due to the naturality of $\iota$. 
In the fourth equality, we put $x:=\left( \eta_{c'}^{M_0,G} \right)^{-1} (\delta \circ M_0(f)) \in G(c')$. 
Thus, we have 
$\varphi_{c'}^{M_0,G}(x)=\mu_{(c',x)}^{M,G} \circ \nu_{(c',{\rm id}_{c'})}^{M,c'}=\eta_{c'}^{M_0,G}(x) \circ \nu_{(c',{\rm id}_{c'})}^{M,c'}$. 
For the remaining equalities, we use the naturality of $\iota$ in the fifth equality, 
the definition of $\varphi_{c'}^{M,G}(x)$ in the sixth equality, the definition of ${\rm Int}_M(f)$ in 
the seventh equality and $\theta_{(c',f)}=\eta_{c'}^{M,G}(x)=\mu_{(c',x)}^{M,G}$ in the last equality. 

Consequently, we obtain a unique morphism $\overline{\delta}: M(c) \to G \otimes M$ such that 
$\theta_{(c,{\rm id}_c)}=\overline{\delta} \circ \nu_{(c,{\rm id}_c)}^{M,c}$ and 
$\theta_{(c',f)}=\overline{\delta} \circ M(f)$ for any morphism $f:c' \to c$ such that $c' \neq c$. 

Now, we show that $\eta_c^{M_0,G} \left( \left( \eta_c^{M,G} \right)^{-1} (\overline{\delta}) \right)=\delta$. 
We put the left-hand side by $\delta'$ in what follows. We appeal to the fact that $M_0$ satisfies the 
gluing condition (Proposition \ref{pro9}). If we can show that 
$\delta' \circ \nu_{(c,{\rm id}_c)}^{M_0,c}=\delta \circ \nu_{(c,{\rm id}_c)}^{M_0,c}=:\zeta_{(c,{\rm id}_c)}$ and 
$\delta' \circ M_0(f)=\delta \circ M_0(f)=:\zeta_{(c',f)}$ for any morphism $f:c' \to c$ such that $c' \neq c$, 
then $\zeta$ is a cocone on $I_{M_0,c}$. Since $M_0$ satisfies the gluing condition by Proposition \ref{pro9}, there is 
a unique morphism $\kappa:M_0(c) \to G \otimes M_0$ such that $\kappa \circ \nu_{(c,{\rm id}_c)}^{M_0,c}=\zeta_{(c,{\rm id}_c)}$ 
and $\kappa \circ M_0(f)=\zeta_{(c',f)}$ for any morphism $f:c' \to c$ such that $c' \neq c$. However, both 
$\delta$ and $\delta'$ satisfy the condition for $\kappa$, we conclude $\delta=\kappa=\delta'$. 

Let $f:c' \to c$ be a morphism such that $c' \neq c$. We want to show that $\delta' \circ M_0(f)=\delta \circ M_0(f)$ holds. 
Putting $x=\left( \eta_c^{M,G} \right)^{-1} (\overline{\delta}) \in G(c)$, we have 
$\delta' \circ M_0(f)=\eta_c^{M_0,G}(x) \circ M_0(f)=\eta_{c'}^{M_0,G}(x \cdot f)$, where 
we use the naturality of $\eta^{M_0,G}$ in the second equality. 
Now, we have 
\begin{align*}
x \cdot f
&= G(f) \left( \left( \eta_c^{M,G} \right)^{-1} (\overline{\delta}) \right) \\
&= \left( \eta_{c'}^{M,G} \right)^{-1} (\overline{\delta} \circ M(f)) \\
&= \left( \eta_{c'}^{M,G} \right)^{-1} (\theta_{(c',f)}) \\
&= \left( \eta_{c'}^{M,G} \right)^{-1} \left( \eta_{c'}^{M,G}\left( \left(\eta_{c'}^{M_0,G}\right)^{-1}(\delta \circ M_0(f))\right) \right) \\
&= \left( \eta_{c'}^{M_0,G} \right)^{-1} (\delta \circ M_0(f)), 
\end{align*}
where we use the naturality of $\eta^{M,G}$ in the second equality. 
Thus, we obtain the desired equality. 

To show the equality $\delta' \circ \nu_{(c,{\rm id}_c)}^{M_0,c}=\delta \circ \nu_{(c,{\rm id}_c)}^{M_0,c}$, 
take a morphism $f:c' \to c$ such that $c' \neq c$ whose existence is guaranteed by the assumption that 
$c$ is not minimal. We have 
\begin{align*}
\delta' \circ \nu_{(c,{\rm id}_c)}^{M_0,c}
&= \eta_c^{M_0,G}(x) \circ \nu_{(c,{\rm id}_c)}^{M_0,c} \\
&= \mu_{(c,x)}^{M_0,G} \circ \nu_{(c,{\rm id}_c)}^{M_0,c} \\
&= \varphi_c^{M_0,G}(x) \\
&= \varphi_{c'}^{M_0,G}(x \cdot f) \circ {\rm Int}_{M_0}(f) \\
&= \mu_{(c',x \cdot f)}^{M_0,G} \circ \nu_{(c',{\rm id}_c')}^{M_0,c'} \circ {\rm Int}_{M_0}(f) \\
&= \eta_{c'}^{M_0,G}(x \cdot f) \circ \nu_{(c',{\rm id}_c')}^{M_0,c'} \circ {\rm Int}_{M_0}(f) \\
&= \eta_{c'}^{M_0,G}\left( \left( \eta_{c'}^{M,G} \right)^{-1} (\delta \circ M_0(f)) \right) \circ \nu_{(c',{\rm id}_c')}^{M_0,c'} \circ {\rm Int}_{M_0}(f) \\
&= \delta \circ M_0(f) \circ \nu_{(c',{\rm id}_c')}^{M_0,c'} \circ {\rm Int}_{M_0}(f) \\
&= \delta \circ \nu_{(c,{\rm id}_c)}^{M_0,c}, 
\end{align*}
where the fourth equality follows from the naturality of $\varphi^{M_0,G}$ and the last equality is due to 
the commutative diagram (\ref{diag1}). This complete the proof that $\eta_c^{M_0,G}$ is a surjection. 

Next, we show that $\eta_c^{M_0,G}$ is injective. 
For any $x,y \in G(c)$, suppose that $\eta_c^{M_0,G}(x)=\eta_c^{M_0,G}(y)$ holds. 
We put the both sides of the equality by $\delta$. If we construct a unique morphism $\overline{\delta}:M(c) \to G \otimes M$ 
such that $\overline{\delta} \circ \nu_{(c,{\rm id}_c)}^{M,c}=\theta_{(c,{\rm id}_c)}$ and 
$\overline{\delta} \circ M(f)=\theta_{(c',f)}$ for any morphism $f:c' \to c$ such that $c' \neq c$ 
by the same manner as in the proof of the surjectivity. 

We shall show that both $\eta_c^{M,G}(x)$ and $\eta_c^{M,G}(y)$ satisfy the condition for $\overline{\delta}$. 
Then, we will obtain $\eta_c^{M,G}(x)=\overline{\delta}=\eta_c^{M,G}(y)$ which implies $x=y$ because $\eta_c^{M,G}$ 
is a bijection. 

For any morphism $f:c' \to c$ such that $c' \neq c$, we have 
\begin{align*}
\eta_c^{M,G}(x) \circ M(f)
&= \eta_{c'}^{M,G}(x \cdot f) \\
&= \eta_{c'}^{M,G} \left( \left( \eta_{c'}^{M_0,G} \right)^{-1} (\delta \circ M_0(f)) \right) \\
&= \theta_{(c',f)}, 
\end{align*}
where we use the naturality of $\eta^{M_0,G}$ in the second equality. By the same manner, we obtain 
$\eta_c^{M,G}(y) \circ M(f)=\overline{\delta} \circ M(f)$. 
We also have 
\begin{align*}
\eta_c^{M,G}(x) \circ \nu_{(c,{\rm id}_c)}^{M,c} 
&= \mu_{(c,x)}^{M,G} \circ \nu_{(c,{\rm id}_c)}^{M,c} \\
&= \varphi_c^{M,G}(x) \\
&= \iota_c \circ \varphi_c^{M_0,G}(x) \\
&= \iota_c \left( \eta_c^{M_0,G}(x) \circ \nu_{(c,{\rm id}_c)}^{M_0,c} \right) \\
&= \iota_c \left( \delta \circ \nu_{(c,{\rm id}_c)}^{M_0,c} \right) \\
&= \theta_{(c,{\rm id}_c)}
\end{align*}
and $\eta_c^{M,G}(y) \circ \nu_{(c,{\rm id}_c)}^{M,c}=\theta_{(c,{\rm id}_c)}$ by the same manner, as desired. 
This completes the proof of the injectivity. 

\hfill $\Box$ \\

The next proposition is a generalization of Proposition \ref{pro4}. 

\begin{Pro}
Let $\mathcal{C}$ be a well-founded small category with a dual structure $\sigma$ and 
$M:\mathcal{C} \to \hat{\mathcal{C}}$ a representation satisfying the gluing condition. 
Then, the tensor product representation $M \otimes N$ also satisfies the gluing condition 
for any representation $N:\mathcal{C} \to \mathcal{D}$ into a bicomplete category $\mathcal{D}$. 
\label{pro10}
\end{Pro}
\textit{Proof. }
We show that the cocone 
\begin{align*}
&\{ \nu_{(c,{\rm id}_c)}^{M \otimes N,c}:{\rm Int}_{M \otimes N}(c) \to M \otimes N(c) \} \\
&\cup \{ M \otimes N(f): M \otimes N(c') \to M \otimes N(c) \}_{(c',f) \in {\rm Elts}(\mathbf{y}(c))^*, (c',f) \neq (c,{\rm id}_c)}
\end{align*}
is a colimit on $I_{M \otimes N,c}$ for any $c \in \mathcal{C}$. 

Let 
\begin{align*}
\{ \theta_{(c,{\rm id}_c)}:{\rm Int}_{M \otimes N}(c) \to X \} 
\cup \{ \theta_{(c',f)}: M \otimes N(c') \to X \}_{(c',f) \in {\rm Elts}(\mathbf{y}(c))^*, (c',f) \neq (c,{\rm id}_c)}
\end{align*}
be a cocone on $I_{M \otimes N,c}$. Since $M$ satisfies the gluing condition, 
\begin{align*}
\{ \nu_{(c,{\rm id}_c)}^{M,c}:{\rm Int}_M(c) \to M(c) \} 
\cup \{ M(f): M(c') \to M(c) \}_{(c',f) \in {\rm Elts}(\mathbf{y}(c))^*, (c',f) \neq (c,{\rm id}_c)}
\end{align*}
is a colimit on $I_{M,c}$. Since $(-) \otimes N$ preserves colimits, 
\begin{align*}
&\{ \nu_{(c,{\rm id}_c)}^{M,c} \otimes N:{\rm Int}_M(c) \otimes N \to M(c) \otimes N \} \\
&\cup \{ M(f) \otimes N : M(c') \otimes N \to M(c) \otimes N \}_{(c',f) \in {\rm Elts}(\mathbf{y}(c))^*, (c',f) \neq (c,{\rm id}_c)}
\end{align*}
is also a colimit on $(-) \otimes N \circ I_{M,c}$. Note that $M(c) \otimes N=M \otimes N(c)$ and $M(f) \otimes N=M \otimes N(f)$. 

Now, for the canonical map $\beta_c:{\rm Int}_M(c) \otimes N \to {\rm Int}_{M \otimes N}(c)$, 
\begin{align*}
&\{ \theta_{(c,{\rm id}_c)} \circ \beta_c :{\rm Int}_M(c) \otimes N \to X \} \\
&\cup \{ \theta_{(c',f)} : M \otimes N(c') \to X \}_{(c',f) \in {\rm Elts}(\mathbf{y}(c))^*, (c',f) \neq (c,{\rm id}_c)}
\end{align*}
is a cocone on $(-) \otimes N \circ I_{M,c}$ because the following diagram commutes for any morphism 
$f^*:(c,{\rm id}_c) \to (c',f)$ in ${\rm Elts}(\mathbf{y}(c))^*$ such that $(c',f) \neq (c,{\rm id}_c)$: 
\[
\xymatrix{
{\rm Int}_{M \otimes N}(c) \ar[r]^-{\theta_{(c,{\rm id}_c)}} \ar[rd]^-{\nu_{(c',f)}^{M \otimes N,c}} & X \\
{\rm Int}_M(c) \otimes N \ar[u]^-{\beta_c} \ar[r]_-{\nu_{(c',f)}^{M,c} \otimes N} & M \otimes N(c'). \ar[u]_-{\theta_{(c',f)}}
}
\]
Consequently, there exists a unique map $\xi_c:M \otimes N \to X$ such that 
(i) $\theta_{(c,{\rm id}_c)} \circ \beta_c=\xi_c \circ \nu_{(c,{\rm id}_c)}^{M,c} \otimes N$ and 
(ii) $\theta_{(c',f)}=\xi_c \circ M(f) \otimes N$ for any morphism $f:c' \to c$ such that $c' \neq c$. 

We shall show that the equality (iii) $\theta_{(c,{\rm id}_c)}=\xi_c \circ \nu_{(c,{\rm id}_c)}^{M \otimes N,c}$ holds. 
Indeed, if $c$ is minimal with respect to $<$ on $|\mathcal{C}|$, then $\beta_c$ is a bijection. 
Thus, we have 
\begin{align*}
\xi_c \circ \nu_{(c,{\rm id}_c)}^{M \otimes N,c}
&= \xi \circ \nu_{(c,{\rm id}_c)}^{M,c} \otimes N \circ \beta_c^{-1} \\
&= \theta_{(c,{\rm id}_c)} \circ \beta_c \circ \beta_c^{-1}=\theta_{(c,{\rm id}_c)}, 
\end{align*}
where we use (i) in the second equality. If $c$ is not minimal, then there exists a morphism $f:c' \to c$ such that $c' \neq c$. 
In this case, (iii) holds because the following diagram commutes: 
\[
\xymatrix@-1pc{
M \otimes N(c) \ar[rr]^-{\xi_c} & & X \\
& M \otimes N (c') \ar[lu]_-{M \otimes N(f)} \ar[ru]^-{\theta_{(c',f)}} & \\
{\rm Int}_{M \otimes N}(c) \ar[uu]^-{\nu_{(c,{\rm id}_c)}^{M \otimes N,c}} \ar[ru]^-{\nu_{(c',f)}^{M \otimes N,c}} \ar@/_2pc/[rruu]_-{\theta_{(c,{\rm id}_c)}}
}
\]
Moreover, $\xi_c$ satisfying both (ii) and (iii) is unique. Suppose a map $\xi'_c:M \otimes N(c) \to X$ satisfies (ii) and (iii). 
By composing (iii) and $\beta_c$, we obtain (i) with $\xi_c$ replaced by $\xi'_c$. Thus, (i) and (ii) also holds for $\xi'_c$. 
Since $\xi_c$ satisfying both (i) and (ii) is unique, we conclude that $\xi_c=\xi'_c$. 
This completes the proof of the claim. 

\hfill $\Box$ \\

Finally, we study the relationship between the stability with respect to a representation 
$M$ satisfying the gluing condition and the stability with respect to a tensor product 
representation $M \otimes N$. 

\begin{Pro}
Let $\mathcal{C}$ be a well-founded small category with a dual structure $\sigma$, 
$M:\mathcal{C} \to \hat{\mathcal{C}}$ a representation satisfying the gluing condition 
and $N:\mathcal{C} \to \mathcal{D}$ a representation into a bicomplete category $\mathcal{D}$. 
For a presheaf $G \in \hat{\mathcal{C}}$, assume that 
the map $(-)\otimes N: \hat{\mathcal{C}}(M(c),G \otimes M) \to \mathcal{D}(M \otimes N(c),(G \otimes M)\otimes N)$ 
is injective for any minimal element $c \in |\mathcal{C}|$ with respect to $<$. 
Then, if $G$ is stable with respect to $M \otimes N$, then $G$ is also stable with respect to $M$. 
\label{pro11}
\end{Pro}
\textit{Proof. }
Let $G$ be stable with respect to $M \otimes N$ and assume that 
the map $(-)\otimes N: \hat{\mathcal{C}}(M(c),G \otimes M) \to \mathcal{D}(M \otimes N(c),(G \otimes M)\otimes N)$ 
is injective for any minimal element $c \in |\mathcal{C}|$ with respect to $<$. 
Define the natural transformation $j:\hat{\mathcal{C}}(M(-),G \otimes M) \to \mathcal{D}(M \otimes N(-),G \otimes (M \otimes N))$ 
by $j_c(\delta):=\alpha_G \circ \delta \otimes N$ for any $c \in \mathcal{C}$ and morphism $\delta:M(c) \to G \otimes M$, 
where $\alpha_G$ is the canonical isomorphism $\alpha_G: \xymatrix{ (G \otimes M) \otimes N \ar[r]^-{\cong} & G \otimes (M \otimes N) }$. 
Then, we have $\eta_c^{M \otimes N,G}=j_c \circ \eta_c^{M,G}$ because 
$\eta_c^{M \otimes N,G}(x)=\mu_{(c,x)}^{M \otimes N,G}=\alpha_G \circ \mu_{(c,x)}^{M,G} \otimes N=j_c \left( \mu_{(c,x)}^{M,G} \right)=j_c \left( \eta_c^{M,G}(x) \right)$ 
for any $x \in G(c)$. 
From this, we can show that $\eta_c^{M,G}$ is an injection for any $c \in \mathcal{C}$. 
Indeed, if $\eta_c^{M,G}(x)=\eta_c^{M,G}(y)$ for $x,y \in G(c)$, then 
$\eta_c^{M \otimes N,G}(x)=j_c \left( \eta_c^{M,G}(x) \right)=j_c \left( \eta_c^{M,G}(y) \right)=\eta_c^{M \otimes N,G}(y)$. 
Since $\eta_c^{M \otimes N,G}$ is assumed to be a bijection, we obtain $x=y$. 

To show that $\eta_c^{M,G}$ is a surjection for any $c \in \mathcal{C}$, we appeal to a well-founded induction on 
$|\mathcal{C}|$ with respect to $<$. 
If $c$ is a minimal element of $|\mathcal{C}|$ with respect to $<$, then $\eta_c^{M,G}$ turns out to be a surjection. 
Indeed, for any morphism $\delta:M(c) \to G \otimes M$, if we put $x:=\left( \eta_c^{M \otimes N,G} \right)^{-1}(j_c(\delta)) \in G(c)$, 
then $j_c(\delta)=\eta_c^{M \otimes N,G}(x)=j_c \left( \eta_c^{M,G}(x) \right)$. By the assumption of the claim, 
$j_c$ is an injection. Hence, we obtain $\delta=\eta_c^{M,G}(x)$. 

When $c$ is not minimal, assume that $\eta_{c'}^{M,G}$ is surjective for any $c'<c$. 
Since we have shown that $\eta_{d}^{M,G}$ is injective for any $d \in \mathcal{C}$, 
$\eta_{c'}^{M,G}$ is a bijection for any $c'<c$. 
For any morphism $\delta:M(c) \to G \otimes M$, put $x:=\left( \eta_c^{M \otimes N,G} \right)^{-1}(j_c(\delta)) \in G(c)$ 
and $\delta'=\eta_c^{M,G}(x)$. 
Suppose that we have $\delta \circ \nu_{(c,{\rm id}_c)}^{M,c}=\delta' \circ \nu_{(c,{\rm id}_c)}^{M,c}=:\zeta_{(c,{\rm id}_c)}$ 
and $\delta \circ M(f)=\delta' \circ M(f)=:\zeta_{(c',f)}$ for any morphism $f:c' \to c$ such that $c' \neq c$. 
Then, $\zeta$ is a cocone on $I_{M,c}$. Since $M$ satisfies the gluing condition, there exists a unique morphism 
$\kappa:M(c) \to G \otimes M$ such that $\zeta_{(c,{\rm id}_c)}=\kappa \circ \nu_{(c,{\rm id}_c)}^{M,c}$ and 
$\zeta_{(c',f)}=\kappa \circ M(f)$ for any morphism $f:c' \to c$ such that $c' \neq c$. Since both 
$\delta$ and $\delta'$ satisfies the condition for $\kappa$, we conclude that $\delta=\kappa=\delta'$ which 
implies that $\eta_c^{M,G}$ is a surjection. 

Thus, our task is to show that 
$\delta \circ \nu_{(c,{\rm id}_c)}^{M,c}=\delta' \circ \nu_{(c,{\rm id}_c)}^{M,c}$ 
and $\delta \circ M(f)=\delta' \circ M(f)$ for any morphism $f:c' \to c$ such that $c' \neq c$. 
For the latter equality, fix any morphism $f:c' \to c$ such that $c' \neq c$. Then, we have 
\begin{align*}
\delta' \circ M(f) 
&= \eta_c^{M,G}(x) \circ M(f) \\
&= \eta_{c'}^{M,G}(x \cdot f) \\
&= \eta_{c'}^{M,G} \left( G(f) \left( \left( \eta_c^{M \otimes N,G} \right)^{-1} (j_c(\delta)) \right) \right) \\
&= \eta_{c'}^{M,G} \left( \left( \eta_{c'}^{M \otimes N,G} \right)^{-1} (j_c(\delta) \circ M \otimes N(f))  \right) \\
&= \eta_{c'}^{M,G} \left( \left( \eta_{c'}^{M \otimes N,G} \right)^{-1} (j_{c'}(\delta \circ M(f))) \right) \\
&= \delta \circ M(f), 
\end{align*}
where we use the naturality of $\eta^{M,G}$ in the second equality, that of $\eta^{M \otimes N,G}$ 
in the fourth equality and that of $j$ in the fifth equality. The last equality follows from the equality 
$\eta_{c'}^{M,G}=j_{c'}^{-1} \circ \eta_{c'}^{M \otimes N,G}$ which in turn obtained from the invertibility of 
$\eta_{c'}^{M,G}$ and $j_{c'}$. 

The former equality also follows from a similar calculation showing 
$\delta \circ \nu_{(c,{\rm id}_c)}^{M_0,c}=\delta' \circ \nu_{(c,{\rm id}_c)}^{M_0,c}$ in the proof of Theorem \ref{thm5}. 

\hfill $\Box$ \\

\begin{Cor}
Let $\mathcal{C}$ be a well-founded small category with a dual structure $\sigma$, 
$M:\mathcal{C} \to \hat{\mathcal{C}}$ a representation satisfying the gluing condition 
and $N:\mathcal{C} \to \mathcal{D}$ a representation into a bicomplete category $\mathcal{D}$. 
For a presheaf $G \in \hat{\mathcal{C}}$, assume that 
the component of the unit 
$\eta_c^{M,G}:G(c) \to \hat{\mathcal{C}}(M(c),G \otimes M)$ 
is an isomorphism for any minimal element $c \in |\mathcal{C}|$ with respect to $<$. 
Then, if $G$ is stable with respect to $M \otimes N$, then $G$ is also stable with respect to $M$. 
\label{cor2}
\end{Cor}

To obtain instances of Corollary \ref{cor2}, let us consider the case when $\mathcal{C}=\upuparrows$. 
If a representation $M:\upuparrows \to \hat{\upuparrows}$ satisfies $M(0)=\mathbf{y}(1)$ and the gluing condition, 
then we have $\eta_0^{M,G}:\xymatrix{ G(0) \ar[r]^-{\cong} & \hat{\mathcal{C}}(M(0),G \otimes M) }$ for any directed network $G$. 
This can be extended to a general case and proved as follows. 

\begin{Pro}
Let $\mathcal{C}$ be a well-founded small category with a dual structure $\sigma$. We assume that 
there exists the minimum element $c$ in $|\mathcal{C}|$ with respect to $<$. Then, if a representation 
$M:\mathcal{C} \to \hat{\mathcal{C}}$ satisfies the gluing condition and $M(c) \cong \mathbf{y}(\sigma(c))$, 
then we have $\eta_c^{M,G}:\xymatrix{ G(c) \ar[r]^-{\cong} & \hat{\mathcal{C}}(M(c),G \otimes M) }$ for any presheaf 
$G \in \hat{\mathcal{C}}$. 
\label{pro12}
\end{Pro}
\textit{Proof. }
For any $G \in \hat{\mathcal{C}}$, we have 
\begin{align*}
\hat{\mathcal{C}}(M(c),G \otimes M) 
&\cong \hat{\mathcal{C}}(\mathbf{y}(\sigma(c)),G \otimes M) \\
&\cong G \otimes M(\sigma(c)) \\
&= \left( \sum_{c' \in \mathcal{C}} G(c') \times M(c')(\sigma(c)) \right) / \sim.
\end{align*}
Now, suppose that $\{M(f):M(c'') \to M(c')\}_{(c'',f) \in {\rm Elts}(\mathbf{y}(c')), c'' \neq c'}$ is an 
epimorphic family for any $c' \neq c$. Then, each component of the morphism $\sum M(f) : \sum M(c'') \to M(c')$ is a surjection, 
where $\sum$ is taken over all $(c'',f) \in {\rm Elts}(\mathbf{y}(c'))$ such that $c'' \neq c'$. 
In particular, for any $a \in M(c')(\sigma(c))$ there exists $f:c'' \to c'$ and $b \in M(c'')(\sigma(c))$ such that 
$c'' \neq c'$ and $f \cdot b=a$. By a well-founded induction, we can show that for any $c' \neq c$ and $a \in M(c')(\sigma(c))$, 
there exists a morphism $f:c \to c'$ such that $f \cdot {\rm id}_{\sigma(c)}=a$ 
because $M(c)(\sigma(c)) \cong \mathbf{y}(\sigma(c))(\sigma(c)) = \{ {\rm id}_{\sigma(c)} \}$. 
Consequently, we obtain $G(c) \cong \hat{\mathcal{C}}(M(c) ,G \otimes M)$. One can see that this isomorphism is given by $\eta_c^{M,G}$. 

Let us show that $\{M(f):M(c'') \to M(c')\}_{(c'',f) \in {\rm Elts}(\mathbf{y}(c')), c'' \neq c'}$ is an 
epimorphic family for any $c' \neq c$. 
In general, if a functor $F:\mathcal{I} \to \mathcal{A}$ has a colimit, then it is an epimorphic family. 
Thus, the cocone 
\begin{align*}
& \{ \nu_{(c',{\rm id}_{c'})}^{M,c'}:I_{M,c'}(c',{\rm id}_{c'})={\rm Int}_M(c') \to M(c') \} \\
&\cup \{ M(f): I_{M,c'}(c'',f)=M(c'') \to M(c') \}_{(c'',f) \in {\rm Elts}(\mathbf{y}(c'))^*, (c'',f) \neq (c',{\rm id}_{c'})}
\end{align*}
is an epimorphic family, since $M$ satisfies the gluing condition. 
This is equivalent to saying that the morphism 
$\nu_{(c',{\rm id}_{c'})}^{M,c'}+\sum M(f):{\rm Int}_M(c')+\sum M(c'') \to M(c')$ is an epimorphism, 
where $\sum$ is taken over all $(c'',f) \in {\rm Elts}(\mathbf{y}(c'))^*$ such that $(c'',f) \neq (c',{\rm id}_{c'})$. 
Since $c'$ is not minimal, there exists a morphism $f:c'' \to c'$ such that $c'' \neq c'$ and 
$\nu_{c',{\rm id}_{c'}}^{M,c'}=M(f) \circ \nu_{(c'',f)}^{M,c'}$. Thus, 
$\sum M(f):\sum M(c'') \to M(c')$ is an epimorphism. 

\hfill $\Box$ \\

\section*{Acknowledgements}
This work was supported by JST PRESTO program. 
The author thanks Professor Tadao Oda, Professor Izumi Ojima, 
Professor Yukio Gunji, Professor Toru Tsujishita, 
Dr. Igor Balaz, Dr. Ichiro Hasuo, Dr. Asaki Nishikawa and Dr. Ken Shiotani 
for their invaluable comments and discussions on the work.

\end{document}